\documentclass[12pt]{article}
\usepackage{amsmath}
\usepackage{graphicx}
\usepackage{enumerate}
\usepackage{natbib}
\usepackage{subfigure}
\usepackage{lipsum}
\usepackage{amsthm}
\usepackage{url} 
\usepackage{xr}
\usepackage{color}
\newcommand{\blind}{1}

\newtheorem{theorem}{Theorem}
\newtheorem{remark}{Remark}
\newtheorem{definition}{Definition}

\newtheorem{lemma}{Lemma}
\newtheorem{example}{Example}
\newtheorem{corollary}{Corollary}

\newtheorem{algorithm}{Algorithm}

\begin{document}

\def\spacingset#1{\renewcommand{\baselinestretch}%
{#1}\small\normalsize} \spacingset{1}


\if1\blind
{
  \title{\bf Simultaneous Inference for Non-Stationary Random Fields, with Application to
Gridded Data Analysis}
  \author{Yunyi Zhang\thanks{
    The authors gratefully acknowledge \textit{National Natural Science Foundation of China, 2023SA0044.}}\hspace{.2cm}\\
    School of Data Science, The Chinese University of Hong Kong, Shenzhen.\\
    and \\
    Zhou Zhou\thanks{
    Corresponding author
    } \\
    Department of Statistical Sciences, University of Toronto}
  \maketitle
} \fi

\if0\blind
{
  \bigskip
  \bigskip
  \bigskip
  \begin{center}
    {\LARGE\bf Kernel smoothing on two-dimensional heterogeneous \& correlated random field: with application to
gridded data analysis}
\end{center}
  \medskip
} \fi

\bigskip
\begin{abstract}
Current statistics literature on statistical inference of random fields typically assumes that the fields are stationary or focuses on models of non-stationary 
Gaussian fields with parametric/semiparametric covariance families,
which may not be sufficiently flexible to tackle complex modern-era random field data. This paper performs simultaneous nonparametric statistical inference for a general class of non-stationary and non-Gaussian random fields by modeling the fields as nonlinear systems with location-dependent transformations of an underlying `shift random field'.
Asymptotic results, including concentration inequalities and
 Gaussian approximation theorems for high dimensional sparse linear forms of the random field, are derived.
A computationally efficient locally weighted multiplier bootstrap algorithm is proposed and theoretically verified as a unified tool for the simultaneous inference of the aforementioned non-stationary non-Gaussian random field. Simulations and real-life data examples demonstrate good performances and broad applications of the proposed algorithm.
\end{abstract}

\noindent%
\vfill

\newpage
\spacingset{1.9} 
\section{Introduction}
\label{section.intro}
An integer-indexed two-dimensional random field consists of a collection of dependent random variables $\{X_i^{(j)}, (i, j)\in\mathbf{Z}^2\}$ and is widely used in statistical applications.
Though much progress has been made in the statistical analyses of such data recently (\cite{MR4430002}, \cite{doi:10.1080/01621459.2023.2218578}), one major limitation of the state-of-the-art literature on the statistical inference of random fields is the assumption
of (close to) strict stationarity, i.e., the joint distribution of
$(X_{i_1}^{(j_1)}, X_{i_2}^{(j_2)},\cdots, X_{i_k}^{(j_k)})$ is the same as the joint distribution of $(X_{i_1 + c}^{(j_1 + d)}, X_{i_2 + c}^{(j_2 + d)},\cdots, X_{i_k + c}^{(j_k + d)})$ for any positive integer $k$ and any $i_1,\cdots, i_k, j_1, \cdots, j_k, c, d\in\mathbf{Z}$. As a result, $Var(X_i^{(j)})$ is identical for all $i,j$ and the covariances have to satisfy the constraint $$
Cov(X_{i_1}^{(j_1)}, X_{i_2}^{(j_2)}) =
Cov(X_{i_1 + c}^{(j_1 + d)}, X_{i_2 + c}^{(j_2 + d)}).
$$ However, as demonstrated in \cite{MR2655659}, \cite{https://doi.org/10.48550/arxiv.2004.06628} and Example \ref{example.hetero} below,
real-life data, such as those in spatial statistics,  often exhibit complicated non-stationary behavior in the mean, covariance and other distributional characteristics.
\begin{example}
Figure \ref{fig.real_life_data} plots the global long--term mean temperature anomaly of February from 1971 to 2000. The dataset is collected from \cite{data_set} \footnote{also see
\noindent \textit{https://psl.noaa.gov/data/gridded/data.noaaglobaltemp.html}}.  In addition, it demonstrates the estimated element--wise variances as well as element--wise row--wise covariances $Cov(X_i^{(j)}, X_i^{(j - 1)})$. According to  Figure \ref{fig.real_life_data}, the data exhibit strong spatial dependence. Besides, the variances and covariances vary significantly among different positions, thus exhibiting spatial non--stationarity. In Figure \ref{fig.real_life_data}, the variances and covariances of data follow complex patterns, which seem to be hard to model using simple parametric models. 
\begin{figure}[htbp]
\centering
    \subfigure[Mean Global Surface Temperature]{
    \includegraphics[width = 1.8in]{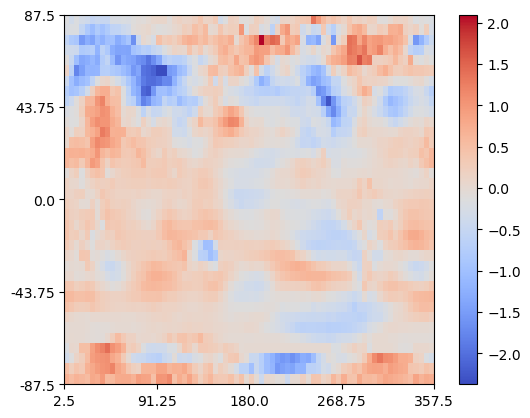}
    \label{oneaX}
  }
    \subfigure[Element--wise variances]{
    \includegraphics[width = 1.8in]{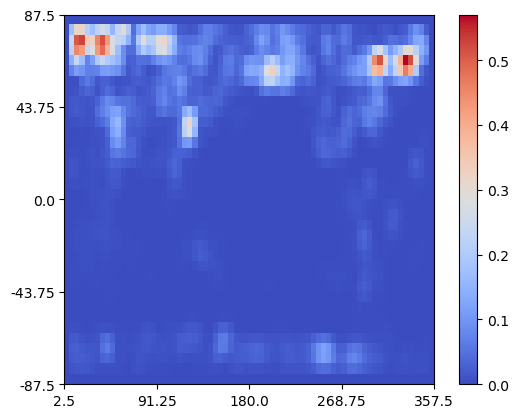}
    \label{oneb}
  }  
    \subfigure[Element--wise row--wise covariance]{
    \includegraphics[width = 1.8in]{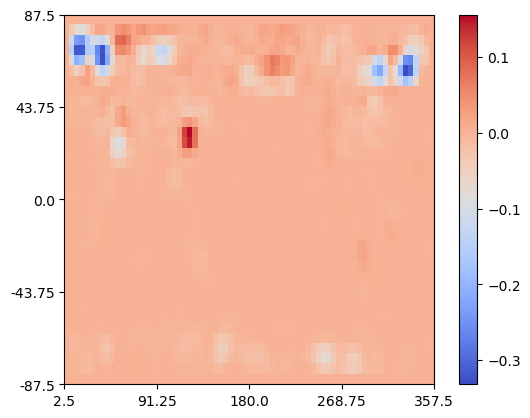}
    \label{oneb}
  }  
  \caption{Long--term mean temperature anomaly
  in February taken from the years 1971 to 2000.
  It also displays the estimated element--wise variances as well as element--wise row–wise covariances. 
  }
    \label{fig.real_life_data}
\end{figure}
\label{example.hetero}
\end{example}
To date, the statistics literature on non-stationary random fields is primarily focused on covariance modelling and estimation through parametric classes such as the Mat\'{e}rn class and Gneiting's model, or locally stationary extensions to those classes; see the literature review section in this introduction for detailed references.  Therefore the focus is mainly on parametric or semiparametric analysis of Gaussian non-stationary fields. However,
the aforementioned parametric/semiparametric models may be inadequate in
describing the complex dependence structure exhibited in modern-era data and consequently could lead to large bias in estimation and inference due to possible model misspecification. In those cases a more flexible nonparametric and non-Gaussian modelling of non-stationary random fields is desirable.  


The primary objective of this paper is to perform simultaneous nonparametric inference of a dependent and non-stationary random field with flexible nonparametric and non-Gaussian modeling of such field. To our knowledge, ours is the first attempt in the literature on simultaneous nonparametric inference of a two-dimensional non-stationary non-Gaussian random field. For presentational clarity and simplicity, we shall focus on the inference of the mean field though the inference of other quantities such as the covariance operator and the distributional field can be carried out using similar techniques. Our techniques originated mainly from highly nontrivial extensions of the literature in non-stationary time series analysis, where there has been a surge of interest in the flexible modeling of temporal heteroscedasticity in the last two decades. To mention a few examples, \cite{MR1429916}, \cite{zhou2009local}, and \cite{MR3920364} among others derived asymptotic results for a type of non-stationary time series known as the ``locally stationary process''; \cite{MR2724865} considered testing weak non-stationarity in a time series;
\cite{zhang2021debiased} adopted non-stationary errors in linear regression; while \cite{MR3798001} and \cite{MR3931381} tested structural changes in non-stationary time series.
Our work leverages concepts from non-stationary time series analysis literature to analyze a non-stationary random field.  In particular, the field is modeled by location-varying nonlinear transformations of an underlying shift random field, which allows very flexible changes (both smoothly and abruptly) in the data generating mechanism across different locations of the random field. As a theoretical contribution, we provide a uniform Gaussian approximation result for high dimensional sparse linear forms of the random field by which simultaneous confidence intervals and hypothesis tests can be implemented.

Despite wide applications of high dimensional Gaussian approximations to various types of data in the literature--such as \cite{MR3161448} for independent observations, \cite{MR3718156} for stationary time series, and \cite{MR4441125}, \cite{zhang2021debiased} for linear regression--Theorem \ref{theorem.Gaussian} in Section \ref{section.statistical_inference} is a highly nontrivial extension of these works to non-stationary random fields. In the time series literature, statisticians typically approximate data with “m-dependent” random variables and use the “big-block small-block” technique to establish Gaussian approximations, as demonstrated in \cite{MR3718156}. However, in our context, the dependence between random variables $X_{i_1}^{(j_1)}$ and $X_{i_2}^{(j_2)}$ depends on the indices $(i_1, j_1)$ and $(i_2, j_2)$. Due to this multidirectional dependence, novel m-dependent and martingale approximations are developed in this paper on two-dimensional neighbourhoods of each location. In particular, the ``blocks'' we construct are two-dimensional rectangles instead of directional line segments, as is common in time series literature.

Methodologically, this paper's main contribution lies in proposing a computationally-efficient locally weighted multiplier bootstrap algorithm as a unified tool for the simultaneous inference of a wide class of non-stationary random fields. To the best of our knowledge, the proposed locally weighted multiplier bootstrap is the first attempt towards bootstrap inference for a general class of non-stationary and non-Gaussian random fields. The proposed bootstrap leverages the ideas of the robust multiplier bootstrap for non-stationary time series introduced in \cite{zhou2013heteroscedasticity} and the dependent wild bootstrap for stationary random fields introduced by \cite{doi:10.1080/01621459.2023.2218578} and uses the convolution of weighted local sums of the random field and i.i.d. standard Gaussian random variables to mimic the joint covariance structure of high dimensional sparse linear forms on the random field. Furthermore, connections between dependent wild bootstrap and weighted local multiplier bootstrap are revealed in this paper via covariance decomposition of dependent Gaussian random fields, which could be of separate interest; see Remark \ref{rem:dwb} in Section \ref{section.methodlogy} for details. The theoretical verification of the bootstrap lies in a novel investigation into Heteroscedasticity and Autocorrelation Consistent (HAC) covariance  operator estimation and decomposition for high dimensional linear forms on non-stationary random fields, which generalizes the corresponding results for sequence data. Computationally, the weights in the proposed bootstrap can be calculated efficiently using Fast Fourier Transformations due to Szego's well-known results on approximate computation of large Toeplitz matrices \citep{gray2006toeplitz}. In particular, for a random field of size $(m,n)$, weight calculation for the proposed bootstrap requires $O(\max(n,m)\log(\max(n,m)))$ computation complexity, which is considerably faster than generating a $mn$-dimensional Gaussian random vector with desired covariance matrix in dependent wild bootstrap where the computation complexity is $O((mn)^3)$ using standard matrix square root algorithms.

\subsection{Literature Review}
There has been a long-lasting interest in the statistical inference of a random field, leading to abundant literature on the central limit theorem of stationary random fields under various dependence conditions.
 To list a few, see the works of \cite{MR672305}, \cite{MR1106280},  \cite{MR1616496}, \cite{MR375410}, \cite{MR2525996},
\cite{MR2988107}, 
and \cite{MR4430002}.
Results on spectral analysis of stationary random fields can be found in   \cite{MR3724221}, \cite{MR3892327}, and \cite{MR3949303}, among others.  
Due to the complex dependence structure of a random field, statistical inference of such data is typically carried out using resampling methods such as the bootstrap, where the current literature mainly focuses on stationary fields. See for instance \cite{MR2001447}, \cite{MR3654797}, \cite{doi:10.1080/01621459.2023.2218578}, and \cite{10.1214/21-EJS1959}. 
 We also refer the readers to book-length treatments of random field and spatial statistics in \cite{spatialstatistics}, \cite{Spatial},  and \cite{Theory_spatial}. 


As previously mentioned, the assumption of stationarity may be too restrictive in many real-life applications. There have been several attempts to overcome this limitation.
One approach involves modeling the covariances of a random field by parametric classes
such as the Mat\'{e}rn class and the Gneiting's model, or locally stationary extensions to those classes. See for instance \cite{MR2755012}, \cite{MR1941475}, \cite{MR2957287}, \cite{MR3328806} and \cite{MR3887154}, among others.
Parametric models are useful not only for discrete random field analysis but also useful for modeling the covariance function of
continuous random fields, as demonstrated by the works of \cite{MR2853727}, \cite{MR1731494}, \cite{MR3674697}, and \cite{MR2594414}. Other contributions include  \cite{MR4483222}
who adopted a Bayesian approach to estimating the covariances of a random field and \cite{MR4361753}  who separated a non-stationary
random field into several weakly stationary random fields before inference.

The rest of this paper is organized as follows: Section \ref{section.methodlogy} introduces the Nadaraya - Watson estimator for a smoothly varying mean field,
and presents the locally weighted multiplier bootstrap algorithm for statistical inference.
Section \ref{section.statistical_inference} introduces a general class of non-Gaussian and non-stationary random fields named the $(M,\alpha)-$ short-range dependent random fields. Several important and general theoretical results for this class of random fields are established, which may be of separate interest.
Section \ref{section.theoretical_justification} provides a theoretical justification for the proposed estimator and the bootstrap algorithm
with the help of the results in Section \ref{section.statistical_inference}.
Section \ref{section.numerical} verifies the validity of the proposed algorithms via simulated data and provides some real-life applications. We postpone the detailed proofs of the proposed theorems to the online supplements.

\section{Methodology}
\label{section.methodlogy}
\subsection{Point estimation}
Aligned with the introduction, this section conducts simultaneous statistical inference on the mean field of a non-stationary random field.  Suppose the observed data $X_i^{(j)}, i = 1,\cdots, n, j = 1,\cdots, m$ are stemmed from a two-dimensional random field,
\begin{equation}
X_i^{(j)} = \mu\left(\frac{i}{n}, \frac{j}{m}\right) + \epsilon_i^{(j)},
\end{equation}
where $\mu(\cdot,\cdot):[0,1]\times [0,1]\to \mathbf{R}$ is the mean field and is assumed to be a continuous function, $\{\epsilon_i^{(j)}\}$ is the non-stationary noise field with $\mathbf{E}\epsilon_i^{(j)}=0$, and $\min(n,m)\to\infty$. In this paper the data generating mechanism of the noise field can be flexibly changing across locations, both smoothly and abruptly. We are interested in the function value $\mu(x_v,y_v)$ for $V:=V(n,m)$  given positions $(x_v, y_v)$, where $x_v,y_v\in[0,1]$.  Our work allows for $V$ diverging to infinity as a function of $n$ and $m$ (possibly much faster than $m\times n)$, in which case the set $\{(x_v, y_v), v=1,2,\cdots, V\}$ forms a dense subset of $[0,1]^2$. We are interested in the fundamental problem of constructing simultaneous confidence regions $L(x_v,y_v)$, $v=1,2,\cdots, V$ for $\mu(\cdot,\cdot)$ such that 
\begin{align*}
    Prob\left(
    \mu(x_v,y_v)\in L(x_v,y_v)\ \text{for }v = 1,\cdots, V
    \right)\to 1 - \alpha\text{ as } n,m\to\infty
\end{align*}
for any given coverage level $1-\alpha$. 

To achieve the goal, suppose $G(\cdot):[-1,1]\to[0,\infty)$ is a given kernel function and $\mathcal{K} = \mathcal{K}(n,m)$ is a chosen bandwidth. Define
$p_v = \lfloor nx_v\rfloor$ and $q_v = \lfloor my_v\rfloor$, where $\lfloor x\rfloor$ stands for the largest integer that is smaller than or equal to $x$. Define the Nadaraya--Watson estimator, as introduced by \cite{doi:10.1137/1109020}, as follows: 
\begin{equation}
\widehat{\mu}(x_v, y_v) = \frac{\sum_{i = p_v - \mathcal{K}}^{p_v + \mathcal{K}}\sum_{j = q_v - \mathcal{K}}^{q_v + \mathcal{K}}X_i^{(j)}\times
G\left(\frac{i - p_v}{\mathcal{K}}\right) G\left(\frac{j - q_v}{\mathcal{K}}\right)} {\sum_{i = p_v - \mathcal{K}}^{p_v + \mathcal{K}}\sum_{j = q_v - \mathcal{K}}^{q_v + \mathcal{K}}G\left(\frac{i - p_v}{\mathcal{K}}\right) G\left(\frac{j - q_v}{\mathcal{K}}\right)}.
\label{eq.define_Nadaraya}
\end{equation}
Note that the estimator $\widehat{\mu}(x_v, y_v)$ is a weighted average of grid points in the neighborhood of $(x_v, y_v)$. 
We assume that the points of interest $(x_v, y_v)$ are not on the boundary of the domain; i.e., we assume that $2\mathcal{K} + 1\leq p_v\leq n - 2\mathcal{K}$ and
$2\mathcal{K} + 1\leq q_v\leq m - 2\mathcal{K}$.

\subsection{The Locally Weighted Multiplier Bootstrap}\label{sec:lwmb}
After deriving $\widehat{\mu}(x_v, y_v)$,  the next step involves constructing simultaneous confidence regions and conducting hypothesis testing for $\mu(x_v,y_v)$. Due to the potentially heterogeneous and dependent noises $\epsilon_i^{(j)}$, the estimators $\widehat{\mu}(x_v, y_v), v = 1,\cdots, V$ may have complex variances and covariances. Directly estimating these terms is difficult; therefore, we propose a bootstrap algorithm,
named ``locally weighted multiplier bootstrap'', to assist statistical inference.  
To describe our algorithm, we introduce several additional notations:  let the row--wise kernel matrix $K_r = \left\{K\left(\frac{i_1 - i_2}{\mathcal{B}}\right)\right\}_{i_1,i_2 = 1,\cdots, n}$ and the column--wise kernel matrix $K_c = \left\{K\left(\frac{j_1 - j_2}{\mathcal{B}}\right)\right\}_{j_1,j_2 = 1,\cdots, m}$, where $K(\cdot):\mathbf{R}\to[0,1]$ is the variance kernel function satisfying Definition \ref{definition.kernel} in Section \ref{section.statistical_inference}. Define the value
\begin{equation}
c_{i, v}^{(j)} = \begin{cases}
\frac{G\left(\frac{i - p_v}{\mathcal{K}}\right) G\left(\frac{j - q_v}{\mathcal{K}}\right)}
{\sqrt{\sum_{i = p_v - \mathcal{K}}^{p_v + \mathcal{K}}\sum_{j = q_v - \mathcal{K}}^{q_v + \mathcal{K}}G^2\left(\frac{i - p_v}{\mathcal{K}}\right) G^2\left(\frac{j - q_v}{\mathcal{K}}\right)}}\ \text{if } \vert i - p_v\vert\leq \mathcal{K}\ \text{and } \vert j - q_v\vert\leq \mathcal{K},\\
0\ \text{otherwise}.
\end{cases}
\label{eq.def_cij}
\end{equation}
According to Remark \ref{Remark.positive}, $K_r$ and $K_c$ are both Toeplitz and symmetric positive semi-definite matrices, so $\exists$ symmetric matrices $Q^{(n)} = (q^{(n)}_{ij})\in\mathbf{R}^{n\times n}$ and $Q^{(m)} = (q^{(m)}_{ij})\in\mathbf{R}^{m\times m}$ such that $K_r = Q^{(n)2}$ and $K_c = Q^{(m)2}$,
i.e., $K\left(\frac{i - j}{\mathcal{B}}\right) = \sum_{p = 1}^n q^{(n)}_{ip}q_{jp}^{(n)}$, and $K\left(\frac{i - j}{\mathcal{B}}\right) = \sum_{s = 1}^m q^{(m)}_{is}q^{(m)}_{js}$. With these notations we can propose our bootstrap algorithm.

\begin{remark}
 For large $m$ and $n$, the square root matrices $Q^{(n)}$ and $Q^{(m)}$ can be computed efficiently using Szego's well-known results on approximate computation of large Toeplitz matrices. See for instance the introduction in \cite{gray2006toeplitz}. In particular, $Q^{(n)}$ ($Q^{(m)}$) can be approximated by calculating Fourier transformations of $K_r$ ($K_c$), taking square roots of the Fourier coefficients, and then performing inverse Fourier transforms. The computational complexity of such operations is $O(n\log n)$ and $O(m\log m)$ for $Q^{(n)}$ and $Q^{(m)}$, respectively.
\end{remark}

\begin{algorithm}[Locally weighted multiplier bootstrap]
$\\$
\noindent
\textbf{Input: } The observed data $X_i^{(j)}, i = 1,\cdots, n, j = 1,\cdots, m$, the kernel function $G(\cdot)$, the variance kernel function $K(\cdot)$, the bandwidth $\mathcal{K}$ and the variance's bandwidth $\mathcal{B}$. The locations of interest $(x_v,y_v), v = 1,\cdots, V$,  the number of bootstrap replicates $B$, the nominal coverage probability $1 - \alpha$.

\noindent \textbf{Extra input for testing: } a prespecified mean field $\mu_0(\cdot, \cdot)$ to be tested.

1. Derive the estimator $\widehat{\mu}(x_v, y_v)$ for $v = 1,\cdots, V$ as in eq.\eqref{eq.define_Nadaraya}. For $i = \mathcal{K} + 1,\cdots, n - \mathcal{K}$ and $j = \mathcal{K} + 1,\cdots, m - \mathcal{K}$, compute $\widehat{\epsilon}_i^{(j)} = X_i^{(j)} - \widehat{\mu}(\frac{i}{n}, \frac{j}{m})$.

2. Derive $q^{(n)}_{ij}$ and $q^{(m)}_{ij}$.

3. Generate i.i.d. standard normal random variables $e_p^{(q)*}$ for $p = 1,\cdots, n$ and $q = 1,\cdots, m$. 
For $v = 1,\cdots, V$, compute the bootstrap statistics
\begin{equation}
\widehat{\mu}^*(x_v, y_v) = \widehat{\mu}(x_v, y_v)  +\frac{1}{T_{n,m}}
    \sum_{p = 1}^n\sum_{q = 1}^m 
    \left(\sum_{i = -\mathcal{K}}^\mathcal{K}\sum_{j = -\mathcal{K}}^\mathcal{K} G\left(\frac{i}{\mathcal{K}}\right) G\left(\frac{j}{\mathcal{K}}\right)\widehat{\epsilon}_{i+p_v}^{(j+q_v)}q_{(i+p_v)p}^{(n)}q_{(j+q_v)q}^{(m)}\right)e_p^{(q)*},
\label{eq.derive_boot_estimator}
\end{equation}
where $T_{n,m}
= \sum_{i = - \mathcal{K}}^{\mathcal{K}}\sum_{j = - \mathcal{K}}^{\mathcal{K}}G\left(\frac{i}{\mathcal{K}}\right) G\left(\frac{j}{\mathcal{K}}\right)$ and $p_v = \lfloor nx_v\rfloor$, $q_v = \lfloor my_v\rfloor$. 

4. Calculate the weighted maximum 
\begin{align*}
    T^* = \max_{v = 1,\cdots, V}\frac{T_{n,m}}{B_{n,m}\times \widehat{\tau}_v}\vert \widehat{\mu}^*(x_v, y_v) - \widehat{\mu}(x_v, y_v)\vert,
\end{align*}
where $B_{n,m} = \sqrt{\sum_{i = p_v - \mathcal{K}}^{p_v + \mathcal{K}}\sum_{j = q_v - \mathcal{K}}^{q_v + \mathcal{K}}G^2\left(\frac{i - p_v}{\mathcal{K}}\right) G^2\left(\frac{j - q_v}{\mathcal{K}}\right)}$.
Our work considers $\widehat{\tau}_v = 1$, for a homogeneous version of the confidence region; and $\widehat{\tau}_v = \widehat{\sigma}_v^{1/3}$, where 
$$
\widehat{\sigma}_v = \sqrt{\sum_{i_1 = p_v - \mathcal{K}}^{p_v + \mathcal{K}}\sum_{j_1 = q_v - \mathcal{K}}^{q_v + \mathcal{K}}\sum_{i_2 = p_v - \mathcal{K}}^{p_v + \mathcal{K}}\sum_{j_2 = q_v - \mathcal{K}}^{q_v + \mathcal{K}} c_{i_1,v}^{(j_1)}c_{i_2,v}^{(j_2)}\widehat{\epsilon}_{i_1}^{(j_1)}\widehat{\epsilon}_{i_2}^{(j_2)}
\times K\left(\frac{i_1 - i_2}{\mathcal{B}}\right)K\left(\frac{j_1 - j_2}{\mathcal{B}}\right)},
$$
for a heterogeneous version.

5. Repeat step 3 - 4 $B$ times. Let $T^*_{(1)}\leq T^*_{(2)}\leq \cdots\leq T^*_{(B)}$ be the ordered statistics of $T^*$ over the $B$ repetitions. Let
$$
C_{1-\alpha}^* = T^*_{(t)},\ \text{where } t = \min\left\{t = 1,\cdots, B: \frac{t}{B}\geq  1 - \alpha\right\}
$$
be the $1-\alpha$ sample quantile.

6.a. (\textbf{for simultaneous confidence region construction}): The $1-\alpha$ simultaneous confidence region for $\mu(x_v, y_v), v = 1,\cdots, V$ is
$$
\left\{z = (z_1,\cdots, z_V)^T\in\mathbf{R}^V: \vert z_v - \widehat{\mu}(x_v, y_v)\vert\leq \frac{B_{n,m}\times \widehat{\tau}_v}{T_{n,m}}\times C^*_{1-\alpha}\ \text{for all }v\right\}.
$$

6.b. (\textbf{for testing}): Reject $H_0$ if
$$
\max_{v = 1,\cdots, V}\frac{T_{n,m}\times \vert \widehat{\mu}(x_v, y_v) - \mu_0(x_v, y_v)\vert}{B_{n,m}\times \widehat{\tau}_v} > C^*_{1 - \alpha}.
$$
\label{algorithm.bootstrap}
\end{algorithm}
We observe from equation \eqref{eq.derive_boot_estimator} that our bootstrap statistic is determined by a convolution of i.i.d. standard normal random variables and weighted local sums of the random field, $\sum_{i = -\mathcal{K}}^\mathcal{K}\sum_{j = -\mathcal{K}}^\mathcal{K} G\left(\frac{i}{\mathcal{K}}\right) G\left(\frac{j}{\mathcal{K}}\right)q_{(i+p_v)p}^{(n)}q_{(j+q_v)q}^{(m)}\widehat{\epsilon}_{i+p_v}^{(j+q_v)}$, where the weights are determined by the kernel $G(\cdot)$ and the square root matrices $Q^{(n)}$ and $Q^{(m)}$. Hence we name the procedure "locally weighted multiplier bootstrap".
Our work leverages the ideas of the robust multiplier bootstrap for non-stationary time series introduced in \cite{zhou2013heteroscedasticity} and the dependent wild bootstrap for stationary random fields introduced by \cite{doi:10.1080/01621459.2023.2218578}. To see why this bootstrap works for non-stationary random fields, define $Cov^*$ as the covariance in the bootstrap world (conditional on observed data), then  
\begin{align*}
    Cov^*(\frac{T_{n,m}}{B_{n,m}}(\widehat{\mu}^*(x_{v_1}, y_{v_1}) - \widehat{\mu}(x_{v_1}, y_{v_1})), \frac{T_{n,m}}{B_{n,m}}(\widehat{\mu}^*(x_{v_2}, y_{v_2}) - \widehat{\mu}(x_{v_2}, y_{v_2})))\\
    = \sum_{i_1 = p_{v_1} -\mathcal{K}}^{p_{v_1} +\mathcal{K}}\sum_{j_1 = q_{v_1} -\mathcal{K}}^{q_{v_1} + \mathcal{K}}\sum_{i_2 = p_{v_2} -\mathcal{K}}^{p_{v_2} + \mathcal{K}}\sum_{j_2 = q_{v_2} -\mathcal{K}}^{q_{v_2} + \mathcal{K}}c_{i_1, v_1}^{(j_1)}c_{i_2, v_2}^{(j_2)}\widehat{\epsilon}_{i_1}^{(j_1)}\widehat{\epsilon}_{i_2}^{(j_2)} K\left(\frac{i_1 - i_2}{\mathcal{B}}\right)K\left(\frac{j_1 - j_2}{\mathcal{B}}\right),
\end{align*}
which is a consistent estimator for the covariance between $\frac{T_{n,m}}{B_{n,m}}(\widehat{\mu}^*(x_{v_1}, y_{v_1}) - \widehat{\mu}(x_{v_1}, y_{v_1}))$ and $ \frac{T_{n,m}}{B_{n,m}}(\widehat{\mu}^*(x_{v_2}, y_{v_2}) - \widehat{\mu}(x_{v_2}, y_{v_2}))$ for a wide range of non-stationary random fields according to Theorem \ref{theorem.covariances} in Section \ref{section.statistical_inference}.   Therefore, the validity of Algorithm \ref{algorithm.bootstrap} arises from simulating joint normal random variables on the random field with heteroscedasticity and autocorrelation consistent covariances that well approximate those of the Nadaraya--Watson estimators. The above result is potentially surprising since intuitively it is difficult to consistently mimic the joint covariance structure of the local nonparametric estimates when there are possibly many abrupt changes in the covariance structure of the random field in every local neighborhood.


\begin{remark}\label{rem:dwb}
In this remark we shall delve into the connection between our locally weighted multiplier bootstrap and the dependent wild bootstrap where the  auxiliary Gaussian random variables are correlated. Define $f_i^{(j)} = \sum_{p = 1}^n\sum_{q = 1}^m q_{ip}^{(n)}q_{jq}^{(m)} e_p^{(q)*}$. Note that $f_i^{(j)}$ are correlated Gaussian random variables with 
\begin{equation}
\begin{aligned}
    \mathbf{E}f_{i_1}^{(j_1)}f_{i_2}^{(j_2)} = 
    \sum_{p = 1}^n\sum_{q = 1}^m q_{i_1p}^{(n)}q_{i_2p}^{(n)}q_{j_1q}^{(m)}q_{j_2q}^{(m)} = K\left(\frac{i_1 - i_2}{\mathcal{B}}\right)K\left(\frac{j_1 - j_2}{\mathcal{B}}\right).
\end{aligned}
\label{eq.generate_kernel}
\end{equation}
Then it is easy to observe that our bootstrap statistics can be written as
\begin{equation}
\begin{aligned}
    \widehat{\mu}^*(x_v, y_v) = \widehat{\mu}(x_v, y_v)  + \frac{1}{T_{n,m}}\sum_{i =  - \mathcal{K}}^{ \mathcal{K}}\sum_{j = - \mathcal{K}}^{\mathcal{K}}\widehat{\epsilon}_{i + p_v}^{(j + q_v)*}
G\left(\frac{i}{\mathcal{K}}\right) G\left(\frac{j}{\mathcal{K}}\right)f_{i + p_v}^{(j + q_v)}.
\end{aligned}
\label{eq.derive_boot_estimator}
\end{equation}
Note that \eqref{eq.derive_boot_estimator} is in the form of the dependent wild bootstrap; that is, a weighted sum of correlated Gaussian random variables multiplied by the observed values of the random field. On the other hand, observe that any correlated Gaussian random field can be written as weighted sums of i.i.d. Gaussian random variables via eigen-decompositions of its covariance operator. Therefore a simple switching of summation order will yield that a dependent wild bootstrap procedure can also be written in the form of weighted sums of i.i.d. standard Gaussian random variables. The above connection is potentially interesting as it reveals that the classic multiplier bootstrap with independent auxiliary random variables and the dependent wild bootstrap can be viewed in a unified way so that strengths can be borrowed from each other. For instance, the heteroscedasticity and autocorrelation robust property of the locally weighted multiplier bootstrap can shed light on whether the dependent wild bootstrap can be made to have such robustness. Reversely, the dependent wild bootstrap is known to be easily implementable for random fields with missing observations \citep{doi:10.1080/01621459.2023.2218578}. The above-mentioned connection could provide some inspiration on how to design weighted multiplier bootstrap in the missing data scenario.

\end{remark}

Note that $\hat\sigma_v$ is an estimate of the standard deviation of $\widehat{\mu}(x_v, y_v)$. The utilization of $\widehat{\sigma}_v$ aligns with Theorem\ref{theorem.covariances} concerning the consistent estimation of the variances of $\widehat{\mu}(x_v, y_v)$. The index $1/3$ is introduced to accommodate possible large fluctuations in the estimation of the variances of $\widehat{\mu}(x_v, y_v)$ which could lead to inferior coverage probabilities and testing results in moderate sample sizes. The heterogeneous version of algorithm \ref{algorithm.bootstrap} allows one to visualize the heterogeneity of the variability of the estimated mean field across locations.

\section{The Non-stationary Random Field Model}
\label{section.statistical_inference}
This section aims to establish a mathematical model for the dependent and non-stationary random noises $\epsilon_i^{(j)}$ as discussed in section \ref{section.methodlogy}. To accomplish this, our work leverages and extends the model of \cite{MR2988107} to account for non--stationary noises. Specifically, suppose $e_i^{(j)}, i,j\in\mathbf{Z}$ are independent (but not necessarily identically distributed) random variables, assume that  $\epsilon_i^{(j)}, i = 1,2,\cdots, n, j = 1,\cdots, m$ satisfy
\begin{equation}
\epsilon_i^{(j)} = G_{i, n}^{(j), m}(e_{i - u}^{(j - v)}; u,v\in\mathbf{Z}).
\label{eq.def_epsilon}
\end{equation}
In other words, we assume that $\epsilon_i^{(j)}$ is generated by a location-varying nonlinear transformation of the shift random field $\{e_{i - u}^{(j - v)}; u,v\in\mathbf{Z}\}$.
Observe that the function $G_{i, n}^{(j), m}$ depends on $i,j,n,m$ which implies that the nonlinear transformation may vary very flexibly with respect to the position $(i,j)$ and the sample sizes $n$ and $m$. Consequently, the joint distribution of the random field $\{\epsilon_i^{(j)}\}$ can be highly non-stationary and non-Gaussian with complex covariances. 

Next, we provide a dependence measure for the above-defined non-stationary random field. Define $e_i^{(j)\dagger},i,j\in\mathbf{Z}$ as independent random variables such that
$e_i^{(j)\dagger}$ is independent of $e_p^{(q)}$ for arbitrary $i,j,p,q\in\mathbf{Z}$; and  $e_i^{(j)\dagger}$ has the same distribution as
$e_i^{(j)}$. For any $p,q\in\mathbf{Z}$, define
\begin{align*}
\epsilon_{i, p}^{(j), (q)}
=
G_{i, n}^{(j), m}\left(e_{i - u}^{(j - v)\Delta}; u,v\in\mathbf{Z}\right), \text{where } e_{i - u}^{(j - v)\Delta} =
\begin{cases}
e_{i - u}^{(j - v)\dagger}\ \text{if } i - u = i + p\ \text{and } j - v = j + q,\\
e_{i - u}^{(j - v)}\ \text{otherwise},
\end{cases}
\end{align*}
that is, replacing the random variable $e_{i+p}^{(j+q)}$ in $\epsilon_i^{(j)}$ with $e_{i+p}^{(j+q)\dagger}$. For a fixed number $M \geq 1$, define
$\delta_{i, p, (M)}^{(j), (q)} = \Vert\epsilon_i^{(j)} - \epsilon_{i, p}^{(j), (q)} \Vert_M$ and the dependence measure
\begin{equation}
\delta_{p}^{(q)} = \max_{i = 1,\cdots, n,\ j = 1,\cdots m}\delta_{i, p, (M)}^{(j), (q)},
\label{eq.def_delta}
\end{equation}
where $\Vert\cdot\Vert_M = (\mathbf{E}\vert\cdot\vert^M)^{1/M}$. We omit `$M$' on the left-hand side of eq.\eqref{eq.def_delta}. However, $M$ is considered to be fixed in our work, so this abuse of notation should not introduce confusion. Note that the dependence measure $\delta_{p}^{(q)}$ measures the maximum impact (in $L^{M}$ norm) the system will experience if an innovation or shock of the system of distance $(p,q)$ away from the current location is changed to an i.i.d. copy. Ideally one would hope that $\delta_p^{(q)}$ would diminish as $\max(p,q)\rightarrow\infty$. With the aforementioned notation, we are now able to introduce a concept, named ``($M,\alpha$)-short range dependent random field'', as follows.

\begin{definition}[($M,\alpha$)-short range dependent random field]
Suppose random variables $\epsilon_i^{(j)}, i = 1,\cdots, n, j = 1,\cdots, m$ satisfy eq.\eqref{eq.def_epsilon}. In addition, suppose
\begin{align*}
\mathbf{E}\epsilon_i^{(j)} = 0\ \text{for } \forall i,j,\
\max_{i = 1,\cdots, n, j = 1,\cdots, m}\Vert\epsilon_i^{(j)}\Vert_M = O(1)\\
 \text{and }
\sup_{p,q\in\mathbf{Z}}(1 + \vert p\vert + \vert q\vert)^\alpha\delta_p^{(q)} = O(1).
\end{align*}
Then we call $\epsilon_i^{(j)}$ a $(M,\alpha)-$short range dependent random field.
\label{definition.Ma}
\end{definition}

\begin{example}[Linear random field]
Suppose the non-stationary random field $\epsilon_i^{(j)}$ has a linear form
$$
\epsilon_i^{(j)} = \sum_{u = -\infty}^\infty\sum_{v = -\infty}^\infty b_{u, i}^{(v), (j)}e_{i - u}^{(j - v)}.
$$
Direct calculations using the definition yield that $\delta_{i, p, (M)}^{(j), (q)}\leq 2\Vert e_{i + p}^{(j + q)}\Vert_M\times \vert b_{-p, i}^{(-q), (j)}\vert$. Therefore Definition \ref{definition.Ma} is ensured if 
$\max_{i = 1,\cdots, n, j = 1,\cdots, m}\vert b_{p, i}^{(q), (j)}\vert\leq \frac{C}{(1 + \vert p\vert + \vert q\vert)^\alpha}$.
In other words, when considering a linear random field, $\delta_p^{(q)}$ in eq.\eqref{eq.def_delta} is proportional to the maximum magnitude of the corresponding coefficients of the linear model. Consequently the $(M,\alpha$)-short range dependent condition is easy to check.
\end{example}

Definition \ref{definition.Ma} requires  $\epsilon_i^{(j)}$ to have $0$ mean, finite $M$-th moment, and satisfy the so--called ``short-range dependent'' condition. The first two are common in the high-dimensional statistics literature, such as \cite{zhang2021debiased}, \cite{MR4441125}, and \cite{MR4206676}. We now elaborate on the third condition. According to Remark \ref{remark.covariance} below, auto-covariances of $\epsilon_i^{(j)}$  decay at a polynomial rate with respect to their distance. Hence, we identify $\alpha$ as an index governing the strength of dependence within the random field. In particular, when $\alpha$ is significantly greater than $3$, the covariances between $\epsilon_{i_1}^{(j_1)}$ and $\epsilon_{i_2}^{(j_2)}$  shrinks rapidly to 0 as $\max(\vert i_1 - i_2\vert, \vert j_1 - j_2\vert)\rightarrow\infty$, indicating a weak dependence. 
\begin{remark}
As demonstrated in corollary \ref{corollary.covariance} in the online supplement, if $\epsilon_i^{(j)}$ are $(M,\alpha)-$short range dependent random variables,
their covariances should satisfy $$
\vert\mathbf{E}\epsilon_{i_1}^{(j_1)}\epsilon_{i_2}^{(j_2)}\vert\leq C\times (1 + \max(\vert i_1 - i_2\vert, \vert j_1 - j_2\vert))^{3-\alpha}.
$$
\label{remark.covariance}
\end{remark}

The remainder parts of this section derive theoretical results concerning linear forms of the non-stationary random field $\epsilon_i^{(j)}$, which is essential for our analysis of the local nonparametric estimates and the locally weighted multiplier bootstrap. In statistical applications such as linear regression, time series analysis, and statistical learning, the consistency results as well as central limit theorems of many statistics rely on the analysis of linear forms of random variables, as demonstrated in  \cite{MR4441125}, \cite{MR1093459}, and \cite{MR1385671}. Therefore, this section could be of independent interest. Our first result involves the concentration inequality for linear combinations of $(M,\alpha)-$short range dependent random variables. \cite{MR0133849} derived similar results for independent random variables, whereas our research extends this analysis to non-stationary random fields.
\begin{lemma}
Suppose $\epsilon_i^{(j)}, i = 1,\cdots, n, j = 1,\cdots, m$ are $(M,\alpha)-$short range dependent random variables with $M>4, \alpha > 3$; and $a_i^{(j)}, i = 1,\cdots, n, j = 1,\cdots, m$ are real numbers. Then there exists a constant $C$ independent of $a_i^{(j)}$ such that
\begin{equation}
\begin{aligned}
\Vert\sum_{i = 1}^n\sum_{j = 1}^m a_i^{(j)}\epsilon_i^{(j)}\Vert_M\leq C\sqrt{\sum_{i = 1}^n\sum_{j = 1}^m a_i^{(j)2}}
\end{aligned}
\label{eq.three_lemma}
\end{equation}
\label{lemma.consistent_linear_combination}
\end{lemma}

The next result derives a Gaussian approximation theorem for high dimensional sparse linear combinations of $\epsilon_i^{(j)}$. Compared to the classical central limit theorems such as those discussed in \cite{MR2002723}, the Gaussian approximation theorem allows for a diverging number of linear combinations, which is favorable in the realm of high-dimensional statistics. Our result extends the high dimensional Gaussian approximation results of \cite{MR3161448}, \cite{MR4441125}, and \cite{MR3718156} to non-stationary random fields. 
\begin{theorem}
Suppose $\epsilon_i^{(j)}, i = 1,\cdots,n, j = 1,\cdots, m$ are $(M,\alpha)-$short range dependent random variables with $M > 4$ and $\alpha > 3$, and $\mathcal{K} = k(n,m)$ is an integer satisfying $\mathcal{K}\to\infty$ and $\mathcal{K}/n\to 0$, $\mathcal{K}/ m\to 0$ as $\min(n,m)\to\infty$. Suppose $V$ is an integer satisfying $V = O\left(\mathcal{K}^{\alpha_V}\right)\ \text{with } 0\leq \alpha_V < \frac{(\alpha - 3)M}{1 + 5(\alpha - 3)}$. For $v = 1,\cdots, V$, suppose $x_v,y_v$ are two integers such that $\mathcal{K} + 1\leq x_v\leq n - \mathcal{K}$ and $\mathcal{K} + 1\leq y_v\leq m -\mathcal{K}$.
Suppose $a_{i, v}^{(j)}, i = 1,\cdots, n, j = 1,\cdots, m, v = 1,\cdots, V$ are real numbers satisfying the following conditions:
\begin{equation}
\begin{aligned}
\sum_{i = 1}^n\sum_{j = 1}^m a_{i,v}^{(j)2}>0\ \text{for all } v,\\
\exists\ \text{real numbers } \lambda_{i,v}, \tau_{j,v}\geq 0\ \text{such that }
\frac{\vert a_{i,v}^{(j)}\vert}{\sqrt{\sum_{i = 1}^n\sum_{j = 1}^m a_{i,v}^{(j)2}}}\leq \lambda_{i,v}\times \tau_{j,v}\ \text{for all }i,j,v,\\ 
\sum_{i = 1}^n \lambda_{i,v}^2\leq C,\ \sum_{j = 1}^m \tau_{j,v}^2\leq C, \vert\lambda_{i,v}\vert\leq C/\sqrt{\mathcal{K}},\ \vert\tau_{j,v}\vert\leq C/\sqrt{\mathcal{K}},\\
a_{i,v}^{(j)} = 0\ \text{if } \vert i - x_v\vert > \mathcal{K}\ \text{or } \vert j - y_v\vert > \mathcal{K},
\end{aligned}
\label{eq.condition_T1}
\end{equation}
where $C>0$ is a fixed constant, and $k(n,m)$ is a function of $n$ and $m$. In addition, suppose there exists a constant $c_0>0$ such that
\begin{equation}
\begin{aligned}
\sum_{i_1 = 1}^n\sum_{j_1 = 1}^m\sum_{i_2 = 1}^n\sum_{j_2  = 1}^m b_{i_1}^{(j_2)}b_{i_2}^{(j_2)}(\mathbf{E}\epsilon_{i_1}^{(j_1)}\epsilon_{i_2}^{(j_2)})\geq c_0\sum_{i = 1}^n\sum_{j = 1}^m b_{i}^{(j)2}
\end{aligned}
\label{eq.covariance_bij}
\end{equation}
for any $b_i^{(j)}\in\mathbf{R}, i = 1,\cdots, n, j = 1,\cdots, m$. Define $R_v = \sum_{i= 1}^n\sum_{j = 1}^ma_{i, v}^{(j)}\epsilon_{i}^{(j)}$ for $v = 1,\cdots, V$, then
\begin{equation}
\begin{aligned}
\sup_{x\in\mathbf{R}}\Big\vert
Prob\left(\max_{v = 1,\cdots, V}\frac{\vert R_v\vert}{\sqrt{\sum_{i = 1}^n\sum_{j = 1}^m a_{i, v}^{(j)2}}}\leq x\right)
- Prob\left(\max_{v = 1,\cdots, V}\vert\xi_v\vert\leq x\right)\Big\vert = o(1),
\end{aligned}
\label{eq.Gaussian_result}
\end{equation}
where $\xi_v, v = 1,\cdots, V$ are joint normal random variables such that $\mathbf{E}\xi_v = 0$ and
\begin{align*}
\mathbf{E}\xi_{v_1}\xi_{v_2} = \frac{1}{\sqrt{\sum_{i = 1}^n\sum_{j = 1}^m a_{i, v_1}^{(j)2}}\times \sqrt{\sum_{i = 1}^n\sum_{j = 1}^m a_{i, v_2}^{(j)2}}}\mathbf{E}R_{v_1}R_{v_2}.
\end{align*}
\label{theorem.Gaussian}
\end{theorem}
\begin{example}
This example provides illustrations of conditions in Theorem \ref{theorem.Gaussian}, utilizing the  Nadaraya--Watson estimator in Section \ref{section.methodlogy} as a demonstrative example. Notice that
\begin{align*}
\frac{T_{n,m}}{B_{n,m}}(\widehat{\mu}(x_v, y_v) - \mu(x_v, y_v)) = bias + \sum_{i = -\mathcal{K}}^\mathcal{K}\sum_{j = -\mathcal{K}}^\mathcal{K}\frac{G\left(\frac{i}{\mathcal{K}}\right)G\left(\frac{j}{\mathcal{K}}\right)}{B_{n,m}}\epsilon_{i + p_v}^{(j + q_v)}.
\end{align*}
Denoting $c_{i,v}^{(j)}$ as in eq.\eqref{eq.def_cij}, then the second term on the right hand side of the above equation becomes $\sum_{i = 1}^n\sum_{j = 1}^m c_{i,v}^{(j)}\epsilon_i^{(j)}$. Because $\sum_{i = 1}^n\sum_{j = 1}^m c^{(j)2}_{i,v} = 1$; $c_{i,v}^{(j)} = \frac{G\left(\frac{i}{\mathcal{K}}\right)}{\sqrt{\sum_{i = -\mathcal{K}}^{\mathcal{K}}G^2\left(\frac{i}{\mathcal{K}}\right)}}\times \frac{G\left(\frac{j}{\mathcal{K}}\right)}{\sqrt{\sum_{i = -\mathcal{K}}^{\mathcal{K}}G^2\left(\frac{i}{\mathcal{K}}\right)}}$; and $G(\cdot)$ is bounded, it is easy to check that conditions in Theorem \ref{theorem.Gaussian} are all satisfied. Therefore, in nonparametric regression, achieving conditions in Theorem \ref{theorem.Gaussian} is straightforward. Notably, $\mathcal{K}$ in Theorem \ref{theorem.Gaussian}, as well as $\mathcal{B}$ in Theorem \ref{theorem.covariances}, coincide with $\mathcal{K}$ and $\mathcal{B}$ in Section \ref{section.methodlogy} after choosing $c_{i,v}^{(j)}$.
\label{example.expand}
\end{example}

\begin{remark}
\label{remark.meaning_cov}
We introduce eq.\eqref{eq.condition_T1} to accommodate the setup of nonparametric regression, where the weights outside a neighborhood of the target positions $(x_v,y_v)$ are set to 0 (so the bias remains small when  applying Theorem \ref{theorem.Gaussian} in Section \ref{section.methodlogy}), and the weights within the neighborhood of the target positions are distributed to avoid concentrating in one or a few small areas. Consider a counter example where $a_{x_v, v}^{(y_v)} = 1$ and $a_{i,v}^{(j)} = 0$ otherwise, then $R_v = \epsilon_{x_v}^{(y_v)}$, and the Gaussian approximation theorem would not hold true. The adoption of $\lambda_{i,v}$ and $\tau_{j,v}$ also prevents the weights from concentrating along a single axis.

Eq.\eqref{eq.covariance_bij} is equivalent to
$
\mathbf{E}\left(\sum_{i = 1}^n\sum_{j = 1}^mb_i^{(j)}\epsilon_i^{(j)}\right)^2\geq c_0\sum_{i = 1}^n\sum_{j = 1}^m b_i^{(j)2}
$. That is, we expand the matrix $\{\epsilon_i^{(j)}\}_{i = 1,\cdots, n, j = 1,\cdots, m}$ into an $n\times m$ column vector, and the minimum eigenvalue of the covariance matrix of this vector should not shrink to 0.  
\end{remark}

Our final key result in this section provides Heteroscedasticity and Autocorrelation Consistent (HAC) covariance estimation for the high dimensional linear forms $\sum_{i = 1}^n\sum_{j = 1}^m a_{i,v}^{(j)}\epsilon_i^{(j)}$, $v=1,\cdots, V$, which extends the corresponding results in time series to non-stationary random fields. Definition \ref{definition.kernel} introduces several requirements for a desirable kernel.
\begin{definition}[Variance kernel function]
Suppose the function $K\left(\cdot\right):\mathbf{R}\to [0,1]$ is symmetric, continuous differentiable, $K(0) = 1$, and $K(x)$ is decreasing on $[0,\infty)$. In addition, define the Fourier transformation of $K$ as $\mathcal{F}K(x) = \int_{\mathbf{R}} K(t)\exp(-2\pi\mathrm{i}tx)\mathrm{d}t$. Suppose $\mathcal{F}K(x)\geq 0$ for all $x\in\mathbf{R}$, $\int_{0}^\infty K^2(x)dx<\infty$, and $\int_{0}^\infty x\times K(x)dx < \infty$, then we call $K$ the variance kernel function.
\label{definition.kernel}
\end{definition}

Since
$
\int_{0}^\infty K(x)dx \leq \int_{0}^1 K(x)dx + \int_{1}^\infty x\times K(x)dx < \infty
$, Definition \ref{definition.kernel} implies $\int_{0}^\infty K(x)dx < \infty$.
\begin{remark}
\label{Remark.positive}
Definition \ref{definition.kernel} can be achieved by selecting $K(x) = \exp\left(-x^2 / 2\right)$. According to \cite{MR2656050} and \cite{zhang2021debiased}, the matrix $\{K\left(\frac{i - j}{k}\right)\}_{i,j = 1,\cdots, n}$ is positive semi--definite for all $k > 0$ if $\mathcal{F}K(x)\geq 0$. 
\end{remark}

\begin{theorem}[HAC Covariance Estimation for Non-stationary Random Fields]
Suppose random variables $\epsilon_i^{(j)}$ are $(M,\alpha)-$short range dependent random variables with $M>4$ and $\alpha > 5$, and the linear combination coefficients $a_{i,v}^{(j)},\ v = 1,\cdots, V$ satisfy conditions in Theorem \ref{theorem.Gaussian}. Suppose $K(\cdot)$ is a variance kernel function and $\mathcal{B} = \mathcal{B}(n,m)$ is a bandwidth that satisfies $\mathcal{B}\to\infty$ and $\frac{V^{4/M}\mathcal{B}^3}{\mathcal{K}}\to 0$ as $\min(n,m)\to\infty$. Then
\begin{equation}
\begin{aligned}
\Big\Vert\ \max_{v_1, v_2 = 1,\cdots, V}
\vert
\sum_{i_1 = 1}^n\sum_{j_1 = 1}^m\sum_{i_2 = 1}^n\sum_{j_2 = 1}^m c_{i_1,v_1}^{(j_1)}c_{i_2,v_2}^{(j_2)}\epsilon_{i_1}^{(j_1)}\epsilon_{i_2}^{(j_2)}
\times K\left(\frac{i_1 - i_2}{\mathcal{B}}\right)K\left(\frac{j_1 - j_2}{\mathcal{B}}\right)\\
- Cov\left(
\sum_{i_1 = 1}^n\sum_{j_1 = 1}^m c_{i_1,v_1}^{(j_1)}\epsilon_{i_1}^{(j_1)}, \sum_{i_2 = 1}^n\sum_{j_2 = 1}^m c_{i_2,v_2}^{(j_2)}\epsilon_{i_2}^{(j_2)}
\right)
\vert\ \Big\Vert_{M/2}
= O\left(S(\mathcal{B}) + \frac{V^{4/M}\mathcal{B}^3}{\mathcal{K}}\right),
\end{aligned}
\label{eq.variance}
\end{equation}
where $c_{i,v}^{(j)} = a_{i,v}^{(j)} / \sqrt{\sum_{i = 1}^n\sum_{j = 1}^m a_{i,v}^{(j)2}}$. $S(\mathcal{B}) = \mathcal{B}^{5 - \alpha}$ if $5 < \alpha < 6$;
$S(\mathcal{B}) = \log(\mathcal{B}) / \mathcal{B}$ if $\alpha = 6$; and $S(\mathcal{B}) = 1 / \mathcal{B}$ if $\alpha > 6$.
\label{theorem.covariances}
\end{theorem}
Estimating the covariances for the linear combinations $\sum_{i = 1}^n\sum_{j = 1}^m c_{i,v}^{(j)}\epsilon_{i}^{(j)}$, where $v = 1,\cdots, V$, is essential for conducting statistical inference in our case. Theorem \ref{theorem.covariances} offers a consistent procedure for estimating these covariances. Note that our HAC covariance estimator, i.e., the first term on the left hand side of \eqref{eq.variance} can be rewritten as
\begin{align*}
\sum_{s = 1 - n}^{n - 1}\sum_{t = 1 - m}^{m - 1} K\left(\frac{s}{\mathcal{B}}\right) K\left(\frac{t}{\mathcal{B}}\right)\sum_{i_1 = \max(1, 1 - s)}^{\min(n, n - s)}\sum_{j_1 = \max(1, 1 - t)}^{\min(m, m - t)}
c_{i_1,v_1}^{(j_1)}c_{i_1 + s,v_2}^{(j_1 + t)}\epsilon_{i_1}^{(j_1)}\epsilon_{i_1 + s}^{(j_1 + t)},
\end{align*}
which is a weighted average of sample covariances across  lags with slowly diverging order. Therefore, the form of eq.\eqref{eq.variance} is inspired by the HAC estimation method for non-stationary time series illustrated in \cite{b76ccb64-7fa5-32f9-a4fe-72298146be7d} and \cite{MR2748557}, and extends their approach to high-dimensional linear statistics of two--dimensional non-stationary random fields. 

According to our discussion in Section \ref{sec:lwmb}, the conditional covariance of our locally weighted multiplier bootstrap is of the form $$\sum_{i_1 = 1}^n\sum_{j_1 = 1}^m\sum_{i_2 = 1}^n\sum_{j_2 = 1}^m c_{i_1,v_1}^{(j_1)}c_{i_2,v_2}^{(j_2)}\epsilon_{i_1}^{(j_1)}\epsilon_{i_2}^{(j_2)}
\times K\left(\frac{i_1 - i_2}{\mathcal{B}}\right)K\left(\frac{j_1 - j_2}{\mathcal{B}}\right).$$
Therefore Theorem \ref{theorem.covariances} provides a key theoretical support for the consistency of our bootstrap procedure.

Notably, if one chooses $a_{i,v}^{(j)} = \frac{1}{T_{n,m}}G\left(\frac{i - p_v}{\mathcal{K}}\right)G\left(\frac{j - q_v}{\mathcal{K}}\right)$ for $\vert i - p_v\vert\leq \mathcal{K} \cap\vert j - q_v\vert\leq \mathcal{K}$, and $0$ otherwise in Section \ref{section.methodlogy}, then the coefficients $c_{i,v}^{(j)}$ in eq.\eqref{eq.def_cij} and $c_{i,v}^{(j)}$ in Theorem \ref{theorem.covariances} are identical. Thus, we employ consistent notation in these instances.

\section{Theory of Estimation and Bootstrap}
\label{section.theoretical_justification}
This section employs theoretical results presented in Section \ref{section.statistical_inference} to analyze the estimators in Section \ref{section.methodlogy}. Before presenting our works, we introduce some technical assumptions regarding the kernels, the mean field $\mu(\cdot,\cdot)$, and the errors $\epsilon_i^{(j)}$.

\textbf{Assumptions: }

1. $\mu(\cdot,\cdot)$ is continuous differentiable on $[0,1]\times[0,1]$, the kernel function $G(\cdot): [-1,1]\to [0,\infty)$ is continuous, symmetric, and is decreasing on $[0,1]$, $G(0) = 1$. $\mathcal{K}>0,\ \mathcal{K}\in\mathbf{Z}$ is the bandwidth satisfying $\mathcal{K}\to\infty$, and
$$
\frac{\mathcal{K}^2\sqrt{\log(\mathcal{K})}}{n}\to 0,\ \frac{\mathcal{K}^2\sqrt{\log(\mathcal{K})}}{m}\to 0\ \text{as $\min(n,m)\to\infty$}
$$

2. $V$ (the number of positions) satisfies $V = O(\mathcal{K}^{\alpha_V})$ with $0\leq \alpha_V < \frac{(\alpha - 3)M}{1 + 5(\alpha - 3)}$. The positions $(x_v,y_v)$ satisfy
$2\mathcal{K} + 1\leq p_v\leq n - 2\mathcal{K}$ and
$2\mathcal{K} + 1\leq q_v\leq m - 2\mathcal{K}$; where $p_v = \lfloor nx_v\rfloor$ and $q_v = \lfloor my_v\rfloor$.

3. $\epsilon_i^{(j)}, i = 1,\cdots, n, j = 1,\cdots, m$ are $(M,\alpha)-$short range dependent random variables with $M > 4$ and $\alpha > 5$. $\mathcal{B}$ is the variance bandwidth satisfying $\mathcal{B} = C\mathcal{K}^{\alpha_B}$ with $0 < \alpha_B < \frac{1}{3} - \frac{4\alpha_V}{3M}$, and $0<C<\infty$ is a constant. $K(\cdot)$ is a variance kernel function, i.e., $K(\cdot)$ satisfies Definition \ref{definition.kernel}.

4. The covariances of $\epsilon_i^{(j)}$ satisfy eq.\eqref{eq.covariance_bij} for any $b_i^{(j)}\in\mathbf{R}$.

5.(for heterogeneous version of Algorithm \ref{algorithm.bootstrap}) $\mathcal{S}(\mathcal{B})\times \mathcal{K}^{\frac{\alpha_V}{M}}\sqrt{\log(\mathcal{K})} = o(1)$ and $0\leq \alpha_B < \frac{1}{3} - \frac{5\alpha_V}{3M}$.

The online supplement proves that the conditions in Theorem \ref{theorem.Gaussian} and Theorem \ref{theorem.covariances} are met with the help of Assumptions 1 and 2. With these assumptions, we are able to derive the consistency and the Gaussian approximation theorem for the proposed estimator \eqref{eq.define_Nadaraya}.


\begin{theorem}
Suppose Assumptions 1 to 4 hold. Then

\begin{equation}
\Vert\ \max_{v = 1,\cdots, V}\vert\widehat{\mu}(x_v, y_v) - \mu(x_v, y_v)\vert\ \Vert_M = O\left(\frac{\mathcal{K}}{n} + \frac{\mathcal{K}}{m} + \frac{V^{1/M}}{\mathcal{K}}\right).
\end{equation}
Moreover, we have
\begin{equation}
\begin{aligned}
\sup_{x\in\mathbf{R}}\vert Prob\left(\max_{v = 1,\cdots, V}\frac{T_{n,m}}{B_{n,m}}\vert\widehat{\mu}(x_v, y_v) - \mu(x_v, y_v)\vert\leq x\right)\\
- Prob\left(\max_{v = 1,\cdots, V}\vert\xi_v\vert\leq x\right)\vert = o(1),
\label{eq.Gaussian_mu}
\end{aligned}
\end{equation}
where $B_{n,m} = \sqrt{\sum_{i =  - \mathcal{K}}^{\mathcal{K}}\sum_{j = - \mathcal{K}}^{\mathcal{K}}G^2\left(\frac{i}{\mathcal{K}}\right) G^2\left(\frac{j}{\mathcal{K}}\right)}$ and $(\xi_1,\cdots, \xi_V)$ are joint normal random variables with $\mathbf{E}\xi_i = 0$ and
\begin{equation}
\begin{aligned}
\mathbf{E}\xi_{v_1}\xi_{v_2} = \frac{1}{B_{n,m}^2}\sum_{i_1 = -\mathcal{K}}^{\mathcal{K}}\sum_{j_1 = -\mathcal{K}}^{\mathcal{K}}\sum_{i_2 = -\mathcal{K}}^{\mathcal{K}}
\sum_{j_2 = -\mathcal{K}}^{\mathcal{K}}\mathbf{E}\epsilon_{i_1 + p_{v_1}}^{(j_1 + q_{v_1})}\epsilon_{i_2 + p_{v_2}}^{(j_2 + q_{v_2})}\\
\times
G\left(\frac{i_1}{\mathcal{K}}\right)G\left(\frac{j_1}{\mathcal{K}}\right)G\left(\frac{i_2}{\mathcal{K}}\right)G\left(\frac{j_2}{\mathcal{K}}\right).
\end{aligned}
\label{eq.covariance_xi}
\end{equation}
In addition, suppose Assumption 5 holds true, then 
\begin{equation}
\begin{aligned}
\sup_{x\in\mathbf{R}}\vert Prob\left(\max_{v = 1,\cdots, V}\frac{T_{n,m}}{\widehat{\tau}_v B_{n,m}}\vert\widehat{\mu}(x_v, y_v) - \mu(x_v, y_v)\vert\leq x\right)\\
- Prob\left(\max_{v = 1,\cdots, V}\vert\frac{\xi_v}{\sigma_v^{1/3}}\vert\leq x\right)\vert = o(1),
\end{aligned}
\label{eq.Gaussian_weight}
\end{equation}
where $\sigma_v = \sqrt{Var(\sum_{i = 1}^n\sum_{j = 1}^m c_{i,v}^{(j)}\epsilon_i^{(j)})}$.
\label{theorem.mu}
\end{theorem}
The existence of joint normal random variables $\xi_v$
satisfying the conditions in Theorem \ref{theorem.mu} is ensured because the right-hand side of eq.\eqref{eq.covariance_xi} constitutes the covariance matrix of the linear forms $\sum_{i = 1}^n\sum_{j = 1}^m c_{i,v}^{(j)}\epsilon_i^{(j)}$ in Example \ref{example.expand}, which is a symmetric positive semi-definite matrix.

The next result covers the consistency of Algorithm \ref{algorithm.bootstrap}. Define the following conditional probability and the conditional expectation:
$Prob^*\left(\cdot\right) = Prob\left(\cdot|X_i^{(j)}, i = 1,\cdots, n, j = 1,\cdots, m\right)$, and $\mathbf{E}^*\cdot = \left(\mathbf{E}\cdot|X_i^{(j)}, i = 1,\cdots, n, j = 1,\cdots, m\right)$, which are also referred to as `probability and expectation in the bootstrap world', like those introduced in \cite{MR1707286}. Notably, $\mathbf{E}^*\cdot = \mathbf{E}\cdot|\epsilon_i^{(j)}, i = 1,\cdots, n, j = 1,\cdots, m$. We propose our result in the following theorem.
\begin{theorem}
Suppose Assumptions 1 to 4 hold, then we have 
\begin{equation}
\begin{aligned}
\sup_{x\in\mathbf{R}}\vert
Prob^*\left(\max_{v = 1,\cdots, V}\frac{T_{n,m}}{B_{n,m}}\vert
\widehat{\mu}^*(x_v, y_v)
- \widehat{\mu}(x_v, y_v)
\vert\leq x\right)\\ - Prob\left(\max_{v = 1,\cdots, V}\vert\xi_v\vert\leq x\right)
\vert = o_p(1).
\end{aligned}
\label{eq.bootstrap_consistency}
\end{equation}
In addition, suppose Assumption 5 holds true, then 
\begin{equation}
\begin{aligned}
\sup_{x\in\mathbf{R}}\vert
Prob^*\left(\max_{v = 1,\cdots, V}\frac{T_{n,m}}{\widehat{\tau}_vB_{n,m}}\vert \widehat{\mu}^*(x_v, y_v) - \widehat{\mu}(x_v, y_v)
\vert\leq x\right)\\
- Prob\left(\max_{v = 1,\cdots, V}\vert\frac{\xi_v}{\tau_v}\vert\leq x\right)
\vert = o_p(1),
\end{aligned}
\label{eq.bootstrap_consistency_hetero}
\end{equation}
where $\tau_v=1$ in the homogeneous case and $\tau_v=\sigma^{1/3}_v$ in the heterogeneous case. 
\label{theorem.bootstrap}
\end{theorem}

\begin{remark}
According to Section 1.2 in \cite{MR1707286}, the consistency of the bootstrap algorithm is ensured if eq.\eqref{eq.bootstrap_consistency} holds.
\end{remark}

\section{Numerical experiments}
\label{section.numerical}
This section evaluates the proposed theoretical results and bootstrap algorithms through numerical experiments. In addition, it presents applications of our work in the analysis of a global climate dataset.

\textbf{Data-driven bandwidth selection: } In practice, bandwidths $\mathcal{K}$ and $\mathcal{B}$ need to be determined prior to estimation. Here we introduce data--driven methods to select the parameters. The selection of $\mathcal{K}$ is based on the cross-validation method in \cite{MR1168165}. Specifically, we pre--choose a largest bandwidth $\mathcal{K}_{\max}$. For each $\mathcal{K} = 1,2,\cdots, \mathcal{K}_{\max}$ and each $2\mathcal{K} + 1\leq i\leq n - 2\mathcal{K}$, $2\mathcal{K} + 1\leq j\leq m - 2\mathcal{K}$, calculate 
$$
\widetilde{\mu}(\frac{i}{n}, \frac{j}{m}) = \frac{1}{T_{n,m}}
\sum_{u,v = -\mathcal{K}, (u,v)\neq (0,0)}^\mathcal{K}
X_{u + i}^{(v + j)}
G\left(\frac{u}{\mathcal{K}}\right) G\left(\frac{v}{\mathcal{K}}\right).
$$
The optimal bandwidth is chosen to minimize the mean square error 
\begin{align*}
    \widetilde{E} = \sum_{i = 2\mathcal{K} + 1}^{n - 2\mathcal{K}}\sum_{j = 2\mathcal{K} + 1}^{m - 2\mathcal{K}} (X_i^{(j)} - \widetilde{\mu}(\frac{i}{n}, \frac{j}{m}))^2.
\end{align*}

Selecting the variance bandwidth $\mathcal{B}$ can be more involved. There has been discussions regarding this problem in the literature. For instance, \cite{MR2041534} developed a bandwidth selection procedure for time series analysis;  while \cite{MR1366282} and \cite{MR2460005} proposed an empirical method that is more suitable for our setting. We also recommend the book by \cite{MR2001447} for more discussion. Algorithm \ref{algorithm.selection} follows the approach of \cite{MR1366282} to select $\mathcal{B}$, but it incorporates randomization instead of considering all small blocks to save computation time.

\begin{algorithm}[Selection of variance bandwidth $\mathcal{B}$]
    \label{algorithm.selection}
    $\\$
\noindent
\textbf{Input: } The observed data $X_i^{(j)}, i = 1,\cdots, n, j = 1,\cdots, m$, the kernel function $G(\cdot)$, the variance kernel function $K(\cdot)$, the bandwidth $\mathcal{K}$, the proportion of small blocks $q\in (0,1)$, the potential variance kernel set $\Gamma$, the number of iterations $H$, the number of bootstrap replicates $B$. A plausible pilot variance bandwidth $\mathcal{B}_p$. Our work selects $\mathcal{B}_p = 5$, the set $\Gamma = \{1,2,\cdots, 10\}$, $q = 0.1$, and $H = 15$.

1. Run Steps 1 to 3 in Algorithm \ref{algorithm.bootstrap} using variance bandwidth $\mathcal{B}_p$, generate the residuals $\epsilon_i^{(j)*}$ for $i = \mathcal{K} + 1,\cdots, n - \mathcal{K}$ and $j = \mathcal{K} + 1,\cdots, m - \mathcal{K}$, then calculate the sample mean 
$
\gamma^*_b = \frac{1}{(n - 2\mathcal{K})\times (m - 2\mathcal{K})}\sum_{i = \mathcal{K} + 1}^{n - \mathcal{K}}\sum_{j = \mathcal{K} + 1}^{n - \mathcal{K}}\epsilon_i^{(j)*}.
$

2. Repeat Step 1 for $b = 1,\cdots, B$, then calculate the sample variance 
$\widehat{\sigma}^{2*} = \frac{1}{B}\sum_{b = 1}^B (\gamma^*_b - \frac{1}{B}\sum_{b = 1}^B \gamma^*_b)^2$.

3. Uniformly draw a location $u\in\{1,2,\cdots, n - \lfloor qn\rfloor\}$ and $v\in\{1,2,\cdots, m - \lfloor qm\rfloor\}$. After that, choose the small block to be $X_{i}^{(j)}$, where $i = u, u + 1,\cdots, u + \lfloor qn\rfloor$, $j = v, v + 1,\cdots, v + \lfloor qm\rfloor$. Denote $n_0 = \lfloor qn\rfloor + 1$ and $m_0 = \lfloor qm\rfloor + 1$.

4. Run Steps 1 to 3 in Algorithm \ref{algorithm.bootstrap} for all $\mathcal{B}\in\Gamma$ on the small block, generate the residuals $\xi_{i}^{(j)*}$ for $i = \mathcal{K} + 1,\cdots, n_0 - \mathcal{K}$ and $j = \mathcal{K} + 1,\cdots, m_0 - \mathcal{K}$, then calculate the sample mean 
$
\eta^*_b = \frac{1}{(n_0 - 2\mathcal{K})\times (m_0 - 2\mathcal{K})}\sum_{i = \mathcal{K} + 1}^{n_0 - \mathcal{K}}\sum_{j = \mathcal{K} + 1}^{m_0 - \mathcal{K}}\xi_i^{(j)*}.
$ 

5. Repeat Step 4 for $b = 1,\cdots, B$, then calculate the sample variance 
$
\widehat{\tau}^{2*}_l = \frac{1}{B}\sum_{b = 1}^B (\eta^*_b - \frac{1}{B}\sum_{b = 1}^B \eta^*_b)^2.
$

6. Repeat Steps 3 to 5 for $l = 1,\cdots, H$. Choose $\mathcal{B}$ that minimize the loss 
$
\sum_{l = 1}^H (\widehat{\tau}^{2*}_l - \widehat{\sigma}^{2*})^2.
$
\end{algorithm}

\begin{figure}[tbp]
\centering
    \subfigure[Standard deviation of AR noise]{
    \includegraphics[width = 2.1in]{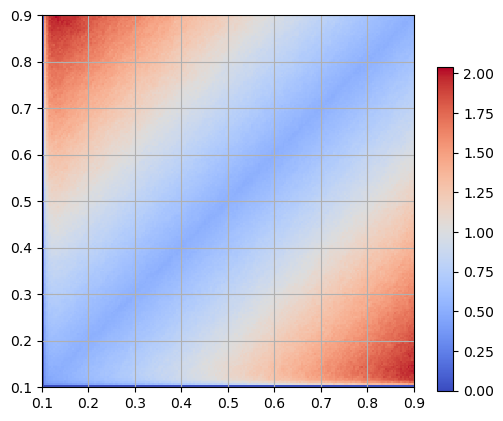}
    \label{oneaX}
  }
    \subfigure[Standard deviation of MA noise]{
    \includegraphics[width = 2.1in]{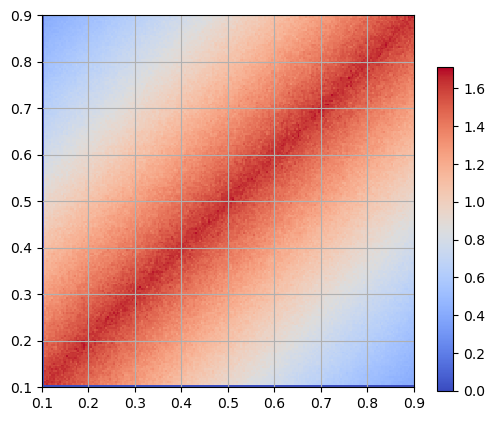}
    \label{oneb}
  }
  \caption{Element--wise standard deviations for AR and MA noises.} 
    \label{datalike}
\end{figure}

\begin{figure}[htbp]
\centering
\subfigure[Elliptical mean field]{
    \includegraphics[width = 1.8in]{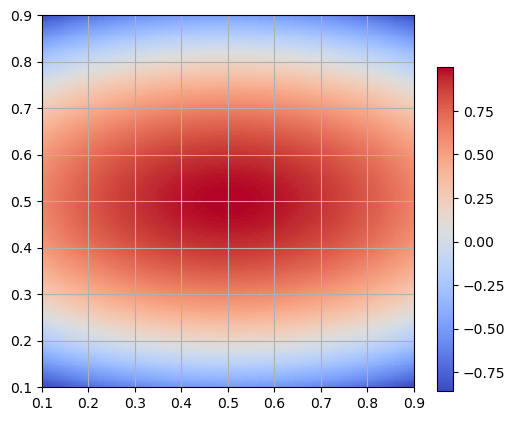}
    \label{oneaX}
}
\subfigure[Elliptical observations]{
    \includegraphics[width = 1.8in]{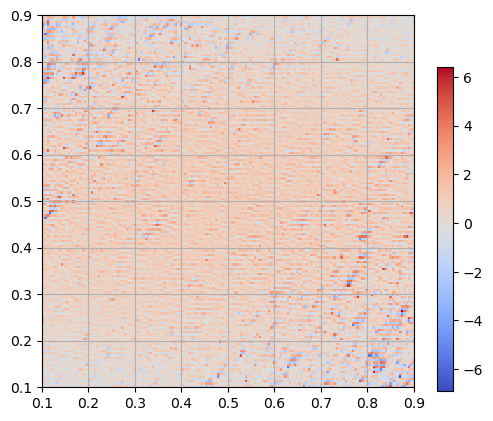}
    \label{oneaX}
}
\subfigure[Elliptical estimated mean]{
    \includegraphics[width = 1.8in]{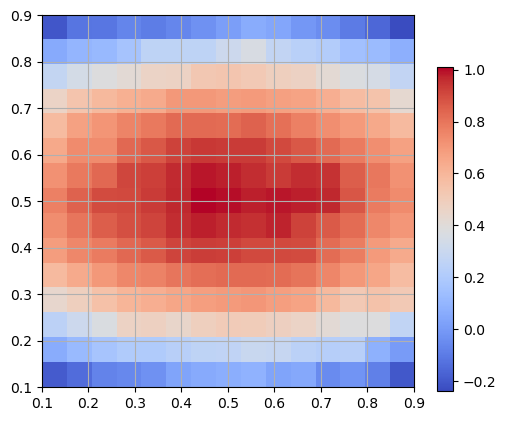}
    \label{oneaX}
}
\subfigure[Sinusoidal mean field]{
    \includegraphics[width = 1.8in]{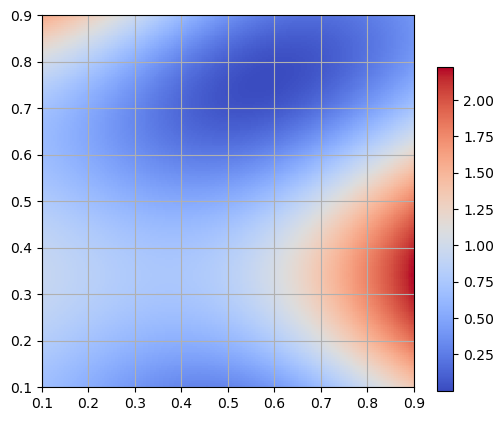}
    \label{oneaX}
}
\subfigure[Sinusoidal observations]{
    \includegraphics[width = 1.8in]{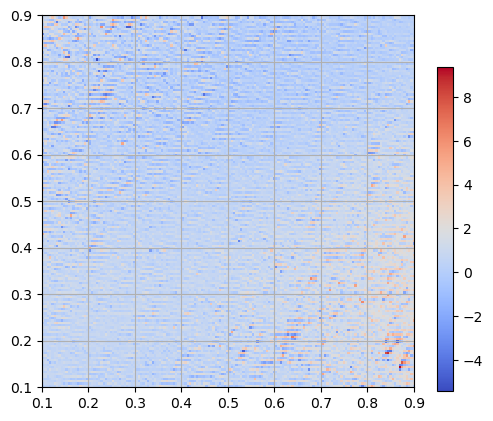}
    \label{oneaX}
}
\subfigure[Sinusoidal estimated  mean]{
    \includegraphics[width = 1.8in]{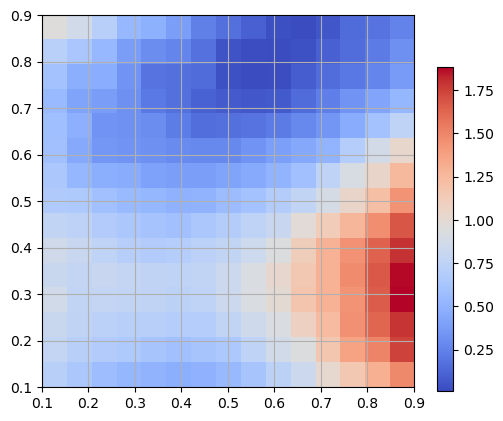}
    \label{oneaX}
}
\caption{Observations and the estimators of the mean field. The grid size is $15\times 15$. }
\label{figure.estimation_res}
\end{figure}

\textbf{Simulated data: } We generate data $X_i^{(j)} = \mu\left(\frac{i}{n}, \frac{j}{m}\right) + \epsilon_i^{(j)}, i = 1,\cdots, n, j = 1,\cdots, m$ with $n = m = 200$. We choose the following mean fields
\begin{align*}
\text{Elliptical: } \mu(x, y) = 1.0 - (\frac{3}{2}(x - 0.5)^2 + 6(y - 0.5)^2),\\
\text{Sinusoidal: } \mu(x, y) = \sin(2(x - 0.6))^2 + \cos(3(y - 0.3))^2 + \sin(2(x - 0.6))\times \cos(3(y - 0.3)),
\end{align*}
and simulate the random errors $\epsilon_i^{(j)}$ from a 2-D autoregressive(AR) model in \cite{app13158877} or a moving average(MA) process, as follows:
\begin{align*}
    \text{AR: } \epsilon_i^{(j)} = 0.3\epsilon_{i - 1}^{(j)} - 0.4\epsilon_i^{(j - 1)} - 0.2\epsilon_{i - 1}^{(j - 1)} + e_i^{(j)},\\
    \text{MA: } \epsilon_i^{(j)} = 0.3f_{i - 1}^{(j)} - 0.4f_i^{(j - 1)} - 0.2f_{i - 1}^{(j - 1)} + f_i^{(j)},
\end{align*}
where $e_{i}^{(j)} = E_{i,1}^{(j)}\times E_{i,2}^{(j)}$, $f_i^{(j)} = F_{i,1}^{(j)}\times F_{i,2}^{(j)}$, and $E_{i,1}^{(j)}$, $E_{i,2}^{(j)}$, $F_{i,1}^{(j)}$, $F_{i,2}^{(j)}$ are mutually independent normal random variables with mean $0$, and 
\begin{align*}
    Var(E_{i,1}^{(j)}) = (0.7 + 0.5\times \vert \frac{i}{n} - \frac{j}{m}\vert)^2,\ Var(E_{i,2}^{(j)}) = (0.5 + 0.7\times \vert \frac{i}{n} - \frac{j}{m}\vert)^2,\\
    Var(F_{i,1}^{(j)}) = (1.2 - 0.5\times \vert \frac{i}{n} - \frac{j}{m}\vert)^2,\ Var(E_{i,2}^{(j)}) = (1.2 - 0.7\times \vert \frac{i}{n} - \frac{j}{m}\vert)^2.
\end{align*}
This construction establishes a correlated and heterogeneous non-Gaussian random field. Figure \ref{datalike} plots the element--wise standard deviations of the noises $\epsilon_i^{(j)}$, illustrating their heterogeneous nature; Figure \ref{figure.estimation_res} illustrates both the efficacy and imperative nature of employing the non-parametric estimator eq.\eqref{eq.define_Nadaraya}  for analyzing a random field: The mean information within the random field is obscured by noises, making it hard for practitioners to discover valuable insights only through direct inspection without resorting to smoothing techniques; 
Table \ref{table.estimation} records the performance of the homogeneous and heterogeneous versions of locally weighted multiplier bootstrap algorithm. It turns out that both versions are able to construct simultaneous confidence regions with desired coverage probabilities, even when the number of positions is large.

\begin{table}[htbp]
\centering
\caption{Simulated coverage probabilities of bootstrap algorithm \ref{algorithm.bootstrap}. The nominal coverage probability is $95\%$. $\mathcal{K}$ and $\mathcal{B}$ are selected via the aforementioned methods. The simulation repeats 200 times, and the number of bootstrap replicates is 200.}
\scriptsize
\begin{tabular}{l l l l l l l l l l l}
\hline\hline
              &               &                    & &      & \multicolumn{2}{l}{Homogeneous} & \multicolumn{2}{l}{Heteroeneous}\\
\hline
Mean &  Error & $\mathcal{K}$ & $\mathcal{B}$ & Grid                     & Coverage & Average width & Coverage & Average width\\
Elliptical  & AR & 10 & 2 &               $20\times 20$            & 93.0\% & 0.266 & 96.5\% & 0.336\\
            & &    &   &               $40\times 40$            & 96.0\% & 0.281 & 94.0\% & 0.365\\
            & &   &    &               $60\times 60$            & 95.5\% & 0.292 & 92.5\% & 0.381\\
\hline 
Sinusoidal  & AR  & 10            & 3             & $20 \times 20$           & 96.5\% & 0.255 & 94.5\% & 0.275\\
            &   & &               & $40 \times 40$           & 95.0\% & 0.274 & 95.5\% & 0.278\\
            &  & &               & $60 \times 60$           & 94.5\% & 0.285 & 96.0\% & 0.309\\
\hline
Elliptical  & MA & 10 & 3 &      $20\times 20$  & $98.5\%$ & 0.365 & 91.5\% &
0.416\\
            &    &    &   &      $40\times 40$  & 98.0\% & 0.388 & 94.5\% & 0.453\\
            &    &    &   &      $60\times 60$  & 97.5\% & 0.400 & 91.5\% & 0.467\\
\hline
Sinusoidal & MA & 10 & 4    & $20\times 20$ & 94.0\% & 0.344 & 94.5\% & 0.451\\
           &    &    &      & $40\times 40$ & 92.5\% & 0.367 & 94.5\% & 0.490\\
           &    &     &     & $60\times 60$ & 97.0\% & 0.383 & 93.5\% & 0.503\\
\hline\hline
\end{tabular}
\label{table.estimation}
\end{table}


\textbf{Testing for signals in a random field: } \cite{MR1951265} proposed a wavelet--based method for detecting signals in a random field with stationary Gaussian noises, which involved testing the null hypothesis $\mu(x,y) = 0$ for all possible $x,y$.
The procedure is implemented in the R--package ``EFDR''\footnote{see \textit{https://github.com/andrewzm/EFDR/.}}. This part demonstrates our method's size and power in this hypothesis testing problem in comparison to the EFDR method. Since $n$ and $m$ have to be powers of $2$ in EFDR, we choose $n = m = 128$. Our work performs hypothesis testing on a $15\times 15$ sub--grid of the random field. In the alternative situation, following \cite{MR1951265}, we consider the alternative mean field
\begin{align*}
    \mu(x, y) = 
    \begin{cases}
    0.3,\ \text{if} (x - 0.5)^2 + (y - 0.5)^2\leq 0.1^2,\\
    0,\ \text{otherwise}.
    \end{cases}
\end{align*}
Since the method proposed by \cite{MR1951265} is tailored for stationary Gaussian errors, we also simulated the situation where $\epsilon_i^{(j)}$ are i.i.d. standard normal. The hypothesis testing result is demonstrated in Table \ref{table.test_mean_field}. From the simulation results we see that when $\epsilon_i^{(j)}$ are i.i.d. Gaussian, the method proposed by \cite{MR1951265} achieves  the desirable size, but is less powerful compared to the locally weighted multiplier bootstrap. When the assumptions of normality and stationarity are violated, their method tends to over-reject $H_0$ under our simulation scenarios, demonstrating that their method does not carry over to the non-stationary case. On the other hand, our methods perform reasonably well across all simulation scenarios.

\begin{table}[htbp]
    \centering
    \caption{Simulated size and power for signal detection. The nominal size is $5\%$. ``FDR'' stands for the test function test.fdr implemented in ``EFDR'' package, and the hyperparameters in test.fdr function are chosen by default. The number of bootstrap replicates and the number of simulations coincide with Table \ref{table.estimation}.}
    \begin{tabular}{l| l l l | l l  l}
    \hline\hline 
    &  \multicolumn{3}{c}{Size under $H_0$}     &  \multicolumn{3}{c}{Power Under $H_1$}\\
    \hline 
Noise  &  Homo & Hetero & FDR       &  Homo & Hetero & FDR\\
\hline 
Normal &   9.0\%     & 0.5\% & 4.5\% & 92.5\% & 96.5\% & 8.0\%\\ 
AR     &   4.5\%    & 2.0\% & 100\%     & 100\%  &         100\%        & 100\%\\
MA     &   2.5\%    & 1.0\%       & 100\%     & 95.0\% &     64.0\%           & 100\%\\
    \hline\hline
    \end{tabular}
    \label{table.test_mean_field}
\end{table}

\begin{figure}[htbp]
\centering
    \subfigure[The difference random field]{
    \includegraphics[width = 1.8in]{real_dataX.png}
    \label{oneaX}
  } 
    \subfigure[Element--wise variances]{
    \includegraphics[width = 2.0in]{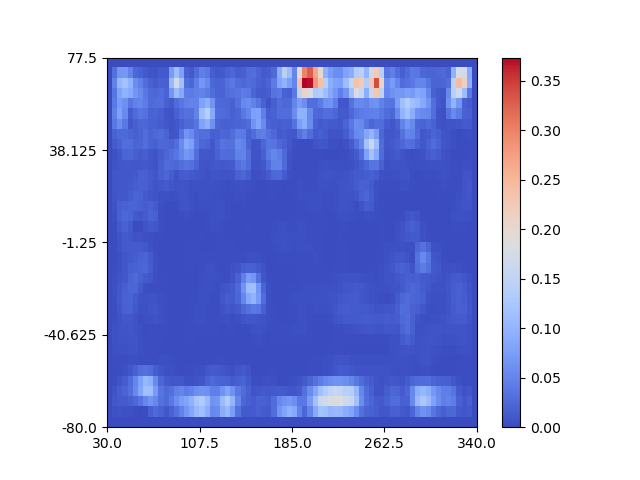}
  }
    \subfigure[Element--wise column--wise covariances]{
    \includegraphics[width = 2.0in]{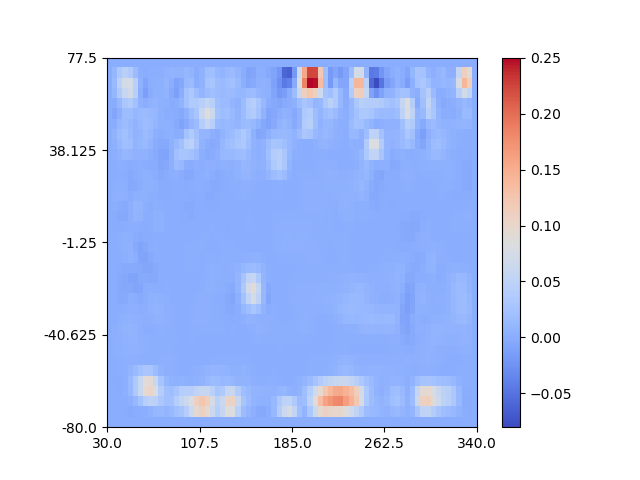}
  }
    \subfigure[Simultaneous confidence region]{
        \includegraphics[width = 2.45in]{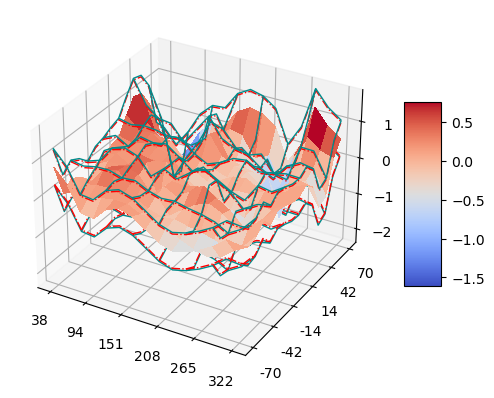}
    }
    \subfigure[Estimated mean temperature anomalies and testing result]{
    \includegraphics[width = 1.8in]{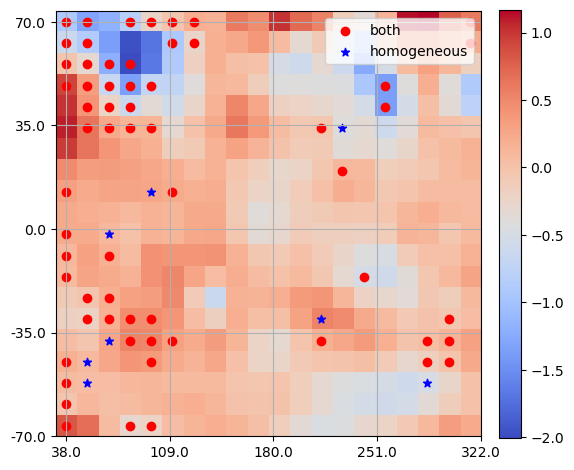}
    \label{oneZZZ}
}   
    \subfigure[Signal detected by EFDR package]{
    \includegraphics[width = 1.7in]{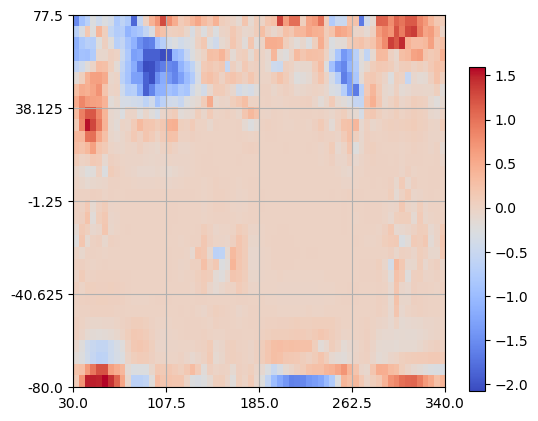}
    \label{oneb}
}
\label{figure.real_temp}
\caption{The difference in temperature anomalies between the years 2008–2010 and 2004–2006, their estimated variances, estimated column--wise covariances, the estimated mean temperature, and the simultaneous confidence region. Red dots in Figure \ref{oneZZZ} stand for positions where the simultaneous confidence regions generated by both homogeneous and heterogeneous methods do not cover 0, while blue dots indicate positions where the homogeneous version of Algorithm \ref{algorithm.bootstrap} do not cover 0. For the EFDR package, non--zero signal is considered to be statistically significant. }
\end{figure}

\textbf{Analysis of global surface temperature. } The NOAA Global Surface Temperature dataset collected by \cite{data_set} records monthly gridded global surface temperature anomalies (departure from the long-term average) from January 1850 to April 2024. Each cell in the dataset represents the average temperature of a $5^{\circ}$ latitude $\times 5^{\circ}$ longitude region, resulting in a grid of  $72 \times 36$ cells. We compare the mean temperature dynamics between the years 2008–2010 and 2004–2006. To achieve the goal, we calculate the 2--year average temperature of two time periods respectively, and subtract the latter from the former, resulting in the difference field. We use Algorithm \ref{algorithm.bootstrap} and the package``EFDR'' to test the null hypothesis $H_0:$ the difference field has 0 mean, and the result is demonstrated in Figure \ref{figure.real_temp}. The variances and column--wise covariances of the difference field vary significantly across positions, highlighting the non--stationary nature of the data. In addition, a normal test (using the `\textit{normaltest}' function from Python's  \textit{scipy} package) applied to the difference field yields a P--value of $2.14\times 10^{-83}$, so the difference field  is highly unlikely to be Gaussian. Therefore, classic results based on stationary or Gaussian assumptions may not be suitable for analyzing this data set. Algorithm \ref{algorithm.bootstrap} confirms a temperature increase in Africa and a temperature decrease in Europe. According to \cite{https://doi.org/10.1002/joc.7295}, a surge of heat waves occurred around 2010 in Africa. On the other hand,  \cite{https://doi.org/10.1029/2010GL044613} mentioned an exceptional Northern Hemisphere mean atmospheric circulation during 2010. These two pieces of evidence support our findings. On the other hand, the testing results of \cite{MR1951265} suggest similar findings but are less satisfactory because they do not detect temperature changes in Africa, possibly due to the non-stationary and non-Gaussian nature of the random field.

\begin{figure}[htbp]
\centering
    \subfigure[The fitted residuals]{
    \includegraphics[width = 1.9in]{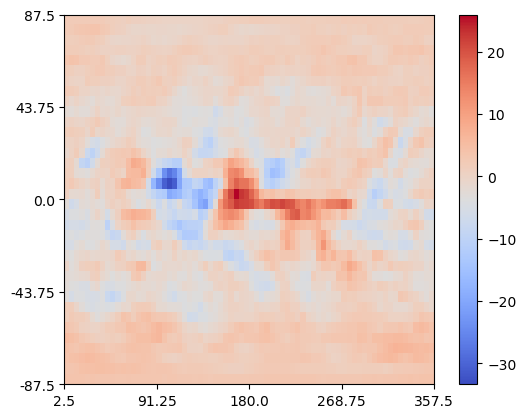}
  }
    \subfigure[Estimation and testing result by Algorithm \ref{algorithm.bootstrap}]{
    \includegraphics[width = 1.9in]{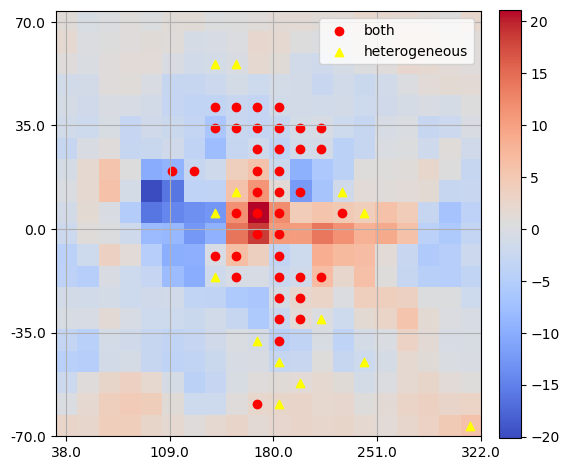}
  } 
    \subfigure[Signal detected by ``EFDR'']{
    \includegraphics[width = 2.0in]{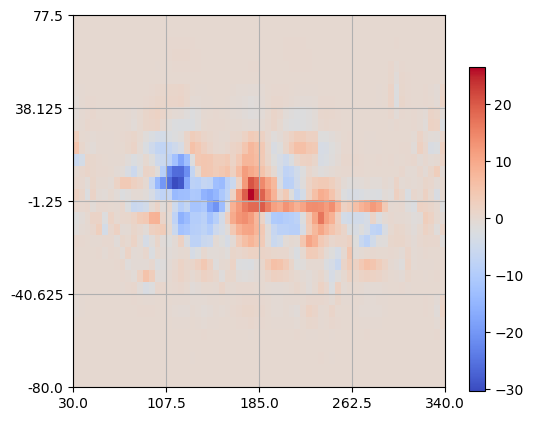}
  }
  \caption{The fitted residuals and their estimated mean for Earth's OLR data. Red dots stand for positions where
the simultaneous confidence regions generated by both homogeneous and heterogeneous methods do not cover 0, and
yellow triangles indicate positions where the simultaneous confidence regions generated by the
heterogeneous version of Algorithm 1 do not cover 0. For the EFDR package, non--zero signal is considered to be statistically significant.}
  \label{figure.test_result2}
\end{figure}

\textbf{Testing for the linearity of Earth’s OLR:} \cite{doi:10.1073/pnas.1809868115} studied the relationship between Earth's outgoing longwave radiation(Earth's OLR) and surface temperature. They found that the relationship remained linear within the temperature range of 220--280K (i.e., $-53.15$--$6.85^\circ C$), contrary to the Stefan--Boltzmann law. We aim to provide statistical justification for this finding.  To achieve the goal, we calculate the difference in the yearly mean temperature anomaly between 2015 and 2016 using the data set proposed by \cite{data_set}, then perform the same calculations for OLR using the dataset in \cite{70af5bdd-0375-35ff-b66e-5a4a3b2f11e2} \footnote{Also see \text{https://psl.noaa.gov/data/gridded/data.olrcdr.interp.html}}. After that, we fit a linear model and derive the fitted residuals. The hypothesis  test we are interested in is $H_0:$ the fitted residuals have mean 0. We adopt Algorithm \ref{algorithm.bootstrap} and the method proposed by \cite{MR1951265} to achieve the goal. Both results, shown in Figure \ref{figure.test_result2}, suggest that the hypothesis tests fail to reject $H_0$ in temperate and frigid zones, where the average temperature is relatively low; and reject $H_0$ in tropical and subtropical areas, where the average temperature is high. These test results support the finding of  \cite{doi:10.1073/pnas.1809868115} that OLR is approximately a linear function of surface temperature for low temperatures(below 280K), but this linearity breaks down at higher temperatures. Compared to \cite{MR1951265}, our work  rejects more positions around the 35th parallel north, which roughly marks the boundary of the subtropical zones. Since subtropical zones are characterized by hot summers and mild winters, this difference seems reasonable.

\appendix
\numberwithin{equation}{section}
\numberwithin{lemma}{section}
\numberwithin{corollary}{section}

\section{Preliminarily}
This section introduces the properties of some special functions that are frequently used in the following proofs.
For a function $f(x_1,\cdots, x_n)$, define the operator $\partial_i f = \frac{\partial f}{\partial x_i}$. For given $\tau,\psi>0$ and $z\in\mathbf{R}$, define the function
$$
G_\tau(x_1,\cdots, x_v) = \frac{1}{\tau}\log\left(\sum_{i = 1}^v\exp(\tau x_i) + \sum_{i = 1}^v \exp(-\tau x_i)\right)
$$
and the function $g_0(x) = (1 - \min(1, \max(x, 0))^4)^4$, $g_{\psi, z}(x) = g_0(\psi(x-  z))$. Define
$h_{\tau, \psi, z}(x_1,\cdots, x_v) = g_{\psi, z}(G_\tau(x_1,\cdots, x_v))$, then from \cite{MR3161448}, \cite{MR3992401} and \cite{MR4441125},
\begin{align*}
\sum_{i = 1}^v \vert\partial_i h_{\tau, \psi, z}(x_1,\cdots, x_v)\vert\leq C\psi\\
\sum_{i = 1}^v\sum_{j = 1}^v\vert\partial_i\partial_j h_{\tau, \psi, z}(x_1,\cdots, x_v)\vert\leq C\psi^2 + C\psi\tau\\
\sum_{i = 1}^v\sum_{j = 1}^v\sum_{k = 1}^v \vert \partial_i\partial_j\partial_k h_{\tau, \psi, z}(x_1,\cdots, x_v)\vert
\leq C(\psi^3 + \psi^2\tau + \psi\tau^2)
\end{align*}
Notably, here $C$ is independent of $x$ and $v$, and the above formula is valid for all positive integer $v$. Another frequently used result is the joint normal random variables' anti-concentration property proposed by \cite{MR3350040}. We quote this result in the following lemma.
\begin{lemma}
(i). Suppose $\xi_1,\cdots, \xi_V$ are joint normal random variables with $\mathbf{E}\xi_i = 0$ and $\mathbf{E}\xi_i\xi_j = \sigma_{ij}$. Suppose $\exists$ two constants $0<c_0\leq C_0<\infty$ such that $c\leq \sigma_{ii}\leq C$ for $i = 1,\cdots, V$. Then for $\forall \delta> 0$,
\begin{equation}
\begin{aligned}
\sup_{x\in\mathbf{R}}\vert Prob\left(\max_{i = 1,\cdots, V}\vert\xi_i\vert\leq x+\delta\right) - Prob\left(\max_{i = 1,\cdots, V}\vert\xi_i\vert\leq x\right)\vert\\
\leq C^\prime\delta(1 + \sqrt{\log(V)} + \sqrt{\vert\log(\delta)\vert})
\end{aligned}
\label{eq.first_normal}
\end{equation}
Here $C^\prime$ only depends on $c_0$ and $C_0$.

\noindent (ii). In addition suppose $\xi_1^\dagger,\cdots, \xi_V^\dagger$ are joint normal random variables with $\mathbf{E}\xi_i^\dagger = 0$ and
$\mathbf{E}\xi_i^\dagger\xi_j^\dagger = \sigma_{ij}^\dagger$, define $\Delta = \max_{i,j = 1,\cdots, V}\vert\sigma_{ij} - \sigma_{ij}^\dagger\vert$.
Suppose $0<\Delta<1$, then
\begin{equation}
\begin{aligned}
\sup_{x\in\mathbf{R}}\vert Prob\left(\max_{i = 1,\cdots, V}\vert\xi_i\vert\leq x\right) -
Prob\left(\max_{i = 1,\cdots, V}\vert\xi_i^\dagger\vert\leq x\right)\vert\\
\leq C^{\prime\prime}\Delta^{1/3}\times (1 + \sqrt{\log(V)} + \vert\log(\Delta)\vert)^2
\end{aligned}
\label{eq.second_normal}
\end{equation}
here $C^{\prime\prime}$ is a constant only depending on $c_0$ and $C_0$.
\label{lemma.Gaussian_property}
\end{lemma}
\begin{proof}
Notice that $\max_{i = 1,\cdots, V}\vert\xi_i\vert = \max(\max_{i = 1,\cdots, V}\xi_i, \max_{i = 1,\cdots, V}(-\xi_i))$, so we have
\begin{align*}
Prob\left(\max_{i = 1,\cdots, V}\vert\xi_i\vert\leq x+\delta\right) - Prob\left(\max_{i = 1,\cdots, V}\vert\xi_i\vert\leq x\right)\\
= Prob\left(0<\max(\max_{i = 1,\cdots, V}\xi_i, \max_{i = 1,\cdots, V}(-\xi_i)) - x\leq \delta\right)\\
\leq Prob\left(\vert \max_{i = 1,\cdots, V}\xi_i - x\vert\leq \delta\right)  + Prob\left(\vert\max_{i = 1,\cdots, V}(-\xi_i) - x\vert\leq x\right)
\end{align*}
From (A.8) in \cite{MR4441125} and theorem 3, eq.(18), eq.(19) in \cite{MR3350040}, we prove \eqref{eq.first_normal}.

Without loss of generality, assume $\xi_i$ and $\xi_j^\dagger$ are independent for any $i,j$. Similar to \cite{MR3350040}, define $Z_i(t) = \sqrt{t}\xi_i  + \sqrt{1 - t}\xi_i^\dagger$, then for any $\tau,\psi>0$ and $z\in\mathbf{R}$,
\begin{align*}
\mathbf{E}(h_{\tau,\psi,z}(\xi_1,\cdots, \xi_V) - h_{\tau,\psi,z}(\xi_1^\dagger,\cdots, \xi_V^\dagger))
= \frac{1}{2}\sum_{i = 1}^V\int_{[0,1]} t^{-1/2}\mathbf{E}\xi_i\partial_i h_{\tau,\psi,z}(Z_1(t),\cdots, Z_V(t))\mathrm{d}t\\
- \frac{1}{2}\sum_{i = 1}^V\int_{[0,1]} (1 - t)^{-1/2}\mathbf{E}\xi_i^\dagger\partial_i h_{\tau,\psi,z}(Z_1(t),\cdots, Z_V(t))\\
= \frac{1}{2}\sum_{i = 1}^V\sum_{j = 1}^V (\sigma_{ij} - \sigma_{ij}^\dagger)\int_{[0,1]} \mathbf{E}\partial_i\partial_j h_{\tau,\psi,z}(Z_1(t),\cdots, Z_V(t))\mathrm{d}t
\end{align*}
therefore
\begin{equation}
\begin{aligned}
\sup_{z\in\mathbf{R}}\vert\mathbf{E}(h_{\tau,\psi,z}(\xi_1,\cdots, \xi_V) - h_{\tau,\psi,z}(\xi_1^\dagger,\cdots, \xi_V^\dagger))\vert\leq C\Delta(\psi^2 + \psi\tau)
\end{aligned}
\label{eq.delta_h}
\end{equation}
In particular, set $t = \frac{1}{\psi} + \frac{\log(2V)}{\tau}$, from eq.(A.17) and eq.(A.18) in \cite{MR4441125}
\begin{equation}
\begin{aligned}
Prob\left(\max_{i = 1,\cdots, V}\vert\xi_i\vert\leq x\right) -
Prob\left(\max_{i = 1,\cdots, V}\vert\xi_i^\dagger\vert\leq x\right)
\leq C^\prime t(1 + \sqrt{\log(V)} + \sqrt{\vert\log(t)\vert})\\
 + \mathbf{E}h_{\tau,\psi,x-\frac{1}{\psi}}(\xi_1,\cdots, \xi_V)
-\mathbf{E}h_{\tau,\psi,x-\frac{1}{\psi}}(\xi_1^\dagger,\cdots, \xi_V^\dagger)\\
\text{and } Prob\left(\max_{i = 1,\cdots, V}\vert\xi_i\vert\leq x\right) -
Prob\left(\max_{i = 1,\cdots, V}\vert\xi_i^\dagger\vert\leq x\right)\\
\geq -C^\prime t(1 + \sqrt{\log(V)} + \sqrt{\vert\log(t)\vert}) + \mathbf{E}h_{\tau,\psi,x+\frac{\log(2V)}{\tau}}(\xi_1,\cdots, \xi_V)\\
- \mathbf{E}h_{\tau,\psi,x+\frac{\log(2V)}{\tau}}(\xi_1,\cdots, \xi_V)
\end{aligned}
\label{eq.prob_to_h}
\end{equation}
which implies
\begin{equation}
\begin{aligned}
\sup_{x\in\mathbf{R}}\vert Prob\left(\max_{i = 1,\cdots, V}\vert\xi_i\vert\leq x\right) -
Prob\left(\max_{i = 1,\cdots, V}\vert\xi_i^\dagger\vert\leq x\right)\vert\\
\leq C^\prime t(1 + \sqrt{\log(V)} + \sqrt{\vert\log(t)\vert}) + C\Delta(\psi^2 + \psi\tau)
\end{aligned}
\end{equation}
Choose $\tau = \psi = \frac{1}{\Delta^{1/3}}\times (1 + \sqrt{\log(V)} + \vert\log(\Delta)\vert)$, then
we prove eq.\eqref{eq.second_normal}.
\end{proof}

\section{Proof of theorems in section \ref{section.statistical_inference}}

We first introduce  the notations for several sub--$\sigma$--fields of the random field, which  will be frequently used in the proofs throughout this and the following sections. For any integers $a\leq b, c\leq d\in\mathbf{Z}$, define $\mathcal{Z}_{(a,b)}^{(c,d)}$ as the $\sigma-$field generated by
$$
\left\{e_x^{(y)}: a\leq x\leq b\ \text{and }  c\leq y\leq d\right\}.
$$
Define $\mathcal{Z} = \mathcal{Z}_{(-\infty,\infty)}^{(-\infty,\infty)}$, indicating the $\sigma-$field generated by all $e_x^{(y)}$. For any given integers $i,j\in\mathbf{Z}$ and $a,b,c,d\geq 0$, define
$\mathcal{H}_{i, (a,b)}^{(j), (c,d)} = \mathcal{Z}_{(i - a, i + b)}^{(j - c, j + d)}$.

We begin the proofs with several lemmas.

\begin{lemma}
Suppose two sets $\mathcal{C}_1,\mathcal{C}_2\subset \{(i,j): i,j\in\mathbf{Z}\}$ such that
$\vert\mathcal{C}_1\vert, \vert\mathcal{C}_2\vert < \infty$ and $\mathcal{C}_1\cap \mathcal{C}_2\neq\emptyset$. Suppose a random variable $X$ is measurable in $\mathcal{Z}$ (see section \ref{section.statistical_inference}) and $\mathbf{E}\vert X\vert<\infty$, then
\begin{equation}
\mathbf{E}\left(\mathbf{E}X|\{e_i^{(j)}:(i,j)\in\mathcal{C}_1\}\right)|\{e_i^{(j)}:(i,j)\in\mathcal{C}_2\}
= \mathbf{E}X|\{e_i^{(j)}: (i,j)\in \mathcal{C}_1\cap\mathcal{C}_2\}
\label{eq.conditional_prob}
\end{equation}
almost surely.
\label{lemma.conditional_prob}
\end{lemma}
\begin{proof}
Define $\mathcal{D}_s, s = 1,2$ as the $\sigma-$field generated by $e_i^{(j)}, (i,j)\in\mathcal{C}_s$, then
$\mathbf{E}Y|\{e_i^{(j)}:(i,j)\in\mathcal{C}_s\} = \mathbf{E}Y|\mathcal{D}_s$ for any random variable $Y$. Define the $\lambda-$system
$$
\left\{A\in\mathcal{D}_2: \mathbf{E}(\mathbf{E}X|\mathcal{D}_1)\times\mathbf{1}_A = \mathbf{E}\left(\mathbf{E}X|\{e_i^{(j)}: (i,j)\in \mathcal{C}_1\cap\mathcal{C}_2\}\right)\times \mathbf{1}_A\right\}
$$
and the $\pi-$system $\prod_{(i,j)\in\mathcal{C}_2}A_i^{(j)}$, here $A_i^{(j)}$ is generate by $e_i^{(j)}$. Since
\begin{align*}
\mathbf{E}\left(\mathbf{E}X|\{e_i^{(j)}: (i,j)\in \mathcal{C}_1\cap\mathcal{C}_2\}\right)\times \prod_{(i,j)\in\mathcal{C}_2}\mathbf{1}_{e_i^{(j)}\in A_i^{(j)}}\\
=\prod_{(i,j)\in(\mathcal{C}_2 - \mathcal{C}_1)}Prob\left(e_i^{(j)}\in A_i^{(j)}\right)\times
\mathbf{E}\left[\left(\mathbf{E}X|\{e_i^{(j)}: (i,j)\in \mathcal{C}_1\cap\mathcal{C}_2\}\right)\times \prod_{(i,j)\in\mathcal{C}_1\cap\mathcal{C}_2}\mathbf{1}_{e_i^{(j)}\in A_i^{(j)}}\right]\\
=\prod_{(i,j)\in(\mathcal{C}_2 - \mathcal{C}_1)}Prob\left(e_i^{(j)}\in A_i^{(j)}\right)\times\mathbf{E}\left[X\times \prod_{(i,j)\in\mathcal{C}_1\cap\mathcal{C}_2}\mathbf{1}_{e_i^{(j)}\in A_i^{(j)}}\right]\\
= \prod_{(i,j)\in(\mathcal{C}_2 - \mathcal{C}_1)}Prob\left(e_i^{(j)}\in A_i^{(j)}\right)\times\mathbf{E}\left[(\mathbf{E}X|\mathcal{D}_1)\times \prod_{(i,j)\in\mathcal{C}_1\cap\mathcal{C}_2}\mathbf{1}_{e_i^{(j)}\in A_i^{(j)}}\right]\\
= \mathbf{E}\left(\mathbf{E}X|\mathcal{D}_1\right)\times \prod_{(i,j)\in\mathcal{C}_2}\mathbf{1}_{e_i^{(j)}\in A_i^{(j)}}
\end{align*}
From $\pi-\lambda$ theorem, we know that
$$
\mathbf{E}(\mathbf{E}X|\mathcal{D}_1)|\mathcal{D}_2 = \mathbf{E}\left(\mathbf{E}X|\{e_i^{(j)}: (i,j)\in \mathcal{C}_1\cap\mathcal{C}_2\}\right)
$$
almost surely, and we prove \eqref{eq.conditional_prob}.
\end{proof}

\begin{lemma}
Suppose $\epsilon_i^{(j)}, i = 1,\cdots, n, j = 1,\cdots, m$ are $(M,\alpha)-$short range dependent random variables with $M \geq 1$ and $\alpha > 1$. In addition suppose two finite sets $\mathcal{C}_1\subset \mathcal{C}_2\subset \{(i,j): i,j\in\mathbf{Z}\}$. Define $\mathcal{D}_k$ as the $\sigma-$field
generated by $e_{i}^{(j)}, (i,j)\in\mathcal{C}_k$ for $k = 1,2$. Then
\begin{equation}
\begin{aligned}
\Vert\mathbf{E}\epsilon_i^{(j)}|\mathcal{D}_1 - \mathbf{E}\epsilon_i^{(j)}|\mathcal{D}_2\Vert_M
\leq \sum_{(s,t)\in\mathcal{C}_2 - \mathcal{C}_1}\delta_{s - i}^{(t - j)}
\label{eq.delta_moment}
\end{aligned}
\end{equation}
\label{lemma.delta_moment}
\end{lemma}
\begin{proof}
Suppose $\mathcal{C}_2 - \mathcal{C}_1 = \{(i_1,j_1),\cdots, (i_K, j_K)\}$ for some integer $K$. Define $\mathcal{T}_0 = \mathcal{D}_1$ and $\mathcal{T}_k$ as the $\sigma-$field generated by $e_i^{(j)}, (i,j)\in \mathcal{C}_1\cup\{(i_1,j_1),\cdots, (i_k,j_k)\}, k = 1,\cdots, K$. Then
$\mathcal{T}_K = \mathcal{D}_2$, and
$$
\Vert\mathbf{E}\epsilon_i^{(j)}|\mathcal{D}_1 - \mathbf{E}\epsilon_i^{(j)}|\mathcal{D}_2\Vert_M\leq \sum_{k = 1}^K
\Vert\mathbf{E}\epsilon_i^{(j)}|\mathcal{T}_k - \mathbf{E}\epsilon_i^{(j)}|\mathcal{T}_{k - 1}\Vert_M
$$
Since $\mathbf{E}\epsilon_i^{(j)}|\mathcal{T}_{k - 1} = \mathbf{E}\epsilon_{i,(i_k - i)}^{(j),(j_k - j)}|\mathcal{T}_k$,
$$
\Vert\mathbf{E}\epsilon_i^{(j)}|\mathcal{T}_k - \mathbf{E}\epsilon_i^{(j)}|\mathcal{T}_{k - 1}\Vert_M
\leq\Vert\epsilon_i^{(j)} - \epsilon_{i, (i_k - i)}^{(j), (j_k - j)}\Vert_M\leq \delta_{i_k - i}^{(j_k - j)}
$$
and we prove \eqref{eq.delta_moment}.
\end{proof}

\begin{lemma}
Suppose $\epsilon_i^{(j)}, i = 1,\cdots, n, j = 1,\cdots, m$ are $(M,\alpha)-$short range dependent random variables with $M>4, \alpha > 3$,
$l\geq 1$ is a given integer, and $a_i^{(j)}, i = 1,\cdots, n, j = 1,\cdots, m$ are real numbers. Then there exists a constant $C$ independent with $a_i^{(j)}$ and $l$ such that
\begin{equation}
\begin{aligned}
\Vert\sum_{i = 1}^n\sum_{j = 1}^m a_i^{(j)}(\mathbf{E}\epsilon_i^{(j)}|\mathcal{H}_{i, (l,l - 1)}^{(j), (l - 1,l - 1)} -
\mathbf{E}\epsilon_i^{(j)}|\mathcal{H}_{i, (l - 1,l - 1)}^{(j), (l - 1,l - 1)})\Vert_M\\
\leq
C\sqrt{\sum_{i = 1}^n\Vert \sum_{j = 1}^m a_i^{(j)}(\mathbf{E}\epsilon_i^{(j)}|\mathcal{H}_{i, (l,l - 1)}^{(j), (l - 1,l - 1)} -
\mathbf{E}\epsilon_i^{(j)}|\mathcal{H}_{i, (l - 1,l - 1)}^{(j), (l - 1,l - 1)})\Vert_M^2}
\\
\Vert\sum_{i = 1}^n\sum_{j = 1}^m a_i^{(j)}(\mathbf{E}\epsilon_i^{(j)}|\mathcal{H}_{i, (l,l)}^{(j), (l - 1,l - 1)} -
\mathbf{E}\epsilon_i^{(j)}|\mathcal{H}_{i, (l,l - 1)}^{(j), (l - 1,l - 1)})\Vert_M\\
\leq
C\sqrt{\sum_{i = 1}^n \Vert\sum_{j = 1}^m a_i^{(j)}(\mathbf{E}\epsilon_i^{(j)}|\mathcal{H}_{i, (l,l)}^{(j), (l - 1,l - 1)} -
\mathbf{E}\epsilon_i^{(j)}|\mathcal{H}_{i, (l,l - 1)}^{(j), (l - 1,l - 1)})\Vert_M^2}
\\
\Vert\sum_{i = 1}^n\sum_{j = 1}^m a_i^{(j)}(\mathbf{E}\epsilon_i^{(j)}|\mathcal{H}_{i, (l,l)}^{(j), (l,l - 1)} -
\mathbf{E}\epsilon_i^{(j)}|\mathcal{H}_{i, (l,l)}^{(j), (l - 1,l - 1)})\Vert_M\\
\leq C\sqrt{\sum_{j = 1}^m \Vert \sum_{i = 1}^n a_i^{(j)}(\mathbf{E}\epsilon_i^{(j)}|\mathcal{H}_{i, (l,l)}^{(j), (l,l - 1)} -
\mathbf{E}\epsilon_i^{(j)}|\mathcal{H}_{i, (l,l)}^{(j), (l - 1,l - 1)})\Vert_M^2}\\
\Vert\sum_{i = 1}^n\sum_{j = 1}^m a_i^{(j)}(\mathbf{E}\epsilon_i^{(j)}|\mathcal{H}_{i, (l,l)}^{(j), (l,l)} -
\mathbf{E}\epsilon_i^{(j)}|\mathcal{H}_{i, (l,l)}^{(j), (l,l - 1)})\Vert_M\\
\leq C\sqrt{\sum_{j = 1}^m\Vert\sum_{i = 1}^na_i^{(j)}(\mathbf{E}\epsilon_i^{(j)}|\mathcal{H}_{i, (l,l)}^{(j), (l,l)} -
\mathbf{E}\epsilon_i^{(j)}|\mathcal{H}_{i, (l,l)}^{(j), (l,l - 1)})\Vert_M^2}
\end{aligned}
\end{equation}
\label{lemma.calculate_martingale_difference}
\end{lemma}
\begin{proof}
Define $P_{s,l} = \sum_{i = n - s}^n\sum_{j = 1}^m a_i^{(j)}(\mathbf{E}\epsilon_i^{(j)}|\mathcal{H}_{i, (l,l - 1)}^{(j), (l - 1,l - 1)} -
\mathbf{E}\epsilon_i^{(j)}|\mathcal{H}_{i, (l - 1,l - 1)}^{(j), (l - 1,l - 1)})$ for $s = 0, 1,\cdots, n - 1$, then $P_{s,l}$ is
$\mathcal{Z}_{(n - s - l, n + l - 1)}^{(2-l, m + l - 1)}-$measurable. Moreover, $\mathcal{Z}_{(n - s - l, n + l - 1)}^{(2-l, m + l - 1)}\subset \mathcal{Z}_{(n - s - l - 1, n + l - 1)}^{(2-l, m + l - 1)}$ and
\begin{align*}
P_{s+1,l} - P_{s,l} = \sum_{j = 1}^m a_{n - s  - 1}^{(j)}(\mathbf{E}\epsilon_{n - s - 1}^{(j)}|\mathcal{H}_{n - s - 1, (l,l - 1)}^{(j), (l - 1,l - 1)} -
\mathbf{E}\epsilon_{n - s - 1}^{(j)}|\mathcal{H}_{n - s - 1, (l - 1,l - 1)}^{(j), (l - 1,l - 1)})
\end{align*}
From lemma \ref{lemma.conditional_prob},
\begin{equation}
\begin{aligned}
\mathbf{E}
\left(\mathbf{E}\epsilon_{n - s - 1}^{(j)}|\mathcal{H}_{n - s - 1, (l,l - 1)}^{(j), (l - 1,l - 1)}\right)|\mathcal{Z}_{(n - s - l, n + l - 1)}^{(2-l, m + l - 1)} = \mathbf{E}\epsilon_{n - s - l}^{(j)}|\mathcal{H}_{n - s - 1, (l - 1, l - 1)}^{(j), (l - 1, l - 1)}\\
\Rightarrow \mathbf{E}(P_{s + 1, l} - P_{s,l})|\mathcal{Z}_{(n - s - l, n + l - 1)}^{(2-l, m + l - 1)}\\
= \sum_{j = 1}^m a_{n - s - 1}^{(j)}\left(\left(\mathbf{E}\epsilon_{n - s - 1}^{(j)}|\mathcal{H}_{n - s - 1, (l,l - 1)}^{(j), (l - 1,l - 1)}\right)|\mathcal{Z}_{(n - s - l, n + l - 1)}^{(2-l, m + l - 1)} - \mathbf{E}\epsilon_{n - s - l}^{(j)}|\mathcal{H}_{n - s - 1, (l - 1, l - 1)}^{(j), (l - 1, l - 1)}\right)\\
= 0
\end{aligned}
\end{equation}
so $P_{s,l}$ form a martingale and from \cite{MR0400380},
\begin{equation}
\begin{aligned}
\Vert\sum_{i = 1}^n\sum_{j = 1}^m a_i^{(j)}(\mathbf{E}\epsilon_i^{(j)}|\mathcal{H}_{i, (l,l - 1)}^{(j), (l - 1,l - 1)} -
\mathbf{E}\epsilon_i^{(j)}|\mathcal{H}_{i, (l - 1,l - 1)}^{(j), (l - 1,l - 1)})\Vert_M\\
\leq C\Vert\sum_{i = 1}^n\left(\sum_{j = 1}^m a_i^{(j)}(\mathbf{E}\epsilon_i^{(j)}|\mathcal{H}_{i, (l,l - 1)}^{(j), (l - 1,l - 1)} -
\mathbf{E}\epsilon_i^{(j)}|\mathcal{H}_{i, (l - 1,l - 1)}^{(j), (l - 1,l - 1)})\right)^2\Vert_{M/2}^{1/2}\\
\leq C\sqrt{\sum_{i = 1}^n\Vert\sum_{j = 1}^m a_i^{(j)}(\mathbf{E}\epsilon_i^{(j)}|\mathcal{H}_{i, (l,l - 1)}^{(j), (l - 1,l - 1)} -
\mathbf{E}\epsilon_i^{(j)}|\mathcal{H}_{i, (l - 1,l - 1)}^{(j), (l - 1,l - 1)})\Vert_M^2}
\end{aligned}
\end{equation}
here $C$ is a constant only depending on $M$.

Similarly define $Q_{s, l} = \sum_{i = 1}^s\sum_{j = 1}^m a_i^{(j)}(\mathbf{E}\epsilon_i^{(j)}|\mathcal{H}_{i, (l,l)}^{(j), (l - 1,l - 1)} -
\mathbf{E}\epsilon_i^{(j)}|\mathcal{H}_{i, (l,l - 1)}^{(j), (l - 1,l - 1)})$ for $s = 1,\cdots, n$, then $Q_{s,l}$ is $\mathcal{Z}_{(1-l,s+l)}^{(2-l,m+l-1)}-$ measurable, $\mathcal{Z}_{(1-l,s+l)}^{(2-l,m+l-1)}\subset \mathcal{Z}_{(1-l,s+l + 1)}^{(2-l,m+l-1)}$, and
$$
Q_{s+1, l} - Q_{s, l} = \sum_{j = 1}^m a_{s + 1}^{(j)}(\mathbf{E}\epsilon_{s + 1}^{(j)}|\mathcal{H}_{s+1, (l,l)}^{(j), (l - 1,l - 1)} -
\mathbf{E}\epsilon_{s+1}^{(j)}|\mathcal{H}_{s+1, (l,l - 1)}^{(j), (l - 1,l - 1)})
$$
From lemma \ref{lemma.conditional_prob},
\begin{equation}
\begin{aligned}
\mathbf{E}\left(\mathbf{E}\epsilon_{s + 1}^{(j)}|\mathcal{H}_{s+1, (l,l)}^{(j), (l - 1,l - 1)}\right)|\mathcal{Z}_{(1-l,s+l)}^{(2-l,m+l-1)}
= \mathbf{E}\epsilon_{s+1}^{(j)}|\mathcal{H}_{s+1, (l,l - 1)}^{(j), (l - 1,l - 1)}\\
\Rightarrow \mathbf{E}\left(Q_{s+1, l} - Q_{s, l}\right)|\mathcal{Z}_{(1-l,s+l)}^{(2-l,m+l-1)}\\
= \sum_{j = 1}^m a_{s + 1}^{(j)}(\mathbf{E}\left(\mathbf{E}\epsilon_{s + 1}^{(j)}|\mathcal{H}_{s+1, (l,l)}^{(j), (l - 1,l - 1)}\right)|\mathcal{Z}_{(1-l,s+l)}^{(2-l,m+l-1)} -
\mathbf{E}\epsilon_{s+1}^{(j)}|\mathcal{H}_{s+1, (l,l - 1)}^{(j), (l - 1,l - 1)}) = 0
\end{aligned}
\end{equation}
so $Q_{s, l}$ form  a martingale and
\begin{equation}
\begin{aligned}
\Vert\sum_{i = 1}^n\sum_{j = 1}^m a_i^{(j)}(\mathbf{E}\epsilon_i^{(j)}|\mathcal{H}_{i, (l,l)}^{(j), (l - 1,l - 1)} -
\mathbf{E}\epsilon_i^{(j)}|\mathcal{H}_{i, (l,l - 1)}^{(j), (l - 1,l - 1)})\Vert_M\\
\leq C\sqrt{\sum_{i = 1}^n\Vert \sum_{j = 1}^m a_i^{(j)}(\mathbf{E}\epsilon_i^{(j)}|\mathcal{H}_{i, (l,l)}^{(j), (l - 1,l - 1)} -
\mathbf{E}\epsilon_i^{(j)}|\mathcal{H}_{i, (l,l - 1)}^{(j), (l - 1,l - 1)})\Vert_M^2}
\end{aligned}
\end{equation}
Define $R_{s, l} = \sum_{j = m - s}^m \sum_{i = 1}^n a_i^{(j)}(\mathbf{E}\epsilon_i^{(j)}|\mathcal{H}_{i, (l,l)}^{(j), (l,l - 1)} -
\mathbf{E}\epsilon_i^{(j)}|\mathcal{H}_{i, (l,l)}^{(j), (l - 1,l - 1)})$ for $s = 0,\cdots , m- 1$, then
$R_{s,l}$ is $\mathcal{Z}_{(1-l, n + l)}^{(m - s - l, m + l - 1)}$ measurable and
$$
R_{s + 1, l} - R_{s, l} = \sum_{i = 1}^n a_i^{(m-s-1)}(\mathbf{E}\epsilon_i^{(m-s-1)}|\mathcal{H}_{i, (l,l)}^{(m-s-1), (l,l - 1)} -
\mathbf{E}\epsilon_i^{(m-s-1)}|\mathcal{H}_{i, (l,l)}^{(m-s-1), (l - 1,l - 1)})
$$
From lemma \ref{lemma.conditional_prob},
\begin{equation}
\begin{aligned}
\mathbf{E}\left(\mathbf{E}\epsilon_i^{(m-s-1)}|\mathcal{H}_{i, (l,l)}^{(m-s-1), (l,l - 1)}\right)|\mathcal{Z}_{(1-l, n + l)}^{(m - s - l, m + l - 1)}
= \mathbf{E}\epsilon_i^{(m-s-1)}|\mathcal{H}_{i, (l, l)}^{(m - s - 1), (l - 1, l - 1)}\\
\Rightarrow \mathbf{E}\left(R_{s + 1, l} - R_{s, l}\right)|\mathcal{Z}_{(1-l, n + l)}^{(m - s - l, m + l - 1)}\\
=\sum_{i = 1}^n a_i^{(m-s-1)}\left(\mathbf{E}\left(\mathbf{E}\epsilon_i^{(m-s-1)}|\mathcal{H}_{i, (l,l)}^{(m-s-1), (l,l - 1)}\right)|\mathcal{Z}_{(1-l, n + l)}^{(m - s - l, m + l - 1)}
- \mathbf{E}\epsilon_i^{(m-s-1)}|\mathcal{H}_{i, (l, l)}^{(m - s - 1), (l - 1, l - 1)}\right)\\
= 0
\end{aligned}
\end{equation}
so $R_{s,l}$ form a martingale. Therefore,
\begin{equation}
\begin{aligned}
\Vert\sum_{i = 1}^n\sum_{j = 1}^m a_i^{(j)}(\mathbf{E}\epsilon_i^{(j)}|\mathcal{H}_{i, (l,l)}^{(j), (l,l - 1)} -
\mathbf{E}\epsilon_i^{(j)}|\mathcal{H}_{i, (l,l)}^{(j), (l - 1,l - 1)})\Vert_M\\
\leq C\sqrt{\sum_{j = 1}^n\Vert\sum_{i = 1}^na_i^{(j)}(\mathbf{E}\epsilon_i^{(j)}|\mathcal{H}_{i, (l,l)}^{(j), (l,l - 1)} -
\mathbf{E}\epsilon_i^{(j)}|\mathcal{H}_{i, (l,l)}^{(j), (l - 1,l - 1)})\Vert_M^2}
\end{aligned}
\end{equation}
Define $T_{s,l} = \sum_{j = 1}^s\sum_{i = 1}^n a_i^{(j)}(\mathbf{E}\epsilon_i^{(j)}|\mathcal{H}_{i, (l,l)}^{(j), (l,l)} -
\mathbf{E}\epsilon_i^{(j)}|\mathcal{H}_{i, (l,l)}^{(j), (l,l - 1)})$ for $s = 1,\cdots, m$, then $T_{s,l}$ is
$\mathcal{Z}_{(1-l,  n + l)}^{(1 - l, s + l)}$ measurable, and
$$
T_{s + 1, l} - T_{s,l} = \sum_{i = 1}^n a_i^{(s+1)}(\mathbf{E}\epsilon_i^{(s+1)}|\mathcal{H}_{i, (l,l)}^{(s+1), (l,l)} -
\mathbf{E}\epsilon_i^{(s+1)}|\mathcal{H}_{i, (l,l)}^{(s+1), (l,l - 1)})
$$
From lemma \ref{lemma.conditional_prob},
\begin{equation}
\begin{aligned}
\mathbf{E}\left(\mathbf{E}\epsilon_i^{(s+1)}|\mathcal{H}_{i, (l,l)}^{(s+1), (l,l)}\right)|\mathcal{Z}_{(1-l,  n + l)}^{(1 - l, s + l)}
= \mathbf{E}\epsilon_i^{(s+1)}|\mathcal{H}_{i, (l,l)}^{(s  + 1),(l,l - 1)}
\end{aligned}
\end{equation}
so $T_{s,l}$ form a martingale, and
\begin{equation}
\begin{aligned}
\Vert\sum_{i = 1}^n\sum_{j = 1}^m a_i^{(j)}(\mathbf{E}\epsilon_i^{(j)}|\mathcal{H}_{i, (l,l)}^{(j), (l,l)} -
\mathbf{E}\epsilon_i^{(j)}|\mathcal{H}_{i, (l,l)}^{(j), (l,l - 1)})\Vert_M\\
\leq C\sqrt{\sum_{j = 1}^m\Vert\sum_{i = 1}^na_i^{(j)}(\mathbf{E}\epsilon_i^{(j)}|\mathcal{H}_{i, (l,l)}^{(j), (l,l)} -
\mathbf{E}\epsilon_i^{(j)}|\mathcal{H}_{i, (l,l)}^{(j), (l,l - 1)})\Vert_M^2}
\end{aligned}
\end{equation}
\end{proof}
\begin{lemma}
Suppose $\epsilon_i^{(j)}, i = 1,\cdots, n, j = 1,\cdots, m$ are $(M,\alpha)-$short range dependent random variables with $M>4, \alpha > 3$,
$l\geq 1$ is a given integer, and $a_i^{(j)}, i = 1,\cdots, n, j = 1,\cdots, m$ are real numbers. Then there exists a constant $C$ independent with $a_i^{(j)}$ and $i,j,l$ such that
\begin{equation}
\begin{aligned}
\Vert \sum_{j = 1}^m a_i^{(j)}(\mathbf{E}\epsilon_i^{(j)}|\mathcal{H}_{i, (l,l - 1)}^{(j), (l - 1,l - 1)} -
\mathbf{E}\epsilon_i^{(j)}|\mathcal{H}_{i, (l - 1,l - 1)}^{(j), (l - 1,l - 1)})\Vert_M
\leq C\sqrt{\sum_{j = 1}^m a_i^{(j)2}}\times (1 + l)^{2 -\alpha}\\
\Vert\sum_{j = 1}^m a_i^{(j)}(\mathbf{E}\epsilon_i^{(j)}|\mathcal{H}_{i, (l,l)}^{(j), (l - 1,l - 1)} -
\mathbf{E}\epsilon_i^{(j)}|\mathcal{H}_{i, (l,l - 1)}^{(j), (l - 1,l - 1)})\Vert_M\leq C\sqrt{\sum_{j = 1}^m a_i^{(j)2}}\times (1 + l)^{2 -\alpha}\\
\Vert \sum_{i = 1}^n a_i^{(j)}(\mathbf{E}\epsilon_i^{(j)}|\mathcal{H}_{i, (l,l)}^{(j), (l,l - 1)} -
\mathbf{E}\epsilon_i^{(j)}|\mathcal{H}_{i, (l,l)}^{(j), (l - 1,l - 1)})\Vert_M\leq C\sqrt{\sum_{i = 1}^n a_i^{(j)2}}\times (1 + l)^{2-\alpha}\\
\Vert\sum_{i = 1}^na_i^{(j)}(\mathbf{E}\epsilon_i^{(j)}|\mathcal{H}_{i, (l,l)}^{(j), (l,l)} -
\mathbf{E}\epsilon_i^{(j)}|\mathcal{H}_{i, (l,l)}^{(j), (l,l - 1)})\Vert_M\leq C\sqrt{\sum_{i = 1}^n a_i^{(j)2}}\times (1 + l)^{2-\alpha}
\end{aligned}
\label{eq.sub_difference}
\end{equation}
\label{lemma.second_dimension_difference}
\end{lemma}
\begin{proof}
Define $p_{i, l}^{(j)} = \mathbf{E}\epsilon_i^{(j)}|\mathcal{H}_{i, (l,l - 1)}^{(j), (l - 1,l - 1)} -
\mathbf{E}\epsilon_i^{(j)}|\mathcal{H}_{i, (l - 1,l - 1)}^{(j), (l - 1,l - 1)}$, then
\begin{align*}
p_{i,l}^{(j)} = \mathbf{E}p_{i,l}^{(j)}|\mathcal{H}_{i, (l,l - 1)}^{(j), (l - 1,l - 1)}
= \mathbf{E}p_{i,l}^{(j)}|\mathcal{H}_{i, (l, l - 1)}^{(j), (0,0)} + \sum_{k = 1}^{l - 1}(\mathbf{E}p_{i,l}^{(j)}|\mathcal{H}_{i, (l,l - 1)}^{(j), (k, k - 1)} - \mathbf{E}p_{i,l}^{(j)}|\mathcal{H}_{i, (l,l - 1)}^{(j), (k - 1, k - 1)})\\
+ \sum_{k = 1}^{l - 1}(\mathbf{E}p_{i,l}^{(j)}|\mathcal{H}_{i, (l,l - 1)}^{(j), (k, k)} - \mathbf{E}p_{i,l}^{(j)}|\mathcal{H}_{i, (l,l - 1)}^{(j), (k, k - 1)})
\end{align*}
therefore
\begin{equation}
\begin{aligned}
\Vert \sum_{j = 1}^m a_i^{(j)}p_{i,l}^{(j)}\Vert_M\leq \Vert\sum_{j = 1}^m a_i^{(j)}\mathbf{E}p_{i,l}^{(j)}|\mathcal{H}_{i, (l, l - 1)}^{(j), (0,0)}\Vert_M\\
+ \sum_{k = 1}^{l - 1}\Vert\sum_{j = 1}^m a_i^{(j)}(\mathbf{E}p_{i,l}^{(j)}|\mathcal{H}_{i, (l,l - 1)}^{(j), (k, k - 1)} - \mathbf{E}p_{i,l}^{(j)}|\mathcal{H}_{i, (l,l - 1)}^{(j), (k - 1, k - 1)})\Vert_M\\
+\sum_{k = 1}^{l - 1}\Vert\sum_{j = 1}^m a_i^{(j)}(\mathbf{E}p_{i,l}^{(j)}|\mathcal{H}_{i, (l,l - 1)}^{(j), (k, k)} - \mathbf{E}p_{i,l}^{(j)}|\mathcal{H}_{i, (l,l - 1)}^{(j), (k, k - 1)})\Vert_M
\label{eq.separate}
\end{aligned}
\end{equation}
For $\mathbf{E}p_{i,l}^{(j)}|\mathcal{H}_{i, (l, l - 1)}^{(j), (0,0)}$ are mutually independent for $j = 1,\cdots, m$, from \cite{MR0133849},
\begin{equation}
\begin{aligned}
\Vert\sum_{j = 1}^m a_i^{(j)}\mathbf{E}p_{i,l}^{(j)}|\mathcal{H}_{i, (l, l - 1)}^{(j), (0,0)}\Vert_M\leq C\sqrt{\sum_{j = 1}^m
a_i^{(j)2}\Vert\mathbf{E}p_{i,l}^{(j)}|\mathcal{H}_{i, (l, l - 1)}^{(j), (0,0)}\Vert_M^2}\\
\leq C\sqrt{\sum_{j = 1}^m a_i^{(j)2}}\times \max_{i = 1,\cdots, n, j = 1,\cdots, m}\Vert \mathbf{E}p_{i,l}^{(j)}|\mathcal{H}_{i, (l, l - 1)}^{(j), (0,0)}\Vert_M
\end{aligned}
\end{equation}
From lemma \ref{lemma.conditional_prob},
\begin{align*}
\mathbf{E}p_{i,l}^{(j)}|\mathcal{H}_{i, (l, l - 1)}^{(j), (0,0)}
= \mathbf{E}\epsilon_i^{(j)}|\mathcal{H}_{i, (l, l - 1)}^{(j), (0, 0)} - \mathbf{E}\epsilon_i^{(j)}|\mathcal{H}_{i, (l - 1, l - 1)}^{(j), (0, 0)}\\
\Rightarrow \Vert\mathbf{E}p_{i,l}^{(j)}|\mathcal{H}_{i, (l, l - 1)}^{(j), (0,0)} \Vert_M =
\Vert \mathbf{E}(\epsilon_i^{(j)} - \epsilon_{i,-l}^{(j), (0)})|\mathcal{H}_{i, (l, l - 1)}^{(j), (0, 0)}\Vert_M\leq \delta_{-l}^{(0)}
\end{align*}
Therefore
\begin{equation}
\Vert\sum_{j = 1}^m a_i^{(j)}\mathbf{E}p_{i,l}^{(j)}|\mathcal{H}_{i, (l, l - 1)}^{(j), (0,0)}\Vert_M\leq C\delta_{-l}^{(0)}\times \sqrt{\sum_{j = 1}^m a_i^{(j)2}}
\label{eq.separateFirst}
\end{equation}
On the other hand, define
$P_{i,l}^{(s)} =  \sum_{j = m - s}^m a_i^{(j)}(\mathbf{E}p_{i,l}^{(j)}|\mathcal{H}_{i, (l,l - 1)}^{(j), (k, k - 1)} - \mathbf{E}p_{i,l}^{(j)}|\mathcal{H}_{i, (l,l - 1)}^{(j), (k - 1, k - 1)})$ for $s = 0,\cdots, m - 1$, then
$P_{i,l}^{(s)}$ is $\mathcal{Z}_{(i-l, i+l-1)}^{(m-s-k, m + k - 1)}$ measurable,
$\mathcal{Z}_{(i-l, i+l-1)}^{(m-s-k, m + k - 1)}\subset \mathcal{Z}_{(i-l, i+l-1)}^{(m-s-1-k, m + k - 1)}$, and
$$
P_{i,l}^{(s+1)} - P_{i,l}^{(s)} = a_i^{(m - s - 1)}(\mathbf{E}p_{i,l}^{(m - s - 1)}|\mathcal{H}_{i, (l,l - 1)}^{(m - s - 1), (k, k - 1)} - \mathbf{E}p_{i,l}^{(m - s - 1)}|\mathcal{H}_{i, (l,l - 1)}^{(m - s - 1), (k - 1, k - 1)})
$$
From lemma \ref{lemma.conditional_prob}
\begin{align*}
\mathbf{E}\left(\mathbf{E}p_{i,l}^{(m - s - 1)}|\mathcal{H}_{i, (l,l - 1)}^{(m - s - 1), (k, k - 1)}\right)|\mathcal{Z}_{(i-l, i+l-1)}^{(m-s-k, m + k - 1)} = \mathbf{E}p_{i,l}^{(m - s - 1)}|\mathcal{H}_{i, (l, l - 1)}^{(m - s - 1), (k - 1, k - 1)}\\
\end{align*}
so $P_{i,l}^{(s)}$ form a martingale, and from \cite{MR0400380},
\begin{equation}
\begin{aligned}
\Vert\sum_{j = 1}^m a_i^{(j)}(\mathbf{E}p_{i,l}^{(j)}|\mathcal{H}_{i, (l,l - 1)}^{(j), (k, k - 1)} - \mathbf{E}p_{i,l}^{(j)}|\mathcal{H}_{i, (l,l - 1)}^{(j), (k - 1, k - 1)})\Vert_M\\
\leq C\sqrt{\sum_{j = 1}^m a_i^{(j)2}\Vert\mathbf{E}p_{i,l}^{(j)}|\mathcal{H}_{i, (l,l - 1)}^{(j), (k, k - 1)} - \mathbf{E}p_{i,l}^{(j)}|\mathcal{H}_{i, (l,l - 1)}^{(j), (k - 1, k - 1)}\Vert_M^2}\\
\leq C\sqrt{\sum_{j = 1}^m a_i^{(j)2}}\times \max_{i = 1,\cdots, n, j = 1,\cdots, m}\Vert\mathbf{E}p_{i,l}^{(j)}|\mathcal{H}_{i, (l,l - 1)}^{(j), (k, k - 1)} - \mathbf{E}p_{i,l}^{(j)}|\mathcal{H}_{i, (l,l - 1)}^{(j), (k - 1, k - 1)}\Vert_M
\end{aligned}
\end{equation}
From lemma \ref{lemma.conditional_prob} and \ref{lemma.delta_moment},
\begin{align*}
\mathbf{E}p_{i,l}^{(j)}|\mathcal{H}_{i, (l,l - 1)}^{(j), (k, k - 1)} - \mathbf{E}p_{i,l}^{(j)}|\mathcal{H}_{i, (l,l - 1)}^{(j), (k - 1, k - 1)}
= \mathbf{E}\epsilon_i^{(j)}|\mathcal{H}_{i, (l, l - 1)}^{(j), (k, k - 1)} - \mathbf{E}\epsilon_i^{(j)}|\mathcal{H}_{i, (l - 1, l- 1)}^{(j), (k, k - 1)}\\
- \mathbf{E}\epsilon_i^{(j)}|\mathcal{H}_{i, (l, l - 1)}^{(j), (k - 1, k - 1)} + \mathbf{E}\epsilon_i^{(j)}|\mathcal{H}_{i, (l - 1, l - 1)}^{(j), (k - 1, k - 1)}\\
\Rightarrow \Vert\mathbf{E}p_{i,l}^{(j)}|\mathcal{H}_{i, (l,l - 1)}^{(j), (k, k - 1)} - \mathbf{E}p_{i,l}^{(j)}|\mathcal{H}_{i, (l,l - 1)}^{(j), (k - 1, k - 1)}\Vert_M
\leq \Vert\mathbf{E}\epsilon_i^{(j)}|\mathcal{H}_{i, (l, l - 1)}^{(j), (k, k - 1)} - \mathbf{E}\epsilon_i^{(j)}|\mathcal{H}_{i, (l - 1, l- 1)}^{(j), (k, k - 1)}\Vert_M\\
+ \Vert\mathbf{E}\epsilon_i^{(j)}|\mathcal{H}_{i, (l, l - 1)}^{(j), (k - 1, k - 1)} - \mathbf{E}\epsilon_i^{(j)}|\mathcal{H}_{i, (l - 1, l - 1)}^{(j), (k - 1, k - 1)}\Vert_M\\
\leq 2\sum_{s = -k}^k\delta_{-l}^{(s)}\leq C(1 + l)^{1-\alpha}
\end{align*}
therefore
\begin{equation}
\Vert\sum_{j = 1}^m a_i^{(j)}(\mathbf{E}p_{i,l}^{(j)}|\mathcal{H}_{i, (l,l - 1)}^{(j), (k, k - 1)} - \mathbf{E}p_{i,l}^{(j)}|\mathcal{H}_{i, (l,l - 1)}^{(j), (k - 1, k - 1)})\Vert_M\leq C\sqrt{\sum_{j = 1}^m a_i^{(j)2}}\times (1 + l)^{1-\alpha}
\label{eq.separate_2}
\end{equation}
for a constant $C$. Similarly define $\widetilde{P}_{i,l}^{(s)} = \sum_{j = 1}^s a_i^{(j)}(\mathbf{E}p_{i,l}^{(j)}|\mathcal{H}_{i, (l,l - 1)}^{(j), (k, k)} - \mathbf{E}p_{i,l}^{(j)}|\mathcal{H}_{i, (l,l - 1)}^{(j), (k, k - 1)})$ for $s = 1,\cdots, m$, then
$\widetilde{P}_{i,l}^{(s)}$ is $\mathcal{Z}_{(i - l, i + l - 1)}^{(1 - k, s + k)}$ measurable,
$\mathcal{Z}_{(i - l, i + l - 1)}^{(1 - k, s + k)}\subset \mathcal{Z}_{(i - l, i + l - 1)}^{(1 - k, s + k + 1)}$, and
$$
\widetilde{P}_{i,l}^{(s + 1)} - \widetilde{P}_{i,l}^{(s)} = a_i^{(s+1)}(\mathbf{E}p_{i,l}^{(s+1)}|\mathcal{H}_{i, (l,l - 1)}^{(s+1), (k, k)} - \mathbf{E}p_{i,l}^{(s+1)}|\mathcal{H}_{i, (l,l - 1)}^{(s+1), (k, k - 1)})
$$
Since
$$
\mathbf{E}\left(\mathbf{E}p_{i,l}^{(s+1)}|\mathcal{H}_{i, (l,l - 1)}^{(s+1), (k, k)}\right)|\mathcal{Z}_{(i - l, i + l - 1)}^{(1 - k, s + k)}
= \mathbf{E}p_{i,l}^{(s+1)}|\mathcal{H}_{i, (l, l - 1)}^{(s+1), (k, k - 1)}
$$
and
\begin{align*}
\mathbf{E}p_{i,l}^{(j)}|\mathcal{H}_{i, (l,l - 1)}^{(j), (k, k)} - \mathbf{E}p_{i,l}^{(j)}|\mathcal{H}_{i, (l,l - 1)}^{(j), (k, k - 1)}
= \mathbf{E}\epsilon_i^{(j)}|\mathcal{H}_{i, (l, l - 1)}^{(j), (k,k)} - \mathbf{E}\epsilon_i^{(j)}|\mathcal{H}_{i, (l - 1, l - 1)}^{(j), (k,k)}\\
- \mathbf{E}\epsilon_i^{(j)}|\mathcal{H}_{i, (l, l - 1)}^{(j), (k,k - 1)} + \mathbf{E}\epsilon_i^{(j)}|\mathcal{H}_{i, (l - 1, l - 1)}^{(j), (k,k - 1)}
\end{align*}
from \cite{MR0400380} and lemma \ref{lemma.delta_moment},
\begin{equation}
\begin{aligned}
\Vert\sum_{j = 1}^m a_i^{(j)}(\mathbf{E}p_{i,l}^{(j)}|\mathcal{H}_{i, (l,l - 1)}^{(j), (k, k)} - \mathbf{E}p_{i,l}^{(j)}|\mathcal{H}_{i, (l,l - 1)}^{(j), (k, k - 1)})\Vert_M\\
\leq C\sqrt{\sum_{j = 1}^m a_i^{(j)2}\Vert\mathbf{E}p_{i,l}^{(j)}|\mathcal{H}_{i, (l,l - 1)}^{(j), (k, k)} - \mathbf{E}p_{i,l}^{(j)}|\mathcal{H}_{i, (l,l - 1)}^{(j), (k, k - 1)}\Vert_M^2}\\
\leq 2C\sqrt{\sum_{j = 1}^m a_i^{(j)2}}\times \sum_{s = -k}^k \delta_{-l}^{(s)}\leq C^\prime\sqrt{\sum_{j = 1}^m a_i^{(j)2}}\times(1 + l)^{1-\alpha}
\end{aligned}
\label{eq.separate_3}
\end{equation}
with $C^\prime$ being a constant. From \eqref{eq.separate}, \eqref{eq.separateFirst}, \eqref{eq.separate_2}, and \eqref{eq.separate_3},
\begin{equation}
\begin{aligned}
\Vert \sum_{j = 1}^m a_i^{(j)}(\mathbf{E}\epsilon_i^{(j)}|\mathcal{H}_{i, (l,l - 1)}^{(j), (l - 1,l - 1)} -
\mathbf{E}\epsilon_i^{(j)}|\mathcal{H}_{i, (l - 1,l - 1)}^{(j), (l - 1,l - 1)})\Vert_M
\leq C\sqrt{\sum_{j = 1}^m a_i^{(j)2}}\times (1 + l)^{2-\alpha}
\end{aligned}
\end{equation}
Similarly, define $q_{i,l}^{(j)} = \mathbf{E}\epsilon_i^{(j)}|\mathcal{H}_{i, (l,l)}^{(j), (l - 1,l - 1)} -
\mathbf{E}\epsilon_i^{(j)}|\mathcal{H}_{i, (l,l - 1)}^{(j), (l - 1,l - 1)}$, then
\begin{equation}
\begin{aligned}
\Vert
\sum_{j = 1}^m a_i^{(j)}q_{i,l}^{(j)}
\Vert_M\leq C\sqrt{\sum_{j = 1}^m a_i^{(j)2}\Vert\mathbf{E}q_{i,l}^{(j)}|\mathcal{H}_{i, (l,l)}^{(j), (0,0)}\Vert_M^2}\\
+ \sum_{k = 1}^{l - 1}\Vert\sum_{j = 1}^m a_i^{(j)}(\mathbf{E}q_{i,l}^{(j)}|\mathcal{H}_{i,(l,l)}^{(j), (k, k - 1)} - \mathbf{E}q_{i,l}^{(j)}|\mathcal{H}_{i, (l,l)}^{(j), (k - 1, k - 1)})\Vert_M\\
+ \sum_{k = 1}^{l - 1}\Vert\sum_{j = 1}^m a_i^{(j)}(\mathbf{E}q_{i,l}^{(j)}|\mathcal{H}_{i,(l,l)}^{(j), (k, k)} - \mathbf{E}q_{i,l}^{(j)}|\mathcal{H}_{i, (l,l)}^{(j), (k, k - 1)})\Vert_M
\end{aligned}
\end{equation}
Notice that $\mathbf{E}q_{i,l}^{(j)}|\mathcal{H}_{i, (l,l)}^{(j), (0,0)} = \mathbf{E}\epsilon_i^{(j)}|\mathcal{H}_{i, (l,l)}^{(j), (0,0)} - \mathbf{E}\epsilon_i^{(j)}|\mathcal{H}_{i, (l,l - 1)}^{(j), (0,0)}$, we have
$$
\sqrt{\sum_{j = 1}^m a_i^{(j)2}\Vert\mathbf{E}q_{i,l}^{(j)}|\mathcal{H}_{i, (l,l)}^{(j), (0,0)}\Vert_M^2}\leq C\delta_{l}^{(0)}\sqrt{\sum_{j = 1}^m a_i^{(j)2}}
$$
Define $Q_{i,l}^{(s)} = \sum_{j = m - s}^m a_i^{(j)}(\mathbf{E}q_{i,l}^{(j)}|\mathcal{H}_{i, (l,l)}^{(j), (k, k - 1)} - \mathbf{E}q_{i,l}^{(j)}|\mathcal{H}_{i, (l,l)}^{(j), (k - 1, k - 1)})$, then $Q_{i,l}^{(s)}$ is $\mathcal{Z}_{(i-l, i + l)}^{(m - s - k, m + k - 1)}$
measurable, and
$$
\mathbf{E}\left(\mathbf{E}q_{i,l}^{(m - s - 1)}|\mathcal{H}_{i,(l,l)}^{(m - s - 1), (k, k - 1)}\right)|\mathcal{Z}_{(i-l, i + l)}^{(m - s - k, m + k - 1)} = \mathbf{E}q_{i,l}^{(m - s - 1)}|\mathcal{H}_{i, (l,l)}^{(m-s-1), (k - 1, k - 1)}
$$
For
\begin{align*}
\mathbf{E}q_{i,l}^{(j)}|\mathcal{H}_{i,(l,l)}^{(j), (k, k - 1)} - \mathbf{E}q_{i,l}^{(j)}|\mathcal{H}_{i, (l,l)}^{(j), (k - 1, k - 1)}
= \mathbf{E}\epsilon_{i}^{(j)}|\mathcal{H}_{i, (l,l)}^{(j), (k, k - 1)} - \mathbf{E}\epsilon_i^{(j)}|\mathcal{H}_{i, (l, l - 1)}^{(j), (k, k - 1)}\\
- \mathbf{E}\epsilon_i^{(j)}|\mathcal{H}_{i, (l,l)}^{(j), (k - 1, k - 1)} + \mathbf{E}\epsilon_i^{(j)}|\mathcal{H}_{i, (l, l - 1)}^{(j), (k - 1, k - 1)}
\end{align*}
we know that $Q_{i,l}^{(s)}$ is a martingale and
\begin{equation}
\begin{aligned}
\Vert\sum_{j = 1}^m a_i^{(j)}(\mathbf{E}q_{i,l}^{(j)}|\mathcal{H}_{i,(l,l)}^{(j), (k, k - 1)} - \mathbf{E}q_{i,l}^{(j)}|\mathcal{H}_{i, (l,l)}^{(j), (k - 1, k - 1)})\Vert_M\\
\leq C\sqrt{\sum_{j = 1}^m a_i^{(j)2}}\times \max_{i = 1,\cdots, n, j = 1,\cdots,m}\Vert\mathbf{E}q_{i,l}^{(j)}|\mathcal{H}_{i,(l,l)}^{(j), (k, k - 1)} - \mathbf{E}q_{i,l}^{(j)}|\mathcal{H}_{i, (l,l)}^{(j), (k - 1, k - 1)}\Vert_M\\
\leq 2C\sqrt{\sum_{j = 1}^m a_i^{(j)2}}\times\sum_{s = -k}^k\delta_l^{(s)}\leq C^\prime\sqrt{\sum_{j = 1}^m a_i^{(j)2}}\times (1 + l)^{1-\alpha}
\end{aligned}
\end{equation}
here $C^\prime$ is a constant. Define $\widetilde{Q}_{i,l}^{(s)} = \sum_{j = 1}^s a_i^{(j)}(\mathbf{E}q_{i,l}^{(j)}|\mathcal{H}_{i,(l,l)}^{(j), (k, k)} - \mathbf{E}q_{i,l}^{(j)}|\mathcal{H}_{i, (l,l)}^{(j), (k, k - 1)})$, then $Q_{i,l}^{(s)}$ is
$\mathcal{Z}_{(i - l, i + l)}^{(1 - k, s + k)}$ measurable, and
\begin{align*}
\mathbf{E}\left(\mathbf{E}q_{i,l}^{(s+1)}|\mathcal{H}_{i,(l,l)}^{(s+1), (k, k)}\right)|\mathcal{Z}_{(i - l, i + l)}^{(1 - k, s + k)}
= \mathbf{E}q_{i,l}^{(s + 1)}|\mathcal{H}_{i,(l,l)}^{(s+1), (k, k - 1)}\\
\mathbf{E}q_{i,l}^{(j)}|\mathcal{H}_{i,(l,l)}^{(j), (k, k)} - \mathbf{E}q_{i,l}^{(j)}|\mathcal{H}_{i, (l,l)}^{(j), (k, k - 1)}
=\mathbf{E}\epsilon_i^{(j)}|\mathcal{H}_{i, (l,l)}^{(j), (k,k)} - \mathbf{E}\epsilon_i^{(j)}|\mathcal{H}_{i, (l,l - 1)}^{(j), (k,k)}\\
-\mathbf{E}\epsilon_i^{(j)}|\mathcal{H}_{i, (l,l)}^{(j), (k,k - 1)} + \mathbf{E}\epsilon_i^{(j)}|\mathcal{H}_{i, (l,l - 1)}^{(j), (k,k - 1)}
\end{align*}
Therefore,
\begin{equation}
\begin{aligned}
\Vert\sum_{j = 1}^m a_i^{(j)}(\mathbf{E}q_{i,l}^{(j)}|\mathcal{H}_{i,(l,l)}^{(j), (k, k)} - \mathbf{E}q_{i,l}^{(j)}|\mathcal{H}_{i, (l,l)}^{(j), (k, k - 1)})\Vert_M\\
\leq C\sqrt{\sum_{j = 1}^m a_i^{(j)2}}\times \sum_{s = -k}^k\delta_{l}^{(s)}\leq C^\prime\sqrt{\sum_{j = 1}^m a_i^{(j)2}}\times(1 + l)^{1-\alpha}\\
\text{and } \Vert
\sum_{j = 1}^m a_i^{(j)}q_{i,l}^{(j)}
\Vert_M\leq C\sqrt{\sum_{j = 1}^m a_i^{(j)2}}\times (1 + l)^{2-\alpha}
\end{aligned}
\end{equation}
Define $r_{i,l}^{(j)} = \mathbf{E}\epsilon_i^{(j)}|\mathcal{H}_{i, (l,l)}^{(j), (l,l - 1)} -
\mathbf{E}\epsilon_i^{(j)}|\mathcal{H}_{i, (l,l)}^{(j), (l - 1,l - 1)}$, then
\begin{equation}
\begin{aligned}
r_{i,l}^{(j)}  = \mathbf{E}r_{i,l}^{(j)}|\mathcal{H}_{i, (0, 0)}^{(j), (l, l - 1)} +
\sum_{k = 1}^{l}\left(\mathbf{E}r_{i,l}^{(j)}|\mathcal{H}_{i, (k, k - 1)}^{(j), (l, l - 1)} - \mathbf{E}r_{i,l}^{(j)}|\mathcal{H}_{i, (k - 1, k - 1)}^{(j), (l, l - 1)}\right)\\
+ \sum_{k = 1}^l\left(\mathbf{E}r_{i,l}^{(j)}|\mathcal{H}_{i, (k, k)}^{(j), (l, l - 1)} - \mathbf{E}r_{i,l}^{(j)}|\mathcal{H}_{i, (k, k - 1)}^{(j), (l, l - 1)}\right)\\
\Rightarrow \Vert\sum_{i = 1}^n a_i^{(j)}r_{i,l}^{(j)}\Vert_M\leq
\Vert\sum_{i = 1}^n a_i^{(j)}\mathbf{E}r_{i,l}^{(j)}|\mathcal{H}_{i, (0, 0)}^{(j), (l, l - 1)}\Vert_M\\
+ \sum_{k = 1}^l \Vert\sum_{i = 1}^n a_i^{(j)}(\mathbf{E}r_{i,l}^{(j)}|\mathcal{H}_{i, (k, k - 1)}^{(j), (l, l - 1)} - \mathbf{E}r_{i,l}^{(j)}|\mathcal{H}_{i, (k - 1, k - 1)}^{(j), (l, l - 1)})\Vert_M\\
+ \sum_{k = 1}^l \Vert\sum_{i = 1}^n a_i^{(j)}(\mathbf{E}r_{i,l}^{(j)}|\mathcal{H}_{i, (k, k)}^{(j), (l, l - 1)} - \mathbf{E}r_{i,l}^{(j)}|\mathcal{H}_{i, (k, k - 1)}^{(j), (l, l - 1)})\Vert_M
\end{aligned}
\end{equation}
and
\begin{align*}
\mathbf{E}r_{i,l}^{(j)}|\mathcal{H}_{i, (k, k - 1)}^{(j), (l, l - 1)} - \mathbf{E}r_{i,l}^{(j)}|\mathcal{H}_{i, (k - 1, k - 1)}^{(j), (l, l - 1)}
= \mathbf{E}\epsilon_i^{(j)}|\mathcal{H}_{i, (k, k - 1)}^{(j), (l, l - 1)} - \mathbf{E}\epsilon_i^{(j)}|\mathcal{H}_{i, (k, k - 1)}^{(j), (l - 1, l - 1)}\\
- \mathbf{E}\epsilon_i^{(j)}|\mathcal{H}_{i, (k - 1, k - 1)}^{(j), (l, l - 1)}
+ \mathbf{E}\epsilon_i^{(j)}|\mathcal{H}_{i, (k - 1, k - 1)}^{(j), (l  - 1, l - 1)}
\end{align*}

From lemma \ref{lemma.delta_moment},
\begin{align*}
\mathbf{E}r_{i,l}^{(j)}|\mathcal{H}_{i, (0, 0)}^{(j), (l, l - 1)}
= \mathbf{E}\epsilon_i^{(j)}|\mathcal{H}_{i, (0, 0)}^{(j), (l, l - 1)} - \mathbf{E}\epsilon_i^{(j)}|\mathcal{H}_{i, (0,0)}^{(j), (l - 1,l - 1)}\\
\Rightarrow\Vert\sum_{i = 1}^n a_i^{(j)}\mathbf{E}r_{i,l}^{(j)}|\mathcal{H}_{i, (0, 0)}^{(j), (l, l - 1)}\Vert_M\leq C\delta_0^{(-l)}\sqrt{\sum_{i = 1}^n a_i^{(j)2}}
\end{align*}
Define $R_{s, l}^{(j)} = \sum_{i = n - s}^n a_i^{(j)}(\mathbf{E}r_{i,l}^{(j)}|\mathcal{H}_{i, (k, k - 1)}^{(j), (l, l - 1)} - \mathbf{E}r_{i,l}^{(j)}|\mathcal{H}_{i, (k - 1, k - 1)}^{(j), (l, l - 1)})$ for $s = 0,\cdots, n - 1$, then $R_{s,l}^{(j)}$ is
$\mathcal{Z}_{(n - s - k, n + k - 1)}^{(j-l, j + l - 1)}$ measurable, and
\begin{align*}
\mathbf{E}\left(\mathbf{E}r_{n - s - 1,l}^{(j)}|\mathcal{H}_{n - s - 1, (k, k - 1)}^{(j), (l, l - 1)}\right)|\mathcal{Z}_{(n - s - k, n + k - 1)}^{(j-l, j + l - 1)} = \mathbf{E}r_{n - s - 1,l}^{(j)}|\mathcal{H}_{n-s-1, (k - 1, k - 1)}^{(j), (l, l - 1)}
\end{align*}
Therefore, $R_{s, l}^{(j)}$ form a martingale, and
\begin{equation}
\begin{aligned}
\Vert\sum_{i = 1}^n a_i^{(j)}(\mathbf{E}r_{i,l}^{(j)}|\mathcal{H}_{i, (k, k - 1)}^{(j), (l, l - 1)} - \mathbf{E}r_{i,l}^{(j)}|\mathcal{H}_{i, (k - 1, k - 1)}^{(j), (l, l - 1)})\Vert_M\\
\leq C\sqrt{\sum_{i = 1}^n a_i^{(j)2}}\times \max_{i = 1,\cdots, n, j = 1,\cdots, m}\Vert\mathbf{E}r_{i,l}^{(j)}|\mathcal{H}_{i, (k, k - 1)}^{(j), (l, l - 1)} - \mathbf{E}r_{i,l}^{(j)}|\mathcal{H}_{i, (k - 1, k - 1)}^{(j), (l, l - 1)}\Vert_M\\
\leq 2C\sqrt{\sum_{i = 1}^n a_i^{(j)2}}\times \sum_{s = -k}^k\delta_s^{(-l)}\leq C^\prime\sqrt{\sum_{i = 1}^n a_i^{(j)2}}\times(1 + l)^{1-\alpha}
\end{aligned}
\end{equation}
Similarly define $\widetilde{R}_{s,l}^{(j)} = \sum_{i = 1}^s a_i^{(j)}(\mathbf{E}r_{i,l}^{(j)}|\mathcal{H}_{i, (k, k)}^{(j), (l, l - 1)} - \mathbf{E}r_{i,l}^{(j)}|\mathcal{H}_{i, (k, k - 1)}^{(j), (l, l - 1)})$, then $\widetilde{R}_{s,l}^{(j)}$ is $\mathcal{Z}_{(1 - k, s+k)}^{(j-l, j + l - 1)}$ measurable, and
\begin{align*}
\mathbf{E}\left(\mathbf{E}r_{s+1,l}^{(j)}|\mathcal{H}_{s+1, (k, k)}^{(j), (l, l - 1)}\right)|\mathcal{Z}_{(1 - k, s+k)}^{(j-l, j + l - 1)}
= \mathbf{E}r_{s+1,l}^{(j)}|\mathcal{H}_{s+1, (k, k - 1)}^{(j), (l, l - 1)}
\end{align*}
For
\begin{align*}
\mathbf{E}r_{i,l}^{(j)}|\mathcal{H}_{i, (k, k)}^{(j), (l, l - 1)} - \mathbf{E}r_{i,l}^{(j)}|\mathcal{H}_{i, (k, k - 1)}^{(j), (l, l - 1)}
= \mathbf{E}\epsilon_i^{(j)}|\mathcal{H}_{i, (k, k)}^{(j), (l, l - 1)} - \mathbf{E}\epsilon_i^{(j)}|\mathcal{H}_{i, (k, k)}^{(j), (l - 1, l - 1)}\\
- \mathbf{E}\epsilon_i^{(j)}|\mathcal{H}_{i, (k, k - 1)}^{(j), (l, l - 1)} + \mathbf{E}\epsilon_i^{(j)}|\mathcal{H}_{i, (k, k - 1)}^{(j), (l - 1, l - 1)}
\end{align*}
we have
\begin{equation}
\begin{aligned}
\Vert\sum_{i = 1}^n a_i^{(j)}(\mathbf{E}r_{i,l}^{(j)}|\mathcal{H}_{i, (k, k)}^{(j), (l, l - 1)} - \mathbf{E}r_{i,l}^{(j)}|\mathcal{H}_{i, (k, k - 1)}^{(j), (l, l - 1)})\Vert_M\\
\leq C\sqrt{\sum_{i = 1}^na_i^{(j)2}}\times \max_{i = 1,\cdots, n, j = 1,\cdots, m}\Vert\mathbf{E}r_{i,l}^{(j)}|\mathcal{H}_{i, (k, k)}^{(j), (l, l - 1)} - \mathbf{E}r_{i,l}^{(j)}|\mathcal{H}_{i, (k, k - 1)}^{(j), (l, l - 1)}\Vert_M\\
\leq 2C\sqrt{\sum_{i = 1}^na_i^{(j)2}}\times\sum_{s = -k}^k \delta_{s}^{(-l)}\leq C^\prime\sqrt{\sum_{i = 1}^na_i^{(j)2}}\times (1 + l)^{1-\alpha}
\end{aligned}
\end{equation}
and
\begin{equation}
\Vert\sum_{i = 1}^n a_i^{(j)}(\mathbf{E}\epsilon_i^{(j)}|\mathcal{H}_{i, (l,l)}^{(j), (l,l - 1)} -
\mathbf{E}\epsilon_i^{(j)}|\mathcal{H}_{i, (l,l)}^{(j), (l - 1,l - 1)})\Vert_M\leq  C\sqrt{\sum_{i = 1}^n a_i^{(j)2}}\times (1 + l)^{2-\alpha}
\end{equation}
Finally, define $t_{i,l}^{(j)} = \mathbf{E}\epsilon_i^{(j)}|\mathcal{H}_{i, (l,l)}^{(j), (l,l)} -
\mathbf{E}\epsilon_i^{(j)}|\mathcal{H}_{i, (l,l)}^{(j), (l,l - 1)}$, then
\begin{align*}
\Vert\sum_{i = 1}^n a_i^{(j)}t_{i,l}^{(j)}\Vert_M\leq \Vert\sum_{i = 1}^n a_i^{(j)}\mathbf{E}t_{i,l}^{(j)}|\mathcal{H}_{i, (0,0)}^{(j), (l,l)}\Vert_M\\
+\sum_{k = 1}^l\Vert\sum_{i = 1}^n a_i^{(j)}(\mathbf{E}t_{i,l}^{(j)}|\mathcal{H}_{i, (k,k - 1)}^{(j), (l,l)} - \mathbf{E}t_{i,l}^{(j)}|\mathcal{H}_{i, (k - 1,k - 1)}^{(j), (l,l)})\Vert_M\\
\sum_{k = 1}^l\Vert\sum_{i = 1}^n a_i^{(j)}(\mathbf{E}t_{i,l}^{(j)}|\mathcal{H}_{i, (k,k)}^{(j), (l,l)} - \mathbf{E}t_{i,l}^{(j)}|\mathcal{H}_{i, (k,k - 1)}^{(j), (l,l)})\Vert_M
\end{align*}
Notice that
\begin{align*}
\mathbf{E}t_{i,l}^{(j)}|\mathcal{H}_{i, (0,0)}^{(j), (l,l)} = \mathbf{E}\epsilon_i^{(j)}|\mathcal{H}_{i, (0,0)}^{(j), (l,l)} - \mathbf{E}\epsilon_i^{(j)}|\mathcal{H}_{i, (0,0)}^{(j), (l,l - 1)}\\
\Rightarrow \Vert\sum_{i = 1}^n a_i^{(j)}\mathbf{E}t_{i,l}^{(j)}|\mathcal{H}_{i, (0,0)}^{(j), (l,l)}\Vert_M\leq C\sqrt{\sum_{i = 1}^n a_i^{(j)2}}\times \delta_0^{(l)}
\end{align*}
Define $T_{s,l}^{(j)} = \sum_{i = n-s}^n a_i^{(j)}(\mathbf{E}t_{i,l}^{(j)}|\mathcal{H}_{i, (k,k - 1)}^{(j), (l,l)} - \mathbf{E}t_{i,l}^{(j)}|\mathcal{H}_{i, (k - 1,k - 1)}^{(j), (l,l)})$, then $T_{s,l}^{(j)}$ is $\mathcal{Z}_{(n - s - k, n + k - 1)}^{(j-l, j + l)}$ measurable, and
\begin{align*}
\mathbf{E}\left(\mathbf{E}t_{n-s-1,l}^{(j)}|\mathcal{H}_{n-s-1, (k,k - 1)}^{(j), (l,l)}\right)|\mathcal{Z}_{(n - s - k, n + k - 1)}^{(j-l, j + l)}
= \mathbf{E}t_{n-s-1,l}^{(j)}|\mathcal{H}_{n-s-1, (k - 1,k - 1)}^{(j), (l,l)}
\end{align*}
For
\begin{align*}
\mathbf{E}t_{i,l}^{(j)}|\mathcal{H}_{i, (k,k - 1)}^{(j), (l,l)} - \mathbf{E}t_{i,l}^{(j)}|\mathcal{H}_{i, (k - 1,k - 1)}^{(j), (l,l)}
= \mathbf{E}\epsilon_i^{(j)}|\mathcal{H}_{i, (k, k - 1)}^{(j), (l, l)} - \mathbf{E}\epsilon_i^{(j)}|\mathcal{H}_{i, (k, k - 1)}^{(j), (l, l - 1)}\\
- \mathbf{E}\epsilon_i^{(j)}|\mathcal{H}_{i, (k - 1, k - 1)}^{(j), (l, l )} + \mathbf{E}\epsilon_i^{(j)}|\mathcal{H}_{i, (k - 1, k - 1)}^{(j), (l, l - 1)}
\end{align*}
we have
\begin{align*}
\Vert\sum_{i = 1}^n a_i^{(j)}(\mathbf{E}t_{i,l}^{(j)}|\mathcal{H}_{i, (k,k - 1)}^{(j), (l,l)} - \mathbf{E}t_{i,l}^{(j)}|\mathcal{H}_{i, (k - 1,k - 1)}^{(j), (l,l)})\Vert_M\\
\leq C\sqrt{\sum_{i = 1}^n a_i^{(j)2}}\times \max_{i = 1,\cdots, n, j = 1,\cdots, m}\Vert\mathbf{E}t_{i,l}^{(j)}|\mathcal{H}_{i, (k,k - 1)}^{(j), (l,l)} - \mathbf{E}t_{i,l}^{(j)}|\mathcal{H}_{i, (k - 1,k - 1)}^{(j), (l,l)}\Vert_M\\
\leq 2C\sqrt{\sum_{i = 1}^n a_i^{(j)2}}\times \sum_{s = - k}^k \delta_s^{(l)}\leq C^\prime \sqrt{\sum_{i = 1}^n a_i^{(j)2}}\times (1 + l)^{1-\alpha}
\end{align*}
Define $\widetilde{T}_{s,l}^{(j)} = \sum_{i = 1}^s a_i^{(j)}(\mathbf{E}t_{i,l}^{(j)}|\mathcal{H}_{i, (k,k)}^{(j), (l,l)} - \mathbf{E}t_{i,l}^{(j)}|\mathcal{H}_{i, (k,k - 1)}^{(j), (l,l)})$, then $\widetilde{T}_{s,l}^{(j)}$ is $\mathcal{Z}_{(1-k, s + k)}^{(j-l, j + l)}$ measurable, and
\begin{align*}
\mathbf{E}\left(\mathbf{E}t_{s+1,l}^{(j)}|\mathcal{H}_{s+1, (k,k)}^{(j), (l,l)}\right)|\mathcal{Z}_{(1-k, s + k)}^{(j-l, j + l)}
= \mathbf{E}t_{s+1,l}^{(j)}|\mathcal{H}_{s+1, (k,k - 1)}^{(j), (l,l)}
\end{align*}
For
\begin{align*}
\mathbf{E}t_{i,l}^{(j)}|\mathcal{H}_{i, (k,k)}^{(j), (l,l)} - \mathbf{E}t_{i,l}^{(j)}|\mathcal{H}_{i, (k,k - 1)}^{(j), (l,l)}
= \mathbf{E}\epsilon_i^{(j)}|\mathcal{H}_{i, (k,k)}^{(j), (l,l)} - \mathbf{E}\epsilon_i^{(j)}|\mathcal{H}_{i, (k,k)}^{(j), (l,l - 1)}\\
- \mathbf{E}\epsilon_i^{(j)}|\mathcal{H}_{i, (k,k - 1)}^{(j), (l,l)} + \mathbf{E}\epsilon_i^{(j)}|\mathcal{H}_{i, (k,k - 1)}^{(j), (l,l - 1)}
\end{align*}
we have
\begin{align*}
\Vert\sum_{i = 1}^n a_i^{(j)}(\mathbf{E}t_{i,l}^{(j)}|\mathcal{H}_{i, (k,k)}^{(j), (l,l)} - \mathbf{E}t_{i,l}^{(j)}|\mathcal{H}_{i, (k,k - 1)}^{(j), (l,l)})\Vert_M\\
\leq C\sqrt{\sum_{i = 1}^n a_i^{(j)2}}\times \max_{i = 1,\cdots, n, j = 1,\cdots, m}\Vert\mathbf{E}t_{i,l}^{(j)}|\mathcal{H}_{i, (k,k)}^{(j), (l,l)} - \mathbf{E}t_{i,l}^{(j)}|\mathcal{H}_{i, (k,k - 1)}^{(j), (l,l)}\Vert_M\\
\leq 2C\sqrt{\sum_{i = 1}^n a_i^{(j)2}}\times\sum_{s = -k}^k \delta_{s}^{(l)}\leq C^\prime \sqrt{\sum_{i = 1}^n a_i^{(j)2}}\times(1 + l)^{1-\alpha}
\end{align*}
and
\begin{equation}
\Vert\sum_{i = 1}^n a_i^{(j)}(\mathbf{E}\epsilon_i^{(j)}|\mathcal{H}_{i, (l,l)}^{(j), (l,l)} -
\mathbf{E}\epsilon_i^{(j)}|\mathcal{H}_{i, (l,l)}^{(j), (l,l - 1)})\Vert_M\leq C(1 + l)^{2-\alpha}
\end{equation}
Then we prove \eqref{eq.sub_difference}.
\end{proof}

\begin{proof}[proof of lemma \ref{lemma.consistent_linear_combination}]
We prove the following result: suppose $\epsilon_i^{(j)}, i = 1,\cdots, n, j = 1,\cdots, m$ are $(M,\alpha)-$short range dependent random variables with $M>4, \alpha > 3$; and $a_i^{(j)}, i = 1,\cdots, n, j = 1,\cdots, m$ are real numbers. Then there exists a constant $C$ independent of $a_i^{(j)}$ such that for $\forall l\geq 0$,
\begin{equation}
\begin{aligned}
\Vert\sum_{i = 1}^n\sum_{j = 1}^m a_i^{(j)}\epsilon_i^{(j)}\Vert_M\leq C\sqrt{\sum_{i = 1}^n\sum_{j = 1}^m a_i^{(j)2}}\\
\Vert\sum_{i = 1}^n\sum_{j = 1}^m a_i^{(j)}(\epsilon_i^{(j)} - \mathbf{E}\epsilon_i^{(j)}|\mathcal{H}_{i, (l,l)}^{(j), (l,l)})\Vert_M
\leq C\sqrt{\sum_{i = 1}^n\sum_{j = 1}^m a_i^{(j)2}}\times (1 + l)^{3-\alpha}\\
\text{and }\Vert\sum_{i = 1}^n\sum_{j = 1}^m a_i^{(j)}\mathbf{E}\epsilon_i^{(j)}|\mathcal{H}_{i, (l,l)}^{(j), (l,l)}\Vert_M\leq C\sqrt{\sum_{i = 1}^n\sum_{j = 1}^m a_i^{(j)2}}.
\end{aligned}
\label{eq.three_lemma}
\end{equation}

From corollary C.9 in \cite{MR2001996}, for
$\mathcal{H}_{i, (l - 1, l - 1)}^{(j), (l - 1, l - 1)}\subset \mathcal{H}_{i, (l, l - 1)}^{(j), (l - 1, l - 1)}\subset
\mathcal{H}_{i, (l, l)}^{(j), (l - 1, l - 1)}\subset \mathcal{H}_{i, (l, l)}^{(j), (l, l - 1)}\subset \mathcal{H}_{i, (l, l)}^{(j), (l, l)}$, we have
\begin{align*}
\epsilon_i^{(j)} = \mathbf{E}\epsilon_i^{(j)}|\mathcal{H}_{i, (0,0)}^{(j), (0,0)}
+ \sum_{l = 1}^\infty((\mathbf{E}\epsilon_i^{(j)}|\mathcal{H}_{i, (l,l - 1)}^{(j), (l - 1,l - 1)} -
\mathbf{E}\epsilon_i^{(j)}|\mathcal{H}_{i, (l - 1,l - 1)}^{(j), (l - 1,l - 1)})\\
+ (\mathbf{E}\epsilon_i^{(j)}|\mathcal{H}_{i, (l,l)}^{(j), (l - 1,l - 1)} -
\mathbf{E}\epsilon_i^{(j)}|\mathcal{H}_{i, (l,l - 1)}^{(j), (l - 1,l - 1)})
+ (\mathbf{E}\epsilon_i^{(j)}|\mathcal{H}_{i, (l,l)}^{(j), (l,l - 1)} -
\mathbf{E}\epsilon_i^{(j)}|\mathcal{H}_{i, (l,l)}^{(j), (l - 1,l - 1)})\\
+ (\mathbf{E}\epsilon_i^{(j)}|\mathcal{H}_{i, (l,l)}^{(j), (l,l)} -
\mathbf{E}\epsilon_i^{(j)}|\mathcal{H}_{i, (l,l)}^{(j), (l,l - 1)})
)
\end{align*}
Correspondingly,
\begin{equation}
\begin{aligned}
\Vert\sum_{i = 1}^n\sum_{j = 1}^m a_i^{(j)}\epsilon_i^{(j)}\Vert_M\leq
\Vert\sum_{i = 1}^n\sum_{j = 1}^m a_i^{(j)}\mathbf{E}\epsilon_i^{(j)}|\mathcal{H}_{i, (0,0)}^{(j), (0,0)}\Vert_M\\
+\sum_{l = 1}^\infty\Vert\sum_{i = 1}^n\sum_{j = 1}^m a_i^{(j)}(\mathbf{E}\epsilon_i^{(j)}|\mathcal{H}_{i, (l,l - 1)}^{(j), (l - 1,l - 1)} -
\mathbf{E}\epsilon_i^{(j)}|\mathcal{H}_{i, (l - 1,l - 1)}^{(j), (l - 1,l - 1)})\Vert_M\\
+\sum_{l = 1}^\infty\Vert\sum_{i = 1}^n\sum_{j = 1}^m a_i^{(j)}(\mathbf{E}\epsilon_i^{(j)}|\mathcal{H}_{i, (l,l)}^{(j), (l - 1,l - 1)} -
\mathbf{E}\epsilon_i^{(j)}|\mathcal{H}_{i, (l,l - 1)}^{(j), (l - 1,l - 1)})\Vert_M\\
+\sum_{l = 1}^\infty\Vert\sum_{i = 1}^n\sum_{j = 1}^m a_i^{(j)}(\mathbf{E}\epsilon_i^{(j)}|\mathcal{H}_{i, (l,l)}^{(j), (l,l - 1)} -
\mathbf{E}\epsilon_i^{(j)}|\mathcal{H}_{i, (l,l)}^{(j), (l - 1,l - 1)})\Vert_M\\
+ \sum_{l = 1}^\infty\Vert\sum_{i = 1}^n\sum_{j = 1}^m a_i^{(j)}(\mathbf{E}\epsilon_i^{(j)}|\mathcal{H}_{i, (l,l)}^{(j), (l,l)} -
\mathbf{E}\epsilon_i^{(j)}|\mathcal{H}_{i, (l,l)}^{(j), (l,l - 1)})\Vert_M
\end{aligned}
\end{equation}
\begin{equation}
\begin{aligned}
\Vert\sum_{i = 1}^n\sum_{j = 1}^m a_i^{(j)}(\epsilon_i^{(j)} - \mathbf{E}\epsilon_i^{(j)}|\mathcal{H}_{i, (l,l)}^{(j), (l,l)})\Vert_M\\
\leq
\sum_{s = l+1}^\infty\Vert\sum_{i = 1}^n\sum_{j = 1}^m a_i^{(j)}(\mathbf{E}\epsilon_i^{(j)}|\mathcal{H}_{i, (s,s - 1)}^{(j), (s - 1,s - 1)} -
\mathbf{E}\epsilon_i^{(j)}|\mathcal{H}_{i, (s - 1,s - 1)}^{(j), (s - 1,s - 1)})\Vert_M\\
+\sum_{s = l + 1}^\infty\Vert\sum_{i = 1}^n\sum_{j = 1}^m a_i^{(j)}(\mathbf{E}\epsilon_i^{(j)}|\mathcal{H}_{i, (s,s)}^{(j), (s - 1,s - 1)} -
\mathbf{E}\epsilon_i^{(j)}|\mathcal{H}_{i, (s,s - 1)}^{(j), (s - 1,s - 1)})\Vert_M\\
+\sum_{s = l + 1}^\infty\Vert\sum_{i = 1}^n\sum_{j = 1}^m a_i^{(j)}(\mathbf{E}\epsilon_i^{(j)}|\mathcal{H}_{i, (s,s)}^{(j), (s,s - 1)} -
\mathbf{E}\epsilon_i^{(j)}|\mathcal{H}_{i, (s,s)}^{(j), (s - 1,s - 1)})\Vert_M\\
+ \sum_{s = l+1}^\infty\Vert\sum_{i = 1}^n\sum_{j = 1}^m a_i^{(j)}(\mathbf{E}\epsilon_i^{(j)}|\mathcal{H}_{i, (s,s)}^{(j), (s,s)} -
\mathbf{E}\epsilon_i^{(j)}|\mathcal{H}_{i, (s,s)}^{(j), (s,s - 1)})\Vert_M
\end{aligned}
\end{equation}
and
\begin{equation}
\begin{aligned}
\Vert\sum_{i = 1}^n\sum_{j = 1}^m a_i^{(j)}\mathbf{E}\epsilon_i^{(j)}|\mathcal{H}_{i, (l,l)}^{(j), (l,l)}\Vert_M
\leq
\Vert\sum_{i = 1}^n\sum_{j = 1}^m a_i^{(j)}\mathbf{E}\epsilon_i^{(j)}|\mathcal{H}_{i, (0,0)}^{(j), (0,0)}\Vert_M\\
+\sum_{s = 1}^l\Vert\sum_{i = 1}^n\sum_{j = 1}^m a_i^{(j)}(\mathbf{E}\epsilon_i^{(j)}|\mathcal{H}_{i, (s,s - 1)}^{(j), (s - 1,s - 1)} -
\mathbf{E}\epsilon_i^{(j)}|\mathcal{H}_{i, (s - 1,s - 1)}^{(j), (s - 1,s - 1)})\Vert_M\\
+\sum_{s = 1}^l\Vert\sum_{i = 1}^n\sum_{j = 1}^m a_i^{(j)}(\mathbf{E}\epsilon_i^{(j)}|\mathcal{H}_{i, (s,s)}^{(j), (s - 1,s - 1)} -
\mathbf{E}\epsilon_i^{(j)}|\mathcal{H}_{i, (s,s - 1)}^{(j), (s - 1,s - 1)})\Vert_M\\
+\sum_{s = 1}^l\Vert\sum_{i = 1}^n\sum_{j = 1}^m a_i^{(j)}(\mathbf{E}\epsilon_i^{(j)}|\mathcal{H}_{i, (s,s)}^{(j), (s,s - 1)} -
\mathbf{E}\epsilon_i^{(j)}|\mathcal{H}_{i, (s,s)}^{(j), (s - 1,s - 1)})\Vert_M\\
+ \sum_{s =1}^l\Vert\sum_{i = 1}^n\sum_{j = 1}^m a_i^{(j)}(\mathbf{E}\epsilon_i^{(j)}|\mathcal{H}_{i, (s,s)}^{(j), (s,s)} -
\mathbf{E}\epsilon_i^{(j)}|\mathcal{H}_{i, (s,s)}^{(j), (s,s - 1)})\Vert_M
\end{aligned}
\end{equation}
From \cite{MR0133849},
\begin{equation}
\begin{aligned}
\Vert\sum_{i = 1}^n\sum_{j = 1}^m a_i^{(j)}\mathbf{E}\epsilon_i^{(j)}|\mathcal{H}_{i, (0,0)}^{(j), (0,0)}\Vert_M\leq
C\sqrt{\sum_{i = 1}^n\sum_{j = 1}^m a_i^{(j)2}}\times\max_{i = 1,\cdots, n, j = 1,\cdots, m}\Vert\mathbf{E}\epsilon_i^{(j)}|\mathcal{H}_{i, (0,0)}^{(j), (0,0)}\Vert_M\\
\leq C^\prime\sqrt{\sum_{i = 1}^n\sum_{j = 1}^m a_i^{(j)2}}
\end{aligned}
\end{equation}
with a constant $C^\prime$. From lemma \ref{lemma.calculate_martingale_difference} and \ref{lemma.second_dimension_difference},
\begin{equation}
\begin{aligned}
\Vert\sum_{i = 1}^n\sum_{j = 1}^m a_i^{(j)}(\mathbf{E}\epsilon_i^{(j)}|\mathcal{H}_{i, (l,l - 1)}^{(j), (l - 1,l - 1)} -
\mathbf{E}\epsilon_i^{(j)}|\mathcal{H}_{i, (l - 1,l - 1)}^{(j), (l - 1,l - 1)})\Vert_M\\
\leq C\sqrt{\sum_{i = 1}^n\Vert\sum_{j = 1}^m a_i^{(j)}(\mathbf{E}\epsilon_i^{(j)}|\mathcal{H}_{i, (l,l - 1)}^{(j), (l - 1,l - 1)} -
\mathbf{E}\epsilon_i^{(j)}|\mathcal{H}_{i, (l - 1,l - 1)}^{(j), (l - 1,l - 1)})\Vert_M^2}\\
\leq C^\prime\sqrt{\sum_{i = 1}^n\sum_{j = 1}^m a_i^{(j)2}}\times (1 + l)^{2-\alpha}
\end{aligned}
\end{equation}
With a similar deduction, we have
\begin{equation}
\begin{aligned}
\Vert\sum_{i = 1}^n\sum_{j = 1}^m a_i^{(j)}\epsilon_i^{(j)}\Vert_M,\ \Vert\sum_{i = 1}^n\sum_{j = 1}^m a_i^{(j)}\mathbf{E}\epsilon_i^{(j)}|\mathcal{H}_{i, (l,l)}^{(j), (l,l)}\Vert_M\\
\leq C\sqrt{\sum_{i = 1}^n\sum_{j = 1}^m a_i^{(j)2}}
\times\left(1 + \sum_{l = 1}^\infty (1 + l)^{2-\alpha}\right)\leq C^\prime \sqrt{\sum_{i = 1}^n\sum_{j = 1}^m a_i^{(j)2}}\\
\text{and } \Vert\sum_{i = 1}^n\sum_{j = 1}^m a_i^{(j)}(\epsilon_i^{(j)} - \mathbf{E}\epsilon_i^{(j)}|\mathcal{H}_{i, (l,l)}^{(j), (l,l)})\Vert_M
\leq C\sum_{s = l+1}^\infty\sqrt{\sum_{i = 1}^n\sum_{j = 1}^m a_i^{(j)2}}\times (1 + l)^{2-\alpha}\\
\leq C^\prime \sqrt{\sum_{i = 1}^n\sum_{j = 1}^m a_i^{(j)2}}\times (1 + l)^{3-\alpha}
\end{aligned}
\end{equation}
and we prove eq.\eqref{eq.three_lemma}.
\end{proof}

We can derive a corollary about the covariances of $\epsilon_i^{(j)}$ based on lemma \ref{lemma.consistent_linear_combination}.
\begin{corollary}
\label{corollary.covariance}
Suppose $\epsilon_i^{(j)}, i = 1,\cdots, n,\ j = 1,\cdots, m$ are $(M,\alpha)-$short range dependent random variables(see definition \ref{definition.Ma}) with $M > 4$ and $\alpha>3$. Then
there exists a constant $C>0$ such that
\begin{equation}
\vert\mathbf{E}\epsilon_{i_1}^{(j_1)}\epsilon_{i_2}^{(j_2)}\vert\leq C\times (1 + \max(\vert i_1 - i_2\vert, \vert j_1 - j_2\vert))^{3-\alpha}
\end{equation}
for any $i = 1,\cdots, n$ and $j = 1,\cdots, m$.
\end{corollary}
\begin{proof}
Suppose the integer $l\geq 2$, define $k = \lfloor\frac{l}{2}\rfloor$, here $\lfloor x\rfloor$ is the largest integer that is smaller than or equal to $x$. Then
\begin{align*}
\mathbf{E}\epsilon_{i_1}^{(j_1)}\epsilon_{i_2}^{(j_1 + l)} = \mathbf{E}\left(\mathbf{E}\epsilon_{i_1}^{(j_1)}|\mathcal{H}_{i_1, (k,k)}^{(j_1), (k,k)}\times\left(\epsilon_{i_2}^{(j_1 + l)} - \mathbf{E}\epsilon_{i_2}^{(j_1 + l)}|\mathcal{H}_{i_2, (k - 1,k - 1)}^{(j_1 + l), (k - 1,k - 1)}\right)\right)\\
+ \mathbf{E}\left(\left(\epsilon_{i_1}^{(j_1)} - \mathbf{E}\epsilon_{i_1}^{(j_1)}|\mathcal{H}_{i_1, (k,k)}^{(j_1), (k, k)}\right)\times \mathbf{E}\epsilon_{i_2}^{(j_1 + l)}|\mathcal{H}_{i_2, (k - 1,k - 1)}^{(j_1 + l), (k - 1,k - 1)}\right)\\
+ \mathbf{E}(
\left(\epsilon_{i_1}^{(j_1)} - \mathbf{E}\epsilon_{i_1}^{(j_1)}|\mathcal{H}_{i_1, (k,k)}^{(j_1), (k, k)}\right)\times \left(\epsilon_{i_2}^{(j_1 + l)} - \mathbf{E}\epsilon_{i_2}^{(j_1 + l)}|\mathcal{H}_{i_2, (k - 1,k - 1)}^{(j_1 + l), (k - 1,k - 1)}\right)
)
\end{align*}
According to lemma \ref{lemma.consistent_linear_combination}, we have $\vert\mathbf{E}\epsilon_{i_1}^{(j_1)}\epsilon_{i_2}^{(j_1 + l)}\vert\leq C(1 + l)^{3-\alpha}$ with a constant $C$. Similarly,
\begin{align*}
\mathbf{E}\epsilon_{i_1}^{(j_1)}\epsilon_{i_1 + l}^{(j_2)} = \mathbf{E}\left(\mathbf{E}\epsilon_{i_1}^{(j_1)}|\mathcal{H}_{i_1, (k,k)}^{(j_1), (k,k)}\times\left(\epsilon_{i_1+l}^{(j_2)} - \mathbf{E}\epsilon_{i_1+l}^{(j_2)}|\mathcal{H}_{i_1+l, (k - 1,k - 1)}^{(j_2), (k - 1,k - 1)}\right)\right)\\
+ \mathbf{E}\left(\left(\epsilon_{i_1}^{(j_1)} - \mathbf{E}\epsilon_{i_1}^{(j_1)}|\mathcal{H}_{i_1, (k,k)}^{(j_1), (k, k)}\right)\times \mathbf{E}\epsilon_{i_1+l}^{(j_2)}|\mathcal{H}_{i_1+l, (k - 1,k - 1)}^{(j_2), (k - 1,k - 1)}\right)\\
+ \mathbf{E}(
\left(\epsilon_{i_1}^{(j_1)} - \mathbf{E}\epsilon_{i_1}^{(j_1)}|\mathcal{H}_{i_1, (k,k)}^{(j_1), (k, k)}\right)\times \left(\epsilon_{i_1+l}^{(j_2)} - \mathbf{E}\epsilon_{i_1+l}^{(j_2)}|\mathcal{H}_{i_1+l, (k - 1,k - 1)}^{(j_2), (k - 1,k - 1)}\right)
)
\end{align*}
so $\vert \mathbf{E}\epsilon_{i_1}^{(j_1)}\epsilon_{i_1 + l}^{(j_2)}\vert\leq C(1 + l)^{3-\alpha}$ with a constant $C$. Therefore,
\begin{equation}
\vert\mathbf{E}\epsilon_{i_1}^{(j_1)}\epsilon_{i_2}^{(j_2)}\vert\leq C\times (1 + \max(\vert i_1 - i_2\vert, \vert j_1 - j_2\vert))^{3-\alpha}
\end{equation}
\end{proof}

\begin{proof}[proof of theorem \ref{theorem.Gaussian}]
Define
$$
c_{i, v}^{(j)} = a_{i,v}^{(j)} / \sqrt{\sum_{i = 1}^n\sum_{j = 1}^m a_{i, v}^{(j)2}}\ \text{and } T_v = \frac{R_v}{\sqrt{\sum_{i = 1}^n\sum_{j = 1}^m a_{i, v}^{(j)2}}}
$$
then $\vert c_{i, v}^{(j)}\vert\leq \lambda_{i,v}\tau_{j,v}$ and $\sum_{i = 1}^n\sum_{j = 1}^m c_{i, v}^{(j)2} = 1$. Moreover,
$ T_v = \sum_{i = 1}^n\sum_{j = 1}^m c_{i, v}^{(j)}\epsilon_i^{(j)}$. Notice that
\begin{align*}
Var(\xi_v) = \sum_{i_1 = 1}^n\sum_{j_1 = 1}^m\sum_{i_2 = 1}^n\sum_{j_2  = 1}^m c_{i_1, v}^{(j_2)}c_{i_2,v }^{(j_2)}(\mathbf{E}\epsilon_{i_1}^{(j_1)}\epsilon_{i_2}^{(j_2)})\geq c_0\\
\text{and } \Vert R_v\Vert_2\leq \Vert R_v\Vert_M\leq C\sqrt{\sum_{i = 1}^n\sum_{j = 1}^m a_{i, v}^{(j)2}}
\end{align*}
so from eq.\eqref{eq.prob_to_h}, for any given $\tau,\psi > 0$,
\begin{equation}
\begin{aligned}
\sup_{x\in\mathbf{R}}\vert
Prob\left(\max_{v = 1,\cdots, V}\vert T_v\vert\leq x\right)
- Prob\left(\max_{v = 1,\cdots, V}\vert\xi_v\vert\leq x\right)\vert\\
\leq Ct(1 + \sqrt{\log(V)} + \sqrt{\vert\log(t)\vert})\\
+ \sup_{x\in\mathbf{R}}\vert \mathbf{E}h_{\tau, \psi, x}\left(T_1, \cdots, T_V\right) -
\mathbf{E}h_{\tau, \psi, x}(\xi_1,\cdots,\xi_V)\vert
\end{aligned}
\label{eq.difference_h0}
\end{equation}
here $t = \frac{1}{\psi} + \frac{\log(2V)}{\tau}$. For any given integer $s>0$ and $s/\mathcal{K}\to 0$, from lemma \ref{lemma.consistent_linear_combination},
\begin{align*}
\Vert\sum_{i = 1}^n\sum_{j = 1}^m c_{i, v}^{(j)}(\epsilon_i^{(j)} - \mathbf{E}\epsilon_i^{(j)}|\mathcal{H}_{i, (s,s)}^{(j), (s, s)})\Vert_M
\leq C(1 + s)^{3-\alpha}\\
\Rightarrow \Vert\ \max_{v = 1,\cdots, V}\vert\sum_{i = 1}^n\sum_{j = 1}^m c_{i, v}^{(j)}(\epsilon_i^{(j)} - \mathbf{E}\epsilon_i^{(j)}|\mathcal{H}_{i, (s,s)}^{(j), (s, s)})\vert\ \Vert_M\leq CV^{1/M}\times (1 + s)^{3-\alpha}
\end{align*}
Therefore,
\begin{equation}
\begin{aligned}
\sup_{x\in\mathbf{R}}\vert
\mathbf{E}h_{\tau, \psi, x}\left(T_1, \cdots, T_V\right)\\
- \mathbf{E}h_{\tau,\psi,x}\left(\sum_{i = 1}^n\sum_{j = 1}^m c_{i, 1}^{(j)}\mathbf{E}\epsilon_i^{(j)}|\mathcal{H}_{i, (s,s)}^{(j), (s, s)},\cdots, \sum_{i = 1}^n\sum_{j = 1}^m c_{i, V}^{(j)}\mathbf{E}\epsilon_i^{(j)}|\mathcal{H}_{i, (s,s)}^{(j), (s, s)}\right)
\vert\\
\leq C\psi V^{1/M}\times (1 + s)^{3-\alpha}
\end{aligned}
\label{eq.difference_h1}
\end{equation}
Notice that $T_v = \sum_{i = x_v - \mathcal{K}}^{x_v + \mathcal{K}}\sum_{j = y_v - \mathcal{K}}^{y_v + \mathcal{K}}c_{i, v}^{(j)}\epsilon_{i}^{(j)}$, for any $l > 2s$, define the block $S_{i, v}^{(j)}$ as
\begin{align*}
S_{i, v}^{(j)} = \sum_{p = (i - 1)\times (l + 2s) + x_v - \mathcal{K}}^{[x_v - \mathcal{K} + (i - 1)\times (l + 2s) + l - 1]\wedge[x_v + \mathcal{K}]}\sum_{q = (j  -  1)\times (l + 2s) + y_v -\mathcal{K}}^{[y_v -\mathcal{K} + (j - 1)\times (l + 2s) + l - 1]\wedge[y_v + \mathcal{K}]}c_{p, v}^{(q)}\mathbf{E}\epsilon_p^{(q)}|\mathcal{H}_{p, (s,s)}^{(q), (s, s)}
\end{align*}
then the vector $(S_{i, 1}^{(j)},\cdots, S_{i, V}^{(j)} )^T$ is independent of $(S_{p, 1}^{(q)},\cdots, S_{p, V}^{(q)} )^T$ if $(i,j)\neq (p,q)$. Define $n_1 = \lceil\frac{2\mathcal{K} + 1}{l + 2s}\rceil$, here $\lceil x\rceil$ denotes the smallest integer that is larger than or equal to $x$. Then from lemma \ref{lemma.consistent_linear_combination}
\begin{equation}
\begin{aligned}
\Vert
\sum_{i = x_v - \mathcal{K}}^{x_v + \mathcal{K}}\sum_{j = y_v - \mathcal{K}}^{y_v + \mathcal{K}} c_{i, v}^{(j)}\mathbf{E}\epsilon_i^{(j)}|\mathcal{H}_{i, (s,s)}^{(j), (s, s)} - \sum_{i = 1}^{n_1}\sum_{j = 1}^{n_1}
S_{i, v}^{(j)}
\Vert_M\\
\leq\Vert
\sum_{i = x_v - \mathcal{K}}^{x_v + \mathcal{K}}\sum_{j = y_v - \mathcal{K}}^{y_v + \mathcal{K}} c_{i, v}^{(j)}\mathbf{E}\epsilon_i^{(j)}|\mathcal{H}_{i, (s,s)}^{(j), (s, s)} -
\sum_{i = x_v - \mathcal{K}}^{x_v + \mathcal{K}}\sum_{p = 1}^{n_1}
\sum_{j = (p  -  1)\times (l + 2s) + y_v -\mathcal{K}}^{[y_v -\mathcal{K} + (p - 1)\times (l + 2s) + l - 1]\wedge[y_v + \mathcal{K}]}
c_{i, v}^{(j)}\mathbf{E}\epsilon_i^{(j)}|\mathcal{H}_{i, (s,s)}^{(j), (s, s)}
\Vert_M\\
+ \Vert\sum_{i = x_v - \mathcal{K}}^{x_v + \mathcal{K}}\sum_{p = 1}^{n_1}\sum_{j = (p  -  1)\times (l + 2s) + y_v -\mathcal{K}}^{[y_v -\mathcal{K} + (p - 1)\times (l + 2s) + l - 1]\wedge[y_v + \mathcal{K}]}
c_{i, v}^{(j)}\mathbf{E}\epsilon_i^{(j)}|\mathcal{H}_{i, (s,s)}^{(j), (s, s)}- \sum_{i = 1}^{n_1}\sum_{j = 1}^{m_1}
S_{i, v}^{(j)}\Vert_M\\
=\Vert
\sum_{i = x_v - \mathcal{K}}^{x_v + \mathcal{K}}\sum_{p = 1}^{n_1}\sum_{j = y_v - \mathcal{K} + (p  -  1)\times (l + 2s) + l}^{[y_v -\mathcal{K}- 1 + p\times (l + 2s)]\wedge[m]}c_{i, v}^{(j)}\mathbf{E}\epsilon_i^{(j)}|\mathcal{H}_{i, (s,s)}^{(j), (s, s)}
\Vert_M\\
+\Vert
\sum_{p = 1}^{n_1}\sum_{j = (p  -  1)\times (l + 2s) + y_v -\mathcal{K}}^{[y_v -\mathcal{K} + (p - 1)\times (l + 2s) + l - 1]\wedge[y_v + \mathcal{K}]}\sum_{q = 1}^{n_1}\sum_{i = x_v - \mathcal{K} + (q - 1)\times (l + 2s) + l}^{[x_v - \mathcal{K} - 1 + q\times(l + 2s)]\wedge [n]} c_{i, v}^{(j)}\mathbf{E}\epsilon_i^{(j)}|\mathcal{H}_{i, (s,s)}^{(j), (s, s)}
\Vert_M\\
\leq C\sqrt{\sum_{i = x_v - \mathcal{K}}^{x_v + \mathcal{K}}\sum_{p = 1}^{n_1}\sum_{j = y_v - \mathcal{K} + (p  -  1)\times (l + 2s) + l}^{[y_v -\mathcal{K}- 1 + p\times (l + 2s)]\wedge[m]}\lambda_{i,v}^2\tau_{j,v}^2}
+ C\sqrt{\sum_{j = y_v - \mathcal{K}}^{y_v + \mathcal{K}} \sum_{q = 1}^{n_1}\sum_{i = x_v - \mathcal{K} + (q - 1)\times (l + 2s) + l}^{[x_v - \mathcal{K} - 1 + q\times(l + 2s)]\wedge [n]} \lambda_{i,v}^2\tau_{j,v}^2}\\
\leq C^\prime\times \sqrt{\frac{s}{l}}
\end{aligned}
\label{eq.difference_S_C}
\end{equation}
In particular,
\begin{equation}
\begin{aligned}
\Vert\
\max_{v = 1,\cdots, V}\vert\sum_{i = 1}^n\sum_{j = 1}^m c_{i, v}^{(j)}\mathbf{E}\epsilon_i^{(j)}|\mathcal{H}_{i, (s,s)}^{(j), (s, s)} - \sum_{i = 1}^{n_1}\sum_{j = 1}^{m_1}
S_{i, v}^{(j)}\vert
\ \Vert_M
\leq CV^{1/M}\times \sqrt{\frac{s}{l}}\\
\text{and }\sup_{x\in\mathbf{R}}\vert\mathbf{E}h_{\tau,\psi,x}\left(\sum_{i = 1}^n\sum_{j = 1}^m c_{i, 1}^{(j)}\mathbf{E}\epsilon_i^{(j)}|\mathcal{H}_{i, (s,s)}^{(j), (s, s)},\cdots, \sum_{i = 1}^n\sum_{j = 1}^m c_{i, V}^{(j)}\mathbf{E}\epsilon_i^{(j)}|\mathcal{H}_{i, (s,s)}^{(j), (s, s)}\right)\\
- \mathbf{E}h_{\tau,\psi,x}\left(\sum_{i = 1}^{n_1}\sum_{j = 1}^{m_1}
S_{i, 1}^{(j)},\cdots, \sum_{i = 1}^{n_1}\sum_{j = 1}^{m_1}
S_{i, V}^{(j)}\right)\vert\\
\leq C\psi\times V^{1/M}\times\sqrt{\frac{s}{l}}
\end{aligned}
\label{eq.difference_h2}
\end{equation}
Define $P_{k,v} = S_{i, v}^{(j)}$ with $i = 1+\lfloor\frac{k - 1}{n_1}\rfloor$ and $j = k - (i - 1)\times n_1$, here $\lfloor x\rfloor$ denotes the largest integer that is smaller than, or equal to $x$. Define joint normal random vectors $(P_{k,1}^*,\cdots, P_{k,V}^*)^T$ such that $(P_{k_1,1}^*,\cdots, P_{k_1,V}^*)^T$ is independent of $(P_{k_2,1}^*,\cdots, P_{k_2,V}^*)^T$ for $k_1\neq k_2$, $P_{k,i}^*$ is independent of $P_{k_2, i_2}$ for any $i,i_2,k,k_2$, $\mathbf{E}P_{k,i}^* = 0$, and $\mathbf{E}P_{k,v_1}^*P_{k,v_2}^* = \mathbf{E}P_{k, v_1}P_{k, v_2}$. Define
$H_{k,v} = \sum_{j = 1}^{k-1}P_{j,v} + \sum_{j = k + 1}^{n_1^2}P_{j,v}^*$, then $H_{k,v} + P_{k,v} = H_{k+1, v} + P_{k+1, v}^*$. For
\begin{align*}
\mathbf{E}\left(h_{\tau,\psi,x}(H_{k,1} + P_{k,1},\cdots, H_{k,V} + P_{k,V}) - h_{\tau,\psi,x}(H_{k,1} + P_{k,1}^*,\cdots, H_{k,V} + P_{k,V}^*)\right)|H_{k,1},\cdots, H_{k,V}\\
= \sum_{v = 1}^V\partial_vh_{\tau,\psi,x}(H_{k,1},\cdots, H_{k,V})(\mathbf{E}P_{k,v} - \mathbf{E}P_{k,v}^*)\\
+ \frac{1}{2}\sum_{v_1 = 1}^V\sum_{v_2 = 1}^V\partial_{v_1}\partial_{v_2}h_{\tau,\psi,x}(H_{k,1},\cdots, H_{k,V})
(\mathbf{E}P_{k,v_1}P_{k,v_2} - \mathbf{E}P_{k,v_1}^*P^*_{k,v_2})\\
+ \frac{1}{6}\sum_{v_1 = 1}^V\sum_{v_2 = 1}^V\sum_{v_3 = 1}^V
\mathbf{E}\partial_{v_1}\partial_{v_2}\partial_{v_3}
h_{\tau,\psi,x}(\xi_1,\cdots, \xi_V)P_{k,v_1}P_{k,v_2}P_{k,v_3}
|H_{k,1},\cdots, H_{k,V}\\
- \frac{1}{6}\sum_{v_1 = 1}^V\sum_{v_2 = 1}^V\sum_{v_3 = 1}^V
\mathbf{E}\partial_{v_1}\partial_{v_2}\partial_{v_3}
h_{\tau,\psi,x}(\xi_1^*,\cdots, \xi_V^*)P_{k,v_1}^*P_{k,v_2}^*P_{k,v_3}^*
|H_{k,1},\cdots, H_{k,V}
\end{align*}
here $\xi_i,\xi_i^*$ are random variables. We have
\begin{equation}
\begin{aligned}
\vert\mathbf{E}\left(h_{\tau,\psi,x}(H_{k,1} + P_{k,1},\cdots, H_{k,V} + P_{k,V}) - h_{\tau,\psi,x}(H_{k,1} + P_{k,1}^*,\cdots, H_{k,V} + P_{k,V}^*)\right)\vert\\
\leq C(\psi^3 + \psi^2\tau + \psi\tau^2)\times (\Vert\ \max_{v = 1,\cdots,V}\vert P_{k,v}\vert\ \Vert_M^3 + \Vert\ \max_{v = 1,\cdots,V}\vert P_{k,v}^*\vert\ \Vert_M^3)\\
\leq C^\prime V^{3/M}(\psi^3 + \psi^2\tau + \psi\tau^2)\times\max_{v = 1,\cdots, V}\Vert P_{k,v}\Vert_M^3
\end{aligned}
\end{equation}
here we use the fact that $P_{k,v}^*$ has normal distribution and $\Vert P_{k,v}^*\Vert_M\leq C\Vert P_{k,v}^*\Vert_2 = C\Vert P_{k,v}\Vert_2$. From lemma \ref{lemma.consistent_linear_combination},
\begin{align*}
\Vert P_{k,v}\Vert_M = \Vert S_{i,v}^{(j)}\Vert_M\leq C\sqrt{\sum_{p = (i - 1)\times (l + 2s) + x_v - \mathcal{K}}^{[x_v - \mathcal{K} + (i - 1)\times (l + 2s) + l - 1]\wedge[x_v + \mathcal{K}]}\sum_{q = (j  -  1)\times (l + 2s) + y_v -\mathcal{K}}^{[y_v -\mathcal{K} + (j - 1)\times (l + 2s) + l - 1]\wedge[y_v + \mathcal{K}]}c_{p, v}^{(q)2}}\\
\leq C^\prime \frac{l}{\mathcal{K}}
\end{align*}
here $i = 1+\lfloor\frac{k - 1}{n_1}\rfloor$ and $j = k - (i - 1)\times n_1$. So
\begin{equation}
\begin{aligned}
\vert
\mathbf{E}h_{\tau,\psi,x}(\sum_{k = 1}^{n_1^2}P_{k,1},\cdots, \sum_{k = 1}^{n_1^2}P_{k,V})
-\mathbf{E}h_{\tau,\psi,x}(\sum_{k = 1}^{n_1^2}P_{k,1}^*,\cdots, \sum_{k = 1}^{n_1^2}P_{k,V}^*)
\vert\\
\leq \sum_{k = 1}^{n_1^2}\vert\mathbf{E}\left(h_{\tau,\psi,x}(H_{k,1} + P_{k,1},\cdots, H_{k,V} + P_{k,V}) - h_{\tau,\psi,x}(H_{k,1} + P_{k,1}^*,\cdots, H_{k,V} + P_{k,V}^*)\right)\vert\\
\leq C V^{3/M}(\psi^3 + \psi^2\tau + \psi\tau^2)\times \frac{l}{\mathcal{K}}
\end{aligned}
\label{eq.difference_h3}
\end{equation}
Notice that $\sum_{k = 1}^{n_1^2}P_{k,1} = \sum_{i = 1}^{n_1}\sum_{j = 1}^{n_1}S_{i,v}^{(j)}$. From eq.\eqref{eq.difference_S_C} and Cauchy inequality
\begin{align*}
\vert\sum_{i_1 = 1}^{n_1}\sum_{j_1 = 1}^{n_1}\sum_{i_2 = 1}^{n_1}\sum_{j_2 = 1}^{n_1} \mathbf{E}S_{i_1,v_1}^{(j_1)}S_{i_2,v_2}^{(j_2)}
-\sum_{i_1 = 1}^n\sum_{j_1 = 1}^m\sum_{i_2 = 1}^n\sum_{j_2 = 1}^m c_{i_1, v_1}^{(j_1)}c_{i_2, v_2}^{(j_2)}\mathbf{E}\left(\mathbf{E}\epsilon_{i_1}^{(j_1)}|\mathcal{H}_{i_1, (s,s)}^{(j_1), (s, s)}\mathbf{E}\epsilon_{i_2}^{(j_2)}|\mathcal{H}_{i_2, (s,s)}^{(j_2), (s, s)}\right)
\vert\\
\leq \Vert \sum_{i_1 = 1}^{n_1}\sum_{j_1 = 1}^{n_1}S_{i_1,v_1}^{(j_1)}\Vert_2
\times \Vert\sum_{i_2 = 1}^{n_1}\sum_{j_2 = 1}^{n_1}S_{i_2,v_2}^{(j_2)} - \sum_{i_2 = 1}^n\sum_{j_2 = 1}^mc_{i_2, v_2}^{(j_2)}\mathbf{E}\epsilon_{i_2}^{(j_2)}|\mathcal{H}_{i_2, (s,s)}^{(j_2), (s, s)}\Vert_2\\
+ \Vert
\sum_{i_1 = 1}^{n_1}\sum_{j_1 = 1}^{n_1}S_{i_1,v_1}^{(j_1)} - \sum_{i_1 = 1}^n\sum_{j_1 = 1}^m c_{i_1, v_1}^{(j_1)}\mathbf{E}\epsilon_{i_1}^{(j_1)}|\mathcal{H}_{i_1, (s,s)}^{(j_1), (s, s)}
\Vert_2\times \Vert\sum_{i_2 = 1}^{n_1}\sum_{j_2 = 1}^{n_1}S_{i_2,v_2}^{(j_2)}\Vert_2\\
+ \Vert\sum_{i_1 = 1}^{n_1}\sum_{j_1 = 1}^{n_1}S_{i_1,v_1}^{(j_1)} - \sum_{i_1 = 1}^n\sum_{j_1 = 1}^m c_{i_1, v_1}^{(j_1)}\mathbf{E}\epsilon_{i_1}^{(j_1)}|\mathcal{H}_{i_1, (s,s)}^{(j_1), (s, s)}\Vert_2\\
\times\Vert\sum_{i_2 = 1}^{n_1}\sum_{j_2 = 1}^{n_1}S_{i_2,v_2}^{(j_2)} - \sum_{i_2 = 1}^n\sum_{j_2 = 1}^mc_{i_2, v_2}^{(j_2)}\mathbf{E}\epsilon_{i_2}^{(j_2)}|\mathcal{H}_{i_2, (s,s)}^{(j_2), (s, s)}\Vert_2\leq C\sqrt{\frac{s}{l}}
\end{align*}
From lemma \ref{lemma.consistent_linear_combination},
\begin{align*}
\vert\sum_{i_1 = 1}^n\sum_{j_1 = 1}^m\sum_{i_2 = 1}^n\sum_{j_2 = 1}^m c_{i_1, v_1}^{(j_1)}c_{i_2, v_2}^{(j_2)}\left[\mathbf{E}\left(\mathbf{E}\epsilon_{i_1}^{(j_1)}|\mathcal{H}_{i_1, (s,s)}^{(j_1), (s, s)}\mathbf{E}\epsilon_{i_2}^{(j_2)}|\mathcal{H}_{i_2, (s,s)}^{(j_2), (s, s)}\right) - \mathbf{E}\epsilon_{i_1}^{(j_1)}\epsilon_{i_2}^{(j_2)}\right]\vert\\
\leq \Vert \sum_{i = 1}^n\sum_{j = 1}^m c_{i, v_1}^{(j)}\mathbf{E}\epsilon_{i}^{(j)}|\mathcal{H}_{i, (s,s)}^{(j), (s, s)}\Vert_2
\times \Vert \sum_{i = 1}^n\sum_{j = 1}^m c_{i, v_2}^{(j)}(\epsilon_i^{(j)} - \mathbf{E}\epsilon_{i}^{(j)}|\mathcal{H}_{i, (s,s)}^{(j), (s, s)})\Vert_2\\
+ \Vert\sum_{i = 1}^n\sum_{j = 1}^m c_{i, v_1}^{(j)}(\epsilon_i^{(j)} - \mathbf{E}\epsilon_{i}^{(j)}|\mathcal{H}_{i, (s,s)}^{(j), (s, s)})\Vert_2\times \Vert\sum_{i = 1}^n\sum_{j = 1}^m c_{i, v_2}^{(j)}\mathbf{E}\epsilon_{i}^{(j)}|\mathcal{H}_{i, (s,s)}^{(j), (s, s)}\Vert_2\\
+ \Vert\sum_{i = 1}^n\sum_{j = 1}^m c_{i, v_1}^{(j)}(\epsilon_i^{(j)} - \mathbf{E}\epsilon_{i}^{(j)}|\mathcal{H}_{i, (s,s)}^{(j), (s, s)})\Vert_2\times\Vert\sum_{i = 1}^n\sum_{j = 1}^m c_{i, v_2}^{(j)}(\epsilon_i^{(j)} - \mathbf{E}\epsilon_{i}^{(j)}|\mathcal{H}_{i, (s,s)}^{(j), (s, s)})\Vert_2\leq C(1 + s)^{3-\alpha}
\end{align*}
From eq.\eqref{eq.delta_h},
\begin{equation}
\begin{aligned}
\vert\sum_{i_1 = 1}^{n_1}\sum_{j_1 = 1}^{n_1}\sum_{i_2 = 1}^{n_1}\sum_{j_2 = 1}^{n_2} \mathbf{E}S_{i_1,v_1}^{(j_1)}S_{i_2,v_2}^{(j_2)} - \sum_{i_1 = 1}^n\sum_{j_1 = 1}^m\sum_{i_2 = 1}^n\sum_{j_2 = 1}^m c_{i_1, v_1}^{(j_1)}c_{i_2, v_2}^{(j_2)}\mathbf{E}\epsilon_{i_1}^{(j_1)}\epsilon_{i_2}^{(j_2)}\vert\\
\leq C\left((1 + s)^{3-\alpha} + \sqrt{\frac{s}{l}}\right)\\
\Rightarrow \sup_{x\in\mathbf{R}}\vert\mathbf{E}h_{\tau,\psi,x}\left(\sum_{k = 1}^{n_1^2}P_{k,1}^*,\cdots, \sum_{k = 1}^{n_1^2}P_{k,V}^*\right)  - \mathbf{E}h_{\tau,\psi,x}\left(\xi_1,\cdots, \xi_V\right)\vert\\
\leq C(\psi^2 + \psi\tau)\times\left((1 + s)^{3-\alpha} + \sqrt{\frac{s}{l}}\right)
\end{aligned}
\label{eq.difference_h4}
\end{equation}
Choose $\tau = \psi = \log^9(\mathcal{K})$, then
$t \leq \frac{C}{\log^8(\mathcal{K})}$ and $Ct(1 + \sqrt{\log(V)} + \sqrt{\vert\log(t)\vert})\to 0$ as $\min(n,m)\to\infty$. Suppose $\frac{\alpha_V}{M} = \frac{\alpha  - 3}{1 + 5(\alpha - 3)} - \delta$ with $\delta > 0$.
Choose $s = \lfloor\mathcal{K}^{\alpha_s}\rfloor$ and $l = \lfloor \mathcal{K}^{\alpha_l}\rfloor$ with
$$
\alpha_s = \frac{1}{1 + 5(\alpha - 3)} \text{and } \alpha_l = \frac{1 + 2(\alpha - 3)}{1 + 5(\alpha - 3)} + \delta
$$
here $\lfloor x\rfloor$ denotes the largest integer that is smaller than or equal to $x$. For
\begin{align*}
\frac{\alpha_V}{M} - (\alpha - 3)\alpha_s = -\delta<0,\ \frac{\alpha_V}{M} + \frac{\alpha_s}{2} - \frac{\alpha_l}{2} = -\frac{3\delta}{2} < 0\\
\text{and } \frac{3\alpha_V}{M} + \alpha_l - 1 = -2\delta < 0
\end{align*}
from \eqref{eq.difference_h0}, \eqref{eq.difference_h1}, \eqref{eq.difference_h2}, \eqref{eq.difference_h3}, and \eqref{eq.difference_h4}, we prove \eqref{eq.Gaussian_result}.
\end{proof}

\begin{proof}[proof of theorem \ref{theorem.covariances}]
For 
\begin{align*}
    Cov\left(\sum_{i_1 = 1}^n\sum_{j_1 = 1}^m c_{i_1,v_1}^{(j_1)}\epsilon_{i_1}^{(j_1)}, \sum_{i_2 = 1}^n\sum_{j_2 = 1}^mc_{i_2,v_2}^{(j_2)}\epsilon_{i_2}^{(j_2)}\right) = \sum_{i_1 = 1}^n\sum_{j_1 = 1}^m\sum_{i_2 = 1}^n\sum_{j_2 = 1}^m c_{i_1,v_1}^{(j_1)}c_{i_2,v_2}^{(j_2)}\mathbf{E}\epsilon_{i_1}^{(j_1)}\epsilon_{i_2}^{(j_2)}
\end{align*}
notice that
\begin{align*}
\vert
\sum_{i_1 = 1}^n\sum_{j_1 = 1}^m\sum_{i_2 = 1}^n\sum_{j_2 = 1}^m c_{i_1,v_1}^{(j_1)}c_{i_2,v_2}^{(j_2)}\epsilon_{i_1}^{(j_1)}\epsilon_{i_2}^{(j_2)}
\times K\left(\frac{i_1 - i_2}{\mathcal{B}}\right)K\left(\frac{j_1 - j_2}{\mathcal{B}}\right)\\
-\sum_{i_1 = 1}^n\sum_{j_1 = 1}^m\sum_{i_2 = 1}^n\sum_{j_2 = 1}^m c_{i_1,v_1}^{(j_1)}c_{i_2,v_2}^{(j_2)}\mathbf{E}\epsilon_{i_1}^{(j_1)}\epsilon_{i_2}^{(j_2)}
\vert\\
\leq \vert \sum_{i_1 = 1}^n\sum_{j_1 = 1}^m\sum_{i_2 = 1}^n\sum_{j_2 = 1}^m c_{i_1,v_1}^{(j_1)}c_{i_2,v_2}^{(j_2)}
\times K\left(\frac{i_1 - i_2}{\mathcal{B}}\right)K\left(\frac{j_1 - j_2}{\mathcal{B}}\right)(\epsilon_{i_1}^{(j_1)}\epsilon_{i_2}^{(j_2)} - \mathbf{E}\epsilon_{i_1}^{(j_1)}\epsilon_{i_2}^{(j_2)})\vert\\
+ \sum_{i_1 = 1}^n\sum_{j_1 = 1}^m\sum_{i_2 = 1}^n\sum_{j_2 = 1}^m \vert c_{i_1,v_1}^{(j_1)}c_{i_2,v_2}^{(j_2)}\vert
\times K\left(\frac{i_1 - i_2}{\mathcal{B}}\right)\times \left(1 - K\left(\frac{j_1 - j_2}{\mathcal{B}}\right)\right) \vert\mathbf{E}\epsilon_{i_1}^{(j_1)}\epsilon_{i_2}^{(j_2)}\vert\\
+ \sum_{i_1 = 1}^n\sum_{j_1 = 1}^m\sum_{i_2 = 1}^n\sum_{j_2 = 1}^m \vert c_{i_1,v_1}^{(j_1)}c_{i_2,v_2}^{(j_2)}\vert
\times \left(1 - K\left(\frac{i_1 - i_2}{\mathcal{B}}\right)\right) \vert\mathbf{E}\epsilon_{i_1}^{(j_1)}\epsilon_{i_2}^{(j_2)}\vert
\end{align*}
Form remark \ref{remark.covariance},
\begin{align*}
\sum_{i_1 = 1}^n\sum_{j_1 = 1}^m\sum_{i_2 = 1}^n\sum_{j_2 = 1}^m \vert c_{i_1,v_1}^{(j_1)}c_{i_2,v_2}^{(j_2)}\vert
\times K\left(\frac{i_1 - i_2}{\mathcal{B}}\right)\times \left(1 - K\left(\frac{j_1 - j_2}{\mathcal{B}}\right)\right) \vert\mathbf{E}\epsilon_{i_1}^{(j_1)}\epsilon_{i_2}^{(j_2)}\vert\\
\leq C\sum_{s = 1 - n}^{n - 1}\sum_{t = 1 - m}^{m - 1}\left(1 - K\left(\frac{t}{\mathcal{B}}\right)\right)\times (1 + \max(\vert s\vert,\vert t\vert))^{3 - \alpha}\times \sum_{i_1 = 1 \vee (1 - s)}^{n\wedge (n - s)}\sum_{j_1 = 1 \vee (1 - t)}^{m\wedge (m - t)}
\vert c_{i_1,v_1}^{(j_1)}c_{i_1 + s,v_2}^{(j_1 + t)}\vert
\end{align*}
From Cauchy inequality
\begin{align*}
\sum_{i_1 = 1 \vee (1 - s)}^{n\wedge (n - s)}\sum_{j_1 = 1 \vee (1 - t)}^{m\wedge (m - t)}
\vert c_{i_1,v_1}^{(j_1)}c_{i_1 + s,v_2}^{(j_1 + t)}\vert
\leq \sqrt{\sum_{i = 1}^n\sum_{j = 1}^m c_{i,v_1}^{(j)2}}\times \sqrt{\sum_{i = 1}^n\sum_{j = 1}^m c_{i,v_2}^{(j)2}} = 1
\end{align*}
Also notice that
\begin{align*}
\sum_{s = 1 - n}^{n - 1}\sum_{t = 1 - m}^{m - 1}\left(1 - K\left(\frac{t}{\mathcal{B}}\right)\right)\times (1 + \max(\vert s\vert,\vert t\vert))^{3 - \alpha}\\
\leq 4\sum_{s = 0}^{n - 1}\sum_{t = 0}^{m - 1}\left(1 - K\left(\frac{t}{\mathcal{B}}\right)\right)\times (1 + \max(\vert s\vert,\vert t\vert))^{3 - \alpha}\\
\leq 4\sum_{t = 0}^{m - 1}\sum_{s = 0}^{t} \left(1 - K\left(\frac{t}{\mathcal{B}}\right)\right)\times (1 + t)^{3 - \alpha} +
4\sum_{t = 0}^{m - 1}\sum_{s = t+1}^\infty\left(1 - K\left(\frac{t}{\mathcal{B}}\right)\right)\times (1 + s)^{3 - \alpha}\\
\leq C\sum_{t = 0}^{m - 1}\left(1 - K\left(\frac{t}{\mathcal{B}}\right)\right)\times (1 + t)^{4 - \alpha}
\end{align*}
Define $k_0 = \max_{x\in[0,1]}\vert K^\prime(x)\vert$, then
\begin{align*}
\sum_{t = 0}^{m - 1}\left(1 - K\left(\frac{t}{\mathcal{B}}\right)\right)\times (1 + t)^{4 - \alpha}
\leq \frac{k_0}{\mathcal{B}}\sum_{t = 0}^{\lfloor\mathcal{B}\rfloor}(1 + t)^{5 - \alpha} + \sum_{t = \lfloor \mathcal{B}\rfloor + 1}^\infty (1 + t)^{4 - \alpha}\\
\leq \frac{C}{\mathcal{B}}\left(2 + \int_{[1,\mathcal{B}]}(1 + t)^{5 - \alpha}\mathrm{d}t\right) + C\int_{[B,\infty)}(1 + t)^{4-\alpha}\mathrm{d}t
\end{align*}
so we have
\begin{equation}
\begin{aligned}
\sum_{i_1 = 1}^n\sum_{j_1 = 1}^m\sum_{i_2 = 1}^n\sum_{j_2 = 1}^m \vert c_{i_1,v_1}^{(j_1)}c_{i_2,v_2}^{(j_2)}\vert
\times K\left(\frac{i_1 - i_2}{\mathcal{B}}\right)\times \left(1 - K\left(\frac{j_1 - j_2}{\mathcal{B}}\right)\right) \vert\mathbf{E}\epsilon_{i_1}^{(j_1)}\epsilon_{i_2}^{(j_2)}\vert\leq CS(\mathcal{B})
\end{aligned}
\end{equation}
Similarly
\begin{align*}
\sum_{i_1 = 1}^n\sum_{j_1 = 1}^m\sum_{i_2 = 1}^n\sum_{j_2 = 1}^m \vert c_{i_1,v_1}^{(j_1)}c_{i_2,v_2}^{(j_2)}\vert
\times \left(1 - K\left(\frac{i_1 - i_2}{\mathcal{B}}\right)\right) \vert\mathbf{E}\epsilon_{i_1}^{(j_1)}\epsilon_{i_2}^{(j_2)}\vert\\
\leq C\sum_{s = 1 - n}^{n - 1}\sum_{t = 1 - m}^{m - 1}\left(1 - K\left(\frac{s}{\mathcal{B}}\right)\right)\times (1 + \max(\vert s\vert,\vert t\vert))^{3-\alpha}\\
\times \sum_{i_1 = 1 \vee (1 - s)}^{n\wedge (n - s)}\sum_{j_1 = 1 \vee (1 - t)}^{m\wedge (m - t)}
\vert c_{i_1,v_1}^{(j_1)}c_{i_1 + s,v_2}^{(j_1 + t)}\vert
\end{align*}
so we have
\begin{equation}
\begin{aligned}
\sum_{i_1 = 1}^n\sum_{j_1 = 1}^m\sum_{i_2 = 1}^n\sum_{j_2 = 1}^m \vert c_{i_1,v_1}^{(j_1)}c_{i_2,v_2}^{(j_2)}\vert
\times \left(1 - K\left(\frac{i_1 - i_2}{\mathcal{B}}\right)\right) \vert\mathbf{E}\epsilon_{i_1}^{(j_1)}\epsilon_{i_2}^{(j_2)}\vert\leq CS(\mathcal{B})
\end{aligned}
\end{equation}
On the other hand,
\begin{align*}
\Vert
\sum_{i_1 = 1}^n\sum_{j_1 = 1}^m\sum_{i_2 = 1}^n\sum_{j_2 = 1}^m c_{i_1,v_1}^{(j_1)}c_{i_2,v_2}^{(j_2)}
\times K\left(\frac{i_1 - i_2}{\mathcal{B}}\right)K\left(\frac{j_1 - j_2}{\mathcal{B}}\right)(\epsilon_{i_1}^{(j_1)}\epsilon_{i_2}^{(j_2)} - \mathbf{E}\epsilon_{i_1}^{(j_1)}\epsilon_{i_2}^{(j_2)})
\Vert_{M/2}\\
\leq \sum_{s = 1 - n}^{n - 1}\sum_{t = 1 - m}^{m - 1}K\left(\frac{s}{\mathcal{B}}\right)K\left(\frac{t}{\mathcal{B}}\right)
\Vert\sum_{i_1 = 1\vee (1 - s)}^{n\wedge (n - s)}\sum_{j_1 = 1\vee (1 - t)}^{m\wedge (m - t)}
c_{i_1,v_1}^{(j_1)}c_{i_1 + s,v_2}^{(j_1 + t)}(\epsilon_{i_1}^{(j_1)}\epsilon_{i_1+s}^{(j_1+t)} - \mathbf{E}\epsilon_{i_1}^{(j_1)}\epsilon_{i_1+s}^{(j_1+t)})
\Vert_{M/2}
\end{align*}
For
\begin{align*}
\Vert\sum_{i_1 = 1\vee (1 - s)}^{n\wedge (n - s)}\sum_{j_1 = 1\vee (1 - t)}^{m\wedge (m - t)}
c_{i_1,v_1}^{(j_1)}c_{i_1 + s,v_2}^{(j_1 + t)}(\epsilon_{i_1}^{(j_1)}\epsilon_{i_1+s}^{(j_1+t)} - \mathbf{E}\epsilon_{i_1}^{(j_1)}\epsilon_{i_1+s}^{(j_1+t)})
\Vert_{M/2}\\
\leq
\Vert\sum_{i_1 = 1\vee (1 - s)}^{n\wedge (n - s)}\sum_{j_1 = 1\vee (1 - t)}^{m\wedge (m - t)}
c_{i_1,v_1}^{(j_1)}c_{i_1 + s,v_2}^{(j_1 + t)}(\mathbf{E}\epsilon_{i_1}^{(j_1)}\epsilon_{i_1+s}^{(j_1+t)}|\mathcal{H}_{i_1, (0, 0)}^{(j_1), (0,0)} - \mathbf{E}\epsilon_{i_1}^{(j_1)}\epsilon_{i_1+s}^{(j_1+t)})
\Vert_{M/2}\\
+\sum_{l = 1}^\infty\Vert
\sum_{i_1 = 1\vee (1 - s)}^{n\wedge (n - s)}\sum_{j_1 = 1\vee (1 - t)}^{m\wedge (m - t)}
c_{i_1,v_1}^{(j_1)}c_{i_1 + s,v_2}^{(j_1 + t)}(\mathbf{E}\epsilon_{i_1}^{(j_1)}\epsilon_{i_1+s}^{(j_1+t)}|\mathcal{H}_{i_1, (l, l - 1)}^{(j_1), (l - 1,l - 1)} - \mathbf{E}\epsilon_{i_1}^{(j_1)}\epsilon_{i_1+s}^{(j_1+t)}|\mathcal{H}_{i_1, (l - 1, l - 1)}^{(j_1), (l - 1,l - 1)})
\Vert_{M/2}\\
+ \sum_{l = 1}^\infty\Vert
\sum_{i_1 = 1\vee (1 - s)}^{n\wedge (n - s)}\sum_{j_1 = 1\vee (1 - t)}^{m\wedge (m - t)}
c_{i_1,v_1}^{(j_1)}c_{i_1 + s,v_2}^{(j_1 + t)}(\mathbf{E}\epsilon_{i_1}^{(j_1)}\epsilon_{i_1+s}^{(j_1+t)}|\mathcal{H}_{i_1, (l, l)}^{(j_1), (l - 1,l - 1)} - \mathbf{E}\epsilon_{i_1}^{(j_1)}\epsilon_{i_1+s}^{(j_1+t)}|\mathcal{H}_{i_1, (l, l - 1)}^{(j_1), (l - 1,l - 1)})
\Vert_{M/2}\\
+ \sum_{l = 1}^\infty\Vert
\sum_{i_1 = 1\vee (1 - s)}^{n\wedge (n - s)}\sum_{j_1 = 1\vee (1 - t)}^{m\wedge (m - t)}
c_{i_1,v_1}^{(j_1)}c_{i_1 + s,v_2}^{(j_1 + t)}(\mathbf{E}\epsilon_{i_1}^{(j_1)}\epsilon_{i_1+s}^{(j_1+t)}|\mathcal{H}_{i_1, (l, l)}^{(j_1), (l,l - 1)} - \mathbf{E}\epsilon_{i_1}^{(j_1)}\epsilon_{i_1+s}^{(j_1+t)}|\mathcal{H}_{i_1, (l, l)}^{(j_1), (l - 1,l - 1)})
\Vert_{M/2}\\
+\sum_{l = 1}^\infty\Vert
\sum_{i_1 = 1\vee (1 - s)}^{n\wedge (n - s)}\sum_{j_1 = 1\vee (1 - t)}^{m\wedge (m - t)}
c_{i_1,v_1}^{(j_1)}c_{i_1 + s,v_2}^{(j_1 + t)}(\mathbf{E}\epsilon_{i_1}^{(j_1)}\epsilon_{i_1+s}^{(j_1+t)}|\mathcal{H}_{i_1, (l, l)}^{(j_1), (l,l)} - \mathbf{E}\epsilon_{i_1}^{(j_1)}\epsilon_{i_1+s}^{(j_1+t)}|\mathcal{H}_{i_1, (l, l)}^{(j_1), (l,l - 1)})
\Vert_{M/2}
\end{align*}
From \cite{MR0133849}, for
\begin{align*}
\Vert\mathbf{E}\epsilon_{i_1}^{(j_1)}\epsilon_{i_1+s}^{(j_1+t)}|\mathcal{H}_{i_1, (0, 0)}^{(j_1), (0,0)} - \mathbf{E}\epsilon_{i_1}^{(j_1)}\epsilon_{i_1+s}^{(j_1+t)}\Vert_{M/2}\leq 2\Vert\epsilon_{i_1}^{(j_1)}\epsilon_{i_1+s}^{(j_1+t)}\Vert_{M/2}
\leq 2\Vert\epsilon_{i_1}^{(j_1)}\Vert_M\Vert\epsilon_{i_1+s}^{(j_1+t)}\Vert_M\leq C
\end{align*}
we have
\begin{equation}
\begin{aligned}
\Vert\sum_{i_1 = 1\vee (1 - s)}^{n\wedge (n - s)}\sum_{j_1 = 1\vee (1 - t)}^{m\wedge (m - t)}
c_{i_1,v_1}^{(j_1)}c_{i_1 + s,v_2}^{(j_1 + t)}(\mathbf{E}\epsilon_{i_1}^{(j_1)}\epsilon_{i_1+s}^{(j_1+t)}|\mathcal{H}_{i_1, (0, 0)}^{(j_1), (0,0)} - \mathbf{E}\epsilon_{i_1}^{(j_1)}\epsilon_{i_1+s}^{(j_1+t)})
\Vert_{M/2}\\
\leq C\sqrt{\sum_{i_1 = 1\vee (1 - s)}^{n\wedge (n - s)}\sum_{j_1 = 1\vee (1 - t)}^{m\wedge (m - t)}c_{i_1,v_1}^{(j_1)2}c_{i_1 + s,v_2}^{(j_1 + t)2}\Vert\mathbf{E}\epsilon_{i_1}^{(j_1)}\epsilon_{i_1+s}^{(j_1+t)}|\mathcal{H}_{i_1, (0, 0)}^{(j_1), (0,0)} - \mathbf{E}\epsilon_{i_1}^{(j_1)}\epsilon_{i_1+s}^{(j_1+t)}\Vert_{M/2}^2}\\
\leq C^\prime \max_{i = 1,\cdots,n,j = 1,\cdots, m, v = 1,\cdots, V}\vert   c_{i,v}^{(j)}\vert\times \sqrt{\sum_{i_1 = 1\vee (1 - s)}^{n\wedge (n - s)}\sum_{j_1 = 1\vee (1 - t)}^{m\wedge (m - t)}c_{i_1 + s,v_2}^{(j_1 + t)2}}
\leq \frac{C^{\prime\prime}}{\mathcal{K}}
\end{aligned}
\end{equation}
here $C^{\prime\prime}$ is a constant. On the other hand, notice that $\epsilon_i^{(j)}\epsilon_{i + s}^{(j + t)}$ is also a function of $\mathcal{Z}$, from lemma \ref{lemma.calculate_martingale_difference}
\begin{align*}
\Vert
\sum_{i_1 = 1\vee (1 - s)}^{n\wedge (n - s)}\sum_{j_1 = 1\vee (1 - t)}^{m\wedge (m - t)}
c_{i_1,v_1}^{(j_1)}c_{i_1 + s,v_2}^{(j_1 + t)}(\mathbf{E}\epsilon_{i_1}^{(j_1)}\epsilon_{i_1+s}^{(j_1+t)}|\mathcal{H}_{i_1, (l, l - 1)}^{(j_1), (l - 1,l - 1)} - \mathbf{E}\epsilon_{i_1}^{(j_1)}\epsilon_{i_1+s}^{(j_1+t)}|\mathcal{H}_{i_1, (l - 1, l - 1)}^{(j_1), (l - 1,l - 1)})
\Vert_{M/2}\\
\leq C\sqrt{\sum_{i_1 = 1\vee (1 - s)}^{n\wedge (n - s)}\Vert\sum_{j_1 = 1\vee (1 - t)}^{m\wedge (m - t)}
c_{i_1,v_1}^{(j_1)}c_{i_1 + s,v_2}^{(j_1 + t)}(\mathbf{E}\epsilon_{i_1}^{(j_1)}\epsilon_{i_1+s}^{(j_1+t)}|\mathcal{H}_{i_1, (l, l - 1)}^{(j_1), (l - 1,l - 1)} - \mathbf{E}\epsilon_{i_1}^{(j_1)}\epsilon_{i_1+s}^{(j_1+t)}|\mathcal{H}_{i_1, (l - 1, l - 1)}^{(j_1), (l - 1,l - 1)})
\Vert_{M/2}^2}
\end{align*}
Define $p_{i,l}^{(j)} = \mathbf{E}\epsilon_{i}^{(j)}\epsilon_{i+s}^{(j+t)}|\mathcal{H}_{i, (l, l - 1)}^{(j), (l - 1,l - 1)} - \mathbf{E}\epsilon_{i}^{(j)}\epsilon_{i+s}^{(j+t)}|\mathcal{H}_{i, (l - 1, l - 1)}^{(j), (l - 1,l - 1)}$, then
$
p_{i,l}^{(j)} = \mathbf{E}p_{i,l}^{(j)}|\mathcal{H}_{i, (l, l - 1)}^{(j), (l - 1, l - 1)}
$, and the set $B = \{(i,j): i = 1\vee (1 - s),\cdots, n\wedge (n - s), j = 1\vee (1 - t),\cdots, m\wedge (m - t)\}$, from lemma \ref{lemma.second_dimension_difference}
\begin{equation}
\begin{aligned}
\Vert\sum_{j_1 = 1\vee (1 - t)}^{m\wedge (m - t)}
c_{i_1,v_1}^{(j_1)}c_{i_1 + s,v_2}^{(j_1 + t)}(\mathbf{E}\epsilon_{i_1}^{(j_1)}\epsilon_{i_1+s}^{(j_1+t)}|\mathcal{H}_{i_1, (l, l - 1)}^{(j_1), (l - 1,l - 1)} - \mathbf{E}\epsilon_{i_1}^{(j_1)}\epsilon_{i_1+s}^{(j_1+t)}|\mathcal{H}_{i_1, (l - 1, l - 1)}^{(j_1), (l - 1,l - 1)})
\Vert_{M/2}\\
\leq C\sqrt{\sum_{j_1 = 1\vee (1 - t)}^{m\wedge (m - t)}c_{i_1,v_1}^{(j_1)2}c_{i_1 + s,v_2}^{(j_1 + t)2}}\times
\max_{(i,j)\in B}\Vert\mathbf{E}p_{i,l}^{(j)}|\mathcal{H}_{i, (l, l - 1)}^{(j), (0, 0)}\Vert_{M/2}\\
+ C\sqrt{\sum_{j_1 = 1\vee (1 - t)}^{m\wedge (m - t)}c_{i_1,v_1}^{(j_1)2}c_{i_1 + s,v_2}^{(j_1 + t)2}}\times
\sum_{k = 1}^{l - 1}\max_{(i,j)\in B}\Vert\mathbf{E}p_{i,l}^{(j)}|\mathcal{H}_{i, (l, l - 1)}^{(j), (k, k - 1)} -
\mathbf{E}p_{i,l}^{(j)}|\mathcal{H}_{i, (l, l - 1)}^{(j), (k - 1, k - 1)}
\Vert_{M/2}\\
+ C\sqrt{\sum_{j_1 = 1\vee (1 - t)}^{m\wedge (m - t)}c_{i_1,v_1}^{(j_1)2}c_{i_1 + s,v_2}^{(j_1 + t)2}}\times
\sum_{k = 1}^{l - 1}
\max_{(i,j)\in B}\Vert
\mathbf{E}p_{i,l}^{(j)}|\mathcal{H}_{i, (l, l - 1)}^{(j), (k, k)} -
\mathbf{E}p_{i,l}^{(j)}|\mathcal{H}_{i, (l, l - 1)}^{(j), (k , k - 1)}
\Vert_{M/2}
\end{aligned}
\end{equation}
From lemma \ref{lemma.conditional_prob}
\begin{align*}
\Vert\mathbf{E}p_{i,l}^{(j)}|\mathcal{H}_{i, (l, l - 1)}^{(j), (0, 0)}\Vert_{M/2}
= \Vert
\mathbf{E}\epsilon_{i}^{(j)}\epsilon_{i+s}^{(j+t)}|\mathcal{H}_{i, (l, l - 1)}^{(j), (0,0)} - \mathbf{E}\epsilon_{i}^{(j)}\epsilon_{i+s}^{(j+t)}|\mathcal{H}_{i, (l - 1, l - 1)}^{(j), (0,0)}
\Vert_{M/2}\\
= \Vert\mathbf{E}(\epsilon_i^{(j)}\epsilon_{i + s}^{(j + t)} - \epsilon_{i, (-l)}^{(j), (0)}\epsilon_{i + s, (-s-l)}^{(j + t), (-t)})|\mathcal{H}_{i, (l, l - 1)}^{(j), (0,0)}\Vert_{M/2}\\
\leq \Vert\epsilon_i^{(j)}\Vert_{M}\times \Vert\epsilon_{i + s}^{(j + t)} - \epsilon_{i + s, (-s-l)}^{(j + t), (-t)}\Vert_M +
\Vert\epsilon_{i + s, (-s-l)}^{(j + t), (-t)}\Vert_M\times \Vert\epsilon_i^{(j)} - \epsilon_{i, (-l)}^{(j), (0)}\Vert_M\\
\leq C(\delta_{-l}^{(0)} + \delta_{-l-s}^{(-t)})\leq \frac{C^\prime}{(1 + \vert l\vert)^\alpha} + \frac{C^\prime}{(1 + \vert l + s\vert + \vert t\vert)^\alpha}
\end{align*}
From lemma \ref{lemma.delta_moment}
\begin{align*}
\Vert\mathbf{E}p_{i,l}^{(j)}|\mathcal{H}_{i, (l, l - 1)}^{(j), (k, k - 1)} -
\mathbf{E}p_{i,l}^{(j)}|\mathcal{H}_{i, (l, l - 1)}^{(j), (k - 1, k - 1)}
\Vert_{M/2}
\leq \Vert\mathbf{E}\epsilon_{i}^{(j)}\epsilon_{i+s}^{(j+t)}|\mathcal{H}_{i, (l, l - 1)}^{(j), (k, k - 1)}\\
- \mathbf{E}\epsilon_{i}^{(j)}\epsilon_{i+s}^{(j+t)}|\mathcal{H}_{i, (l - 1, l - 1)}^{(j), (k, k - 1)}\Vert_{M/2} +  \Vert\mathbf{E}\epsilon_{i}^{(j)}\epsilon_{i+s}^{(j+t)}|\mathcal{H}_{i, (l, l - 1)}^{(j), (k - 1, k - 1)}
- \mathbf{E}\epsilon_{i}^{(j)}\epsilon_{i+s}^{(j+t)}|\mathcal{H}_{i, (l - 1, l - 1)}^{(j), (k - 1,k - 1)}\Vert_{M / 2}\\
\leq C\sum_{v = -k}^{k}\delta_{-l}^{(v)} + C\sum_{v = -k}^{k}\delta_{-l-s}^{v - t}
\leq C^\prime\sum_{v = -\infty}^{\infty}\frac{1}{(1 + \vert l\vert + \vert v\vert)^\alpha} + C^\prime\sum_{v = -\infty}^{\infty}
\frac{1}{(1 + \vert l + s\vert + \vert v\vert)^\alpha}\\
\leq C^{\prime\prime}(1 + l)^{1-\alpha} + C^{\prime\prime}(1 + \vert l + s\vert)^{1-\alpha}
\end{align*}
and
\begin{align*}
\Vert
\mathbf{E}p_{i,l}^{(j)}|\mathcal{H}_{i, (l, l - 1)}^{(j), (k, k)} -
\mathbf{E}p_{i,l}^{(j)}|\mathcal{H}_{i, (l, l - 1)}^{(j), (k , k - 1)}
\Vert_{M/2}
\leq
\Vert \mathbf{E}\epsilon_{i}^{(j)}\epsilon_{i+s}^{(j+t)}|\mathcal{H}_{i, (l, l - 1)}^{(j), (k,k)}\\
- \mathbf{E}\epsilon_{i}^{(j)}\epsilon_{i+s}^{(j+t)}|\mathcal{H}_{i, (l - 1, l - 1)}^{(j), (k,k)}\Vert_{M/2}
+ \Vert\mathbf{E}\epsilon_{i}^{(j)}\epsilon_{i+s}^{(j+t)}|\mathcal{H}_{i, (l, l - 1)}^{(j), (k,k - 1)} - \mathbf{E}\epsilon_{i}^{(j)}\epsilon_{i+s}^{(j+t)}|\mathcal{H}_{i, (l - 1, l - 1)}^{(j), (k, k - 1)}\Vert_{M/2}\\
\leq C\sum_{v = -k}^k\delta_{-l}^{(v)} + C\sum_{v = -k}^k\delta_{-l - s}^{(v - t)}
\leq C^{\prime\prime}(1 + l)^{1-\alpha} + C^{\prime\prime}(1 + \vert l + s\vert)^{1-\alpha}
\end{align*}
here $C^{\prime\prime}$ is a constant. Therefore,
\begin{equation}
\begin{aligned}
\Vert\sum_{j_1 = 1\vee (1 - t)}^{m\wedge (m - t)}
c_{i_1,v_1}^{(j_1)}c_{i_1 + s,v_2}^{(j_1 + t)}(\mathbf{E}\epsilon_{i_1}^{(j_1)}\epsilon_{i_1+s}^{(j_1+t)}|\mathcal{H}_{i_1, (l, l - 1)}^{(j_1), (l - 1,l - 1)} - \mathbf{E}\epsilon_{i_1}^{(j_1)}\epsilon_{i_1+s}^{(j_1+t)}|\mathcal{H}_{i_1, (l - 1, l - 1)}^{(j_1), (l - 1,l - 1)})
\Vert_{M/2}\\
\leq C\sqrt{\sum_{j_1 = 1\vee (1 - t)}^{m\wedge (m - t)}c_{i_1,v_1}^{(j_1)2}c_{i_1 + s,v_2}^{(j_1 + t)2}}\times ((1 + l)^{2-\alpha} + (1 + \vert l + s\vert)^{1-\alpha}\times (1 + l))\\
\text{and }
\Vert\sum_{i_1 = 1\vee (1 - s)}^{n\wedge (n - s)}\sum_{j_1 = 1\vee (1 - t)}^{m\wedge (m - t)}
c_{i_1,v_1}^{(j_1)}c_{i_1 + s,v_2}^{(j_1 + t)}(\mathbf{E}\epsilon_{i_1}^{(j_1)}\epsilon_{i_1+s}^{(j_1+t)}|\mathcal{H}_{i_1, (0, 0)}^{(j_1), (0,0)} - \mathbf{E}\epsilon_{i_1}^{(j_1)}\epsilon_{i_1+s}^{(j_1+t)})
\Vert_{M/2}\\
\leq C\sqrt{\sum_{i_1 = 1\vee (1 - s)}^{n\wedge (n - s)}\sum_{j_1 = 1\vee (1 - t)}^{m\wedge (m - t)}c_{i_1,v_1}^{(j_1)2}c_{i_1 + s,v_2}^{(j_1 + t)2}}\times((1 + l)^{2-\alpha} + (1 + \vert l + s\vert)^{1-\alpha}\times(1 + l))\\
\leq \frac{C^\prime}{\mathcal{K}}\times((1 + l)^{2-\alpha} + (1 + \vert l + s\vert)^{2-\alpha} + \vert s\vert\times (1 + \vert l + s\vert)^{1-\alpha})
\end{aligned}
\end{equation}
In particular,
\begin{equation}
\begin{aligned}
\sum_{l = 1}^\infty\Vert
\sum_{i_1 = 1\vee (1 - s)}^{n\wedge (n - s)}\sum_{j_1 = 1\vee (1 - t)}^{m\wedge (m - t)}
c_{i_1,v_1}^{(j_1)}c_{i_1 + s,v_2}^{(j_1 + t)}(\mathbf{E}\epsilon_{i_1}^{(j_1)}\epsilon_{i_1+s}^{(j_1+t)}|\mathcal{H}_{i_1, (l, l - 1)}^{(j_1), (l - 1,l - 1)} - \mathbf{E}\epsilon_{i_1}^{(j_1)}\epsilon_{i_1+s}^{(j_1+t)}|\mathcal{H}_{i_1, (l - 1, l - 1)}^{(j_1), (l - 1,l - 1)})
\Vert_{M/2}\\
\leq \frac{C}{\mathcal{K}}\sum_{l = 1}^\infty((1 + l)^{2-\alpha} + (1 + \vert l + s\vert)^{2-\alpha} + \vert s\vert\times (1 + \vert l + s\vert)^{1-\alpha})\\
\leq \frac{C}{\mathcal{K}}\sum_{l = 1}^\infty(1 + l)^{2-\alpha} + \frac{C}{\mathcal{K}}\sum_{l = -\infty}^\infty(1 + \vert l\vert)^{2-\alpha}
+ \frac{C\vert s\vert}{\mathcal{K}}\sum_{l = -\infty}^\infty(1 + \vert l\vert)^{1-\alpha}
\leq \frac{C^\prime\times(1 + \vert s\vert)}{\mathcal{K}}
\end{aligned}
\end{equation}
Similarly,
\begin{align*}
\Vert
\sum_{i_1 = 1\vee (1 - s)}^{n\wedge (n - s)}\sum_{j_1 = 1\vee (1 - t)}^{m\wedge (m - t)}
c_{i_1,v_1}^{(j_1)}c_{i_1 + s,v_2}^{(j_1 + t)}(\mathbf{E}\epsilon_{i_1}^{(j_1)}\epsilon_{i_1+s}^{(j_1+t)}|\mathcal{H}_{i_1, (l, l)}^{(j_1), (l - 1,l - 1)} - \mathbf{E}\epsilon_{i_1}^{(j_1)}\epsilon_{i_1+s}^{(j_1+t)}|\mathcal{H}_{i_1, (l, l - 1)}^{(j_1), (l - 1,l - 1)})
\Vert_{M/2}\\
\leq C\sqrt{\sum_{i_1 = 1\vee (1 - s)}^{n\wedge (n - s)}\Vert\sum_{j_1 = 1\vee (1 - t)}^{m\wedge (m - t)}
c_{i_1,v_1}^{(j_1)}c_{i_1 + s,v_2}^{(j_1 + t)}(\mathbf{E}\epsilon_{i_1}^{(j_1)}\epsilon_{i_1+s}^{(j_1+t)}|\mathcal{H}_{i_1, (l, l)}^{(j_1), (l - 1,l - 1)} - \mathbf{E}\epsilon_{i_1}^{(j_1)}\epsilon_{i_1+s}^{(j_1+t)}|\mathcal{H}_{i_1, (l, l - 1)}^{(j_1), (l - 1,l - 1)})\Vert_{M/2}^2}
\end{align*}
Define $q_{i,l}^{(j)} = \mathbf{E}\epsilon_{i}^{(j)}\epsilon_{i+s}^{(j+t)}|\mathcal{H}_{i, (l, l)}^{(j), (l - 1,l - 1)} - \mathbf{E}\epsilon_{i}^{(j)}\epsilon_{i+s}^{(j+t)}|\mathcal{H}_{i, (l, l - 1)}^{(j), (l - 1,l - 1)}$, then
$q_{i,l}^{(j)} = \mathbf{E}q_{i,l}^{(j)}|\mathcal{H}_{i, (l,l)}^{(j), (l - 1, l - 1)}$, so lemma \ref{lemma.second_dimension_difference} implies
\begin{align*}
\Vert\sum_{j_1 = 1\vee (1 - t)}^{m\wedge (m - t)}
c_{i_1,v_1}^{(j_1)}c_{i_1 + s,v_2}^{(j_1 + t)}(\mathbf{E}\epsilon_{i_1}^{(j_1)}\epsilon_{i_1+s}^{(j_1+t)}|\mathcal{H}_{i_1, (l, l)}^{(j_1), (l - 1,l - 1)} - \mathbf{E}\epsilon_{i_1}^{(j_1)}\epsilon_{i_1+s}^{(j_1+t)}|\mathcal{H}_{i_1, (l, l - 1)}^{(j_1), (l - 1,l - 1)})\Vert_{M/2}\\
\leq C\sqrt{\sum_{j_1 = 1\vee (1 - t)}^{m\wedge (m - t)}
c_{i_1,v_1}^{(j_1)}c_{i_1 + s,v_2}^{(j_1 + t)2}}\times \max_{(i,j)\in B}\Vert\mathbf{E}q_{i,l}^{(j)}|\mathcal{H}_{i, (l,l)}^{(j), (0,0)} \Vert_{M/2}\\
+C\sqrt{\sum_{j_1 = 1\vee (1 - t)}^{m\wedge (m - t)}
c_{i_1,v_1}^{(j_1)}c_{i_1 + s,v_2}^{(j_1 + t)2}}\times \sum_{k = 1}^{l - 1}\max_{(i,j)\in B}\Vert
\mathbf{E}q_{i,l}^{(j)}|\mathcal{H}_{i, (l,l)}^{(j), (k,k - 1)} - \mathbf{E}q_{i,l}^{(j)}|\mathcal{H}_{i, (l,l)}^{(j), (k - 1,k - 1)}
\Vert_{M/2}\\
+C\sqrt{\sum_{j_1 = 1\vee (1 - t)}^{m\wedge (m - t)}
c_{i_1,v_1}^{(j_1)}c_{i_1 + s,v_2}^{(j_1 + t)2}}\times \sum_{k = 1}^{l - 1}\max_{(i,j)\in B}\Vert
\mathbf{E}q_{i,l}^{(j)}|\mathcal{H}_{i, (l,l)}^{(j), (k,k)} - \mathbf{E}q_{i,l}^{(j)}|\mathcal{H}_{i, (l,l)}^{(j), (k,k - 1)}
\Vert_{M/2}
\end{align*}
From lemma \ref{lemma.conditional_prob},
\begin{align*}
\Vert\mathbf{E}q_{i,l}^{(j)}|\mathcal{H}_{i, (l,l)}^{(j), (0,0)} \Vert_{M/2} = \Vert \mathbf{E}\epsilon_i^{(j)}\epsilon_{i + s}^{(j + t)}|\mathcal{H}_{i, (l,l)}^{(j), (0, 0)} - \mathbf{E}\epsilon_i^{(j)}\epsilon_{i + s}^{(j + t)}|\mathcal{H}_{i, (l,l - 1)}^{(j), (0, 0)}\Vert_{M / 2}\\
\leq C\delta_{l}^{(0)} + C\delta_{l - s}^{(-t)}\leq \frac{C^\prime}{(1 + l)^\alpha} + \frac{C^\prime}{(1 + \vert l - s\vert + \vert t\vert)^\alpha}
\end{align*}
From lemma \ref{lemma.delta_moment},
\begin{align*}
\Vert
\mathbf{E}q_{i,l}^{(j)}|\mathcal{H}_{i, (l,l)}^{(j), (k,k - 1)} - \mathbf{E}q_{i,l}^{(j)}|\mathcal{H}_{i, (l,l)}^{(j), (k - 1,k - 1)}
\Vert_{M / 2} \leq \Vert\mathbf{E}\epsilon_i^{(j)}\epsilon_{i + s}^{(j + t)}|\mathcal{H}_{i, (l, l)}^{(j), (k, k - 1)}\\
- \mathbf{E}\epsilon_i^{(j)}\epsilon_{i + s}^{(j + t)}|\mathcal{H}_{i, (l, l - 1)}^{(j),(k, k - 1)} \Vert_{M / 2}
+ \Vert
\mathbf{E}\epsilon_i^{(j)}\epsilon_{i + s}^{(j + t)}|\mathcal{H}_{i, (l, l)}^{(j), (k - 1, k - 1)}
- \mathbf{E}\epsilon_i^{(j)}\epsilon_{i + s}^{(j + t)}|\mathcal{H}_{i, (l, l - 1)}^{(j), (k - 1, k - 1)}
\Vert_{M / 2}\\
\leq C\sum_{v = -k}^k\delta_{l}^{(v)}  + C\sum_{v = -k}^k \delta_{l - s}^{(v - t)}\leq C^\prime(1 + l)^{1 - \alpha} + C^\prime(1 + \vert l - s\vert)^{1 - \alpha}
\end{align*}
and
\begin{align*}
\Vert
\mathbf{E}q_{i,l}^{(j)}|\mathcal{H}_{i, (l,l)}^{(j), (k,k)} - \mathbf{E}q_{i,l}^{(j)}|\mathcal{H}_{i, (l,l)}^{(j), (k,k - 1)}
\Vert_{M/2}
\leq \Vert \mathbf{E}\epsilon_i^{(j)}\epsilon_{i + s}^{(j + t)}|\mathcal{H}_{i, (l, l)}^{(j), (k, k)}\\
- \mathbf{E}\epsilon_i^{(j)}\epsilon_{i + s}^{(j + t)}|\mathcal{H}_{i, (l, l - 1)}^{(j), (k, k)}\Vert_{M / 2} +
\Vert \mathbf{E}\epsilon_i^{(j)}\epsilon_{i + s}^{(j + t)}|\mathcal{H}_{i, (l, l)}^{(j), (k, k - 1)}
- \mathbf{E}\epsilon_i^{(j)}\epsilon_{i + s}^{(j + t)}|\mathcal{H}_{i, (l, l - 1)}^{(j), (k, k - 1)}\Vert_{M / 2}\\
\leq C\sum_{v = -k}^k\delta_{l}^{(v)} + C\sum_{v = -k}^k \delta_{l - s}^{(v - t)}
\leq C^\prime(1 + l)^{1 - \alpha}  + C^\prime (1 + \vert l - s\vert)^{1-\alpha}
\end{align*}
Therefore
\begin{align*}
\Vert\sum_{j_1 = 1\vee (1 - t)}^{m\wedge (m - t)}
c_{i_1,v_1}^{(j_1)}c_{i_1 + s,v_2}^{(j_1 + t)}(\mathbf{E}\epsilon_{i_1}^{(j_1)}\epsilon_{i_1+s}^{(j_1+t)}|\mathcal{H}_{i_1, (l, l)}^{(j_1), (l - 1,l - 1)} - \mathbf{E}\epsilon_{i_1}^{(j_1)}\epsilon_{i_1+s}^{(j_1+t)}|\mathcal{H}_{i_1, (l, l - 1)}^{(j_1), (l - 1,l - 1)})\Vert_{M/2}\\
\leq C\sqrt{\sum_{j_1 = 1\vee (1 - t)}^{m\wedge (m - t)}
c_{i_1,v_1}^{(j_1)}c_{i_1 + s,v_2}^{(j_1 + t)2}}\times ((1 + l)^{2-\alpha} + (1 + \vert l - s\vert)^{2 - \alpha} + \vert s\vert(1 + \vert l - s\vert)^{1 - \alpha})\\
\text{and } \sum_{l = 1}^\infty\Vert
\sum_{i_1 = 1\vee (1 - s)}^{n\wedge (n - s)}\sum_{j_1 = 1\vee (1 - t)}^{m\wedge (m - t)}
c_{i_1,v_1}^{(j_1)}c_{i_1 + s,v_2}^{(j_1 + t)}(\mathbf{E}\epsilon_{i_1}^{(j_1)}\epsilon_{i_1+s}^{(j_1+t)}|\mathcal{H}_{i_1, (l, l)}^{(j_1), (l - 1,l - 1)} - \mathbf{E}\epsilon_{i_1}^{(j_1)}\epsilon_{i_1+s}^{(j_1+t)}|\mathcal{H}_{i_1, (l, l - 1)}^{(j_1), (l - 1,l - 1)})
\Vert_{M/2}\\
\leq \frac{C}{\mathcal{K}}\times \sum_{l = -\infty}^\infty((1 + l)^{2 - \alpha} + (1 + \vert l - s\vert)^{2-\alpha} + \vert s\vert(1 + \vert l - s\vert)^{1-\alpha})\leq \frac{C^\prime\times (1 + \vert s\vert)}{\mathcal{K}}
\end{align*}
Also,
\begin{align*}
\Vert
\sum_{i_1 = 1\vee (1 - s)}^{n\wedge (n - s)}\sum_{j_1 = 1\vee (1 - t)}^{m\wedge (m - t)}
c_{i_1,v_1}^{(j_1)}c_{i_1 + s,v_2}^{(j_1 + t)}(\mathbf{E}\epsilon_{i_1}^{(j_1)}\epsilon_{i_1+s}^{(j_1+t)}|\mathcal{H}_{i_1, (l, l)}^{(j_1), (l,l - 1)} - \mathbf{E}\epsilon_{i_1}^{(j_1)}\epsilon_{i_1+s}^{(j_1+t)}|\mathcal{H}_{i_1, (l, l)}^{(j_1), (l - 1,l - 1)})
\Vert_{M/2}\\
\leq C\sqrt{\sum_{j_1 = 1\vee (1 - t)}^{m\wedge (m - t)}\Vert\sum_{i_1 = 1\vee (1 - s)}^{n\wedge (n - s)}c_{i_1,v_1}^{(j_1)}c_{i_1 + s,v_2}^{(j_1 + t)}(\mathbf{E}\epsilon_{i_1}^{(j_1)}\epsilon_{i_1+s}^{(j_1+t)}|\mathcal{H}_{i_1, (l, l)}^{(j_1), (l,l - 1)} - \mathbf{E}\epsilon_{i_1}^{(j_1)}\epsilon_{i_1+s}^{(j_1+t)}|\mathcal{H}_{i_1, (l, l)}^{(j_1), (l - 1,l - 1)})
\Vert_{M / 2}^2}
\end{align*}
Define $r_{i, l}^{(j)} = \mathbf{E}\epsilon_{i}^{(j)}\epsilon_{i+s}^{(j+t)}|\mathcal{H}_{i, (l, l)}^{(j), (l,l - 1)} - \mathbf{E}\epsilon_{i}^{(j)}\epsilon_{i+s}^{(j+t)}|\mathcal{H}_{i, (l, l)}^{(j), (l - 1,l - 1)}$, then from lemma \ref{lemma.second_dimension_difference}
\begin{align*}
\Vert\sum_{i_1 = 1\vee (1 - s)}^{n\wedge (n - s)}c_{i_1,v_1}^{(j_1)}c_{i_1 + s,v_2}^{(j_1 + t)}(\mathbf{E}\epsilon_{i_1}^{(j_1)}\epsilon_{i_1+s}^{(j_1+t)}|\mathcal{H}_{i_1, (l, l)}^{(j_1), (l,l - 1)} - \mathbf{E}\epsilon_{i_1}^{(j_1)}\epsilon_{i_1+s}^{(j_1+t)}|\mathcal{H}_{i_1, (l, l)}^{(j_1), (l - 1,l - 1)})\Vert_{M/2}\\
\leq C\sqrt{\sum_{i_1 = 1\vee (1 - s)}^{n\wedge (n - s)}c_{i_1,v_1}^{(j_1)2}c_{i_1 + s,v_2}^{(j_1 + t)2}}\times
\max_{(i,j)\in B}\Vert\mathbf{E}r_{i, l}^{(j)}|\mathcal{H}_{i, (0, 0)}^{(j), (l, l - 1)}\Vert_{M / 2}\\
+ C\sqrt{\sum_{i_1 = 1\vee (1 - s)}^{n\wedge (n - s)}c_{i_1,v_1}^{(j_1)2}c_{i_1 + s,v_2}^{(j_1 + t)2}}\times\sum_{k = 1}^{l}\max_{(i,j)\in B }\Vert
\mathbf{E}r_{i, l}^{(j)}|\mathcal{H}_{i, (k, k - 1)}^{(j), (l, l - 1)} - \mathbf{E}r_{i, l}^{(j)}|\mathcal{H}_{i, (k  - 1, k - 1)}^{(j), (l, l - 1)}
\Vert_{M / 2}\\
+ C\sqrt{\sum_{i_1 = 1\vee (1 - s)}^{n\wedge (n - s)}c_{i_1,v_1}^{(j_1)2}c_{i_1 + s,v_2}^{(j_1 + t)2}}\times\sum_{k = 1}^{l}\max_{(i,j)\in B }\Vert
\mathbf{E}r_{i, l}^{(j)}|\mathcal{H}_{i, (k, k)}^{(j), (l, l - 1)} - \mathbf{E}r_{i, l}^{(j)}|\mathcal{H}_{i, (k, k - 1)}^{(j), (l, l - 1)}
\Vert_{M / 2}
\end{align*}
From lemma \ref{lemma.conditional_prob} and lemma \ref{lemma.delta_moment}
\begin{align*}
\Vert\mathbf{E}r_{i, l}^{(j)}|\mathcal{H}_{i, (0, 0)}^{(j), (l, l - 1)}\Vert_{M / 2}
= \Vert
\mathbf{E}\epsilon_{i}^{(j)}\epsilon_{i+s}^{(j+t)}|\mathcal{H}_{i, (0, 0)}^{(j), (l,l - 1)} - \mathbf{E}\epsilon_{i}^{(j)}\epsilon_{i+s}^{(j+t)}|\mathcal{H}_{i, (0, 0)}^{(j), (l - 1,l - 1)}
\Vert_{M / 2}\\
\leq C\delta_{0}^{(-l)} + C\delta_{-s}^{(-l-t)}\leq C^\prime(1 + l)^{-\alpha} + C^\prime(1 + \vert s\vert + \vert l + t\vert)^{-\alpha}\\
\text{and }\Vert
\mathbf{E}r_{i, l}^{(j)}|\mathcal{H}_{i, (k, k - 1)}^{(j), (l, l - 1)} - \mathbf{E}r_{i, l}^{(j)}|\mathcal{H}_{i, (k  - 1, k - 1)}^{(j), (l, l - 1)}
\Vert_{M / 2}\leq \Vert
\mathbf{E}\epsilon_i^{(j)}\epsilon_{i + s}^{(j + t)}|\mathcal{H}_{i, (k, k - 1)}^{(j), (l, l - 1)}\\
- \mathbf{E}\epsilon_i^{(j)}\epsilon_{i + s}^{(j + t)}|\mathcal{H}_{i, (k, k - 1)}^{(j), (l - 1, l - 1)}
\Vert_{M / 2} + \Vert \mathbf{E}\epsilon_i^{(j)}\epsilon_{i + s}^{(j + t)}|\mathcal{H}_{i, (k - 1, k - 1)}^{(j), (l, l - 1)}
- \mathbf{E}\epsilon_i^{(j)}\epsilon_{i + s}^{(j + t)}|\mathcal{H}_{i, (k - 1, k - 1)}^{(j), (l - 1, l - 1)}
\Vert_{M / 2}\\
\leq C\sum_{v = -k}^k\delta_{v}^{(-l)} + C\sum_{v = -k}^k \delta_{v - s}^{(-l-t)}\leq C^\prime(1 + l)^{1-\alpha} + C^\prime(1 + \vert l + t\vert)^{1-\alpha}\\
\text{and } \Vert
\mathbf{E}r_{i, l}^{(j)}|\mathcal{H}_{i, (k, k)}^{(j), (l, l - 1)} - \mathbf{E}r_{i, l}^{(j)}|\mathcal{H}_{i, (k, k - 1)}^{(j), (l, l - 1)}
\Vert_{M / 2}\leq \Vert \mathbf{E}\epsilon_i^{(j)}\epsilon_{i + s}^{(j + t)}|\mathcal{H}_{i, (k, k)}^{(j), (l, l - 1)}\\
- \mathbf{E}\epsilon_i^{(j)}\epsilon_{i + s}^{(j + t)}|\mathcal{H}_{i, (k, k)}^{(j), (l - 1, l - 1)}\Vert_{M/2} +
\Vert\mathbf{E}\epsilon_i^{(j)}\epsilon_{i + s}^{(j + t)}|\mathcal{H}_{i, (k, k - 1)}^{(j), (l, l - 1)}
- \mathbf{E}\epsilon_i^{(j)}\epsilon_{i + s}^{(j + t)}|\mathcal{H}_{i, (k, k - 1)}^{(j), (l - 1, l - 1)}
\Vert_{M / 2}\\
\leq C\sum_{v = -k}^k\delta_{v}^{(-l)} + C\sum_{v = -k}^k \delta_{v - s}^{(-l-t)}\leq C^\prime(1 + l)^{1-\alpha} + C^\prime(1 + \vert l + t\vert)^{1-\alpha}
\end{align*}
Therefore
\begin{align*}
\Vert
\sum_{i_1 = 1\vee (1 - s)}^{n\wedge (n - s)}\sum_{j_1 = 1\vee (1 - t)}^{m\wedge (m - t)}
c_{i_1,v_1}^{(j_1)}c_{i_1 + s,v_2}^{(j_1 + t)}(\mathbf{E}\epsilon_{i_1}^{(j_1)}\epsilon_{i_1+s}^{(j_1+t)}|\mathcal{H}_{i_1, (l, l)}^{(j_1), (l,l - 1)} - \mathbf{E}\epsilon_{i_1}^{(j_1)}\epsilon_{i_1+s}^{(j_1+t)}|\mathcal{H}_{i_1, (l, l)}^{(j_1), (l - 1,l - 1)})
\Vert_{M/2}\\
\leq C\sqrt{\sum_{j_1 = 1\vee (1 - t)}^{m\wedge (m - t)}\sum_{i_1 = 1\vee (1 - s)}^{n\wedge (n - s)}c_{i_1,v_1}^{(j_1)2}c_{i_1 + s,v_2}^{(j_1 + t)2}}
\times ((1 + l)^{2-\alpha} + (1 + \vert l + t\vert)^{2-\alpha} + \vert t\vert(1 + \vert l + t\vert)^{1-\alpha})\\
\text{and } \sum_{l = 1}^\infty\Vert
\sum_{i_1 = 1\vee (1 - s)}^{n\wedge (n - s)}\sum_{j_1 = 1\vee (1 - t)}^{m\wedge (m - t)}
c_{i_1,v_1}^{(j_1)}c_{i_1 + s,v_2}^{(j_1 + t)}(\mathbf{E}\epsilon_{i_1}^{(j_1)}\epsilon_{i_1+s}^{(j_1+t)}|\mathcal{H}_{i_1, (l, l)}^{(j_1), (l,l - 1)} - \mathbf{E}\epsilon_{i_1}^{(j_1)}\epsilon_{i_1+s}^{(j_1+t)}|\mathcal{H}_{i_1, (l, l)}^{(j_1), (l - 1,l - 1)})
\Vert_{M/2}\\
\leq \frac{C}{\mathcal{K}}\times(1 + \vert t\vert)
\end{align*}
Finally, notice that
\begin{align*}
\Vert
\sum_{i_1 = 1\vee (1 - s)}^{n\wedge (n - s)}\sum_{j_1 = 1\vee (1 - t)}^{m\wedge (m - t)}
c_{i_1,v_1}^{(j_1)}c_{i_1 + s,v_2}^{(j_1 + t)}(\mathbf{E}\epsilon_{i_1}^{(j_1)}\epsilon_{i_1+s}^{(j_1+t)}|\mathcal{H}_{i_1, (l, l)}^{(j_1), (l,l)} - \mathbf{E}\epsilon_{i_1}^{(j_1)}\epsilon_{i_1+s}^{(j_1+t)}|\mathcal{H}_{i_1, (l, l)}^{(j_1), (l,l - 1)})
\Vert_{M/2}\\
\leq C\sqrt{\sum_{j_1 = 1\vee (1 - t)}^{m\wedge (m - t)}\Vert\sum_{i_1 = 1\vee (1 - s)}^{n\wedge (n - s)}c_{i_1,v_1}^{(j_1)}c_{i_1 + s,v_2}^{(j_1 + t)}(\mathbf{E}\epsilon_{i_1}^{(j_1)}\epsilon_{i_1+s}^{(j_1+t)}|\mathcal{H}_{i_1, (l, l)}^{(j_1), (l,l)} - \mathbf{E}\epsilon_{i_1}^{(j_1)}\epsilon_{i_1+s}^{(j_1+t)}|\mathcal{H}_{i_1, (l, l)}^{(j_1), (l,l - 1)})
\Vert_{M/2}^2}
\end{align*}
Define $t_{i,l}^{(j)} = \mathbf{E}\epsilon_{i}^{(j)}\epsilon_{i+s}^{(j+t)}|\mathcal{H}_{i, (l, l)}^{(j), (l,l)} - \mathbf{E}\epsilon_{i}^{(j)}\epsilon_{i+s}^{(j+t)}|\mathcal{H}_{i, (l, l)}^{(j), (l,l - 1)}$, then from lemma \ref{lemma.second_dimension_difference}
\begin{align*}
\Vert\sum_{i_1 = 1\vee (1 - s)}^{n\wedge (n - s)}c_{i_1,v_1}^{(j_1)}c_{i_1 + s,v_2}^{(j_1 + t)}(\mathbf{E}\epsilon_{i_1}^{(j_1)}\epsilon_{i_1+s}^{(j_1+t)}|\mathcal{H}_{i_1, (l, l)}^{(j_1), (l,l)} - \mathbf{E}\epsilon_{i_1}^{(j_1)}\epsilon_{i_1+s}^{(j_1+t)}|\mathcal{H}_{i_1, (l, l)}^{(j_1), (l,l - 1)})
\Vert_{M/2}\\
\leq C\sqrt{\sum_{i_1 = 1\vee (1 - s)}^{n\wedge (n - s)}c_{i_1,v_1}^{(j_1)2}c_{i_1 + s,v_2}^{(j_1 + t)2}}\times\max_{(i,j)\in B}
\Vert\mathbf{E}t_{i,l}^{(j)}|\mathcal{H}_{i, (0,0)}^{(j), (l,l)}\Vert_{M / 2}\\
+ C\sqrt{\sum_{i_1 = 1\vee (1 - s)}^{n\wedge (n - s)}c_{i_1,v_1}^{(j_1)2}c_{i_1 + s,v_2}^{(j_1 + t)2}}\times\sum_{k = 1}^l\max_{(i,j)\in B}
\Vert\mathbf{E}t_{i,l}^{(j)}|\mathcal{H}_{i, (k,k - 1)}^{(j), (l,l)} - \mathbf{E}t_{i,l}^{(j)}|\mathcal{H}_{i, (k - 1,k - 1)}^{(j), (l,l)}\Vert_{M / 2}\\
+ C\sqrt{\sum_{i_1 = 1\vee (1 - s)}^{n\wedge (n - s)}c_{i_1,v_1}^{(j_1)2}c_{i_1 + s,v_2}^{(j_1 + t)2}}\times\sum_{k = 1}^l\max_{(i,j)\in B}
\Vert\mathbf{E}t_{i,l}^{(j)}|\mathcal{H}_{i, (k,k)}^{(j), (l,l)} - \mathbf{E}t_{i,l}^{(j)}|\mathcal{H}_{i, (k,k - 1)}^{(j), (l,l)}\Vert_{M / 2}
\end{align*}
From lemma \ref{lemma.conditional_prob} and lemma \ref{lemma.delta_moment}
\begin{align*}
\mathbf{E}t_{i,l}^{(j)}|\mathcal{H}_{i, (0,0)}^{(j), (l,l)} = \mathbf{E}\epsilon_{i}^{(j)}\epsilon_{i+s}^{(j+t)}|\mathcal{H}_{i, (0, 0)}^{(j), (l,l)} - \mathbf{E}\epsilon_{i}^{(j)}\epsilon_{i+s}^{(j+t)}|\mathcal{H}_{i, (0, 0)}^{(j), (l,l - 1)}\\
\Rightarrow \Vert\mathbf{E}\epsilon_{i}^{(j)}\epsilon_{i+s}^{(j+t)}|\mathcal{H}_{i, (0, 0)}^{(j), (l,l)} - \mathbf{E}\epsilon_{i}^{(j)}\epsilon_{i+s}^{(j+t)}|\mathcal{H}_{i, (0, 0)}^{(j), (l,l - 1)}\Vert_{M / 2}\\
\leq C\delta_{0}^{(l)} + C\delta_{-s}^{(l-t)}\leq C^\prime(1 + l)^{-\alpha} + C^\prime (1 + \vert s\vert + \vert l - t\vert)^{-\alpha}\\
\Vert
\mathbf{E}t_{i,l}^{(j)}|\mathcal{H}_{i, (k,k - 1)}^{(j), (l,l)} - \mathbf{E}t_{i,l}^{(j)}|\mathcal{H}_{i, (k - 1,k - 1)}^{(j), (l,l)}
\Vert_{M / 2} = \Vert \mathbf{E}\epsilon_{i}^{(j)}\epsilon_{i+s}^{(j+t)}|\mathcal{H}_{i, (k, k - 1)}^{(j), (l,l)}\\
- \mathbf{E}\epsilon_{i}^{(j)}\epsilon_{i+s}^{(j+t)}|\mathcal{H}_{i, (k, k - 1)}^{(j), (l,l - 1)}\Vert_{M / 2}
+ \Vert \mathbf{E}\epsilon_{i}^{(j)}\epsilon_{i+s}^{(j+t)}|\mathcal{H}_{i, (k - 1, k - 1)}^{(j), (l,l)} - \mathbf{E}\epsilon_{i}^{(j)}\epsilon_{i+s}^{(j+t)}|\mathcal{H}_{i, (k - 1, k - 1)}^{(j), (l,l - 1)}\Vert_{M / 2}\\
\leq C\sum_{v = -k}^k \delta_{v}^{(l)} + C\sum_{v = -k}^k \delta_{v - s}^{(l - t)}\leq C^\prime(1 + l)^{1 - \alpha} + C^\prime(1 + \vert l - t\vert)^{1 - \alpha}\\
\text{and } \Vert \mathbf{E}t_{i,l}^{(j)}|\mathcal{H}_{i, (k,k)}^{(j), (l,l)} - \mathbf{E}t_{i,l}^{(j)}|\mathcal{H}_{i, (k,k - 1)}^{(j), (l,l)}\Vert_{M / 2}\leq \Vert \mathbf{E}\epsilon_{i}^{(j)}\epsilon_{i+s}^{(j+t)}|\mathcal{H}_{i, (k, k)}^{(j), (l,l)}\\
- \mathbf{E}\epsilon_{i}^{(j)}\epsilon_{i+s}^{(j+t)}|\mathcal{H}_{i, (k, k)}^{(j), (l,l - 1)}\Vert_{M/2}
+ \Vert \mathbf{E}\epsilon_{i}^{(j)}\epsilon_{i+s}^{(j+t)}|\mathcal{H}_{i, (k, k - 1)}^{(j), (l,l)}
- \mathbf{E}\epsilon_{i}^{(j)}\epsilon_{i+s}^{(j+t)}|\mathcal{H}_{i, (k, k - 1)}^{(j), (l,l - 1)}\Vert_{M/2}\\
\leq C\sum_{v = -k}^k \delta_{v}^{(l)} + C\sum_{v = -k}^k \delta_{v - s}^{(l - t)}\leq C^\prime(1 + l)^{1 - \alpha} + C^\prime(1 + \vert l - t\vert)^{1 - \alpha}\\
\end{align*}
Therefore
\begin{align*}
\Vert
\sum_{i_1 = 1\vee (1 - s)}^{n\wedge (n - s)}\sum_{j_1 = 1\vee (1 - t)}^{m\wedge (m - t)}
c_{i_1,v_1}^{(j_1)}c_{i_1 + s,v_2}^{(j_1 + t)}(\mathbf{E}\epsilon_{i_1}^{(j_1)}\epsilon_{i_1+s}^{(j_1+t)}|\mathcal{H}_{i_1, (l, l)}^{(j_1), (l,l - 1)} - \mathbf{E}\epsilon_{i_1}^{(j_1)}\epsilon_{i_1+s}^{(j_1+t)}|\mathcal{H}_{i_1, (l, l)}^{(j_1), (l - 1,l - 1)})
\Vert_{M/2}\\
\leq C\sqrt{\sum_{i_1 = 1\vee (1 - s)}^{n\wedge (n - s)}\sum_{j_1 = 1\vee (1 - t)}^{m\wedge (m - t)}c_{i_1,v_1}^{(j_1)2}c_{i_1 + s,v_2}^{(j_1 + t)2}}
\times((1 + l)^{2 - \alpha} +(1 + \vert l - t\vert)^{2 - \alpha} + \vert t\vert(1 + \vert l - t\vert)^{1 - \alpha})\\
\text{and } \sum_{l = 1}^\infty\Vert
\sum_{i_1 = 1\vee (1 - s)}^{n\wedge (n - s)}\sum_{j_1 = 1\vee (1 - t)}^{m\wedge (m - t)}
c_{i_1,v_1}^{(j_1)}c_{i_1 + s,v_2}^{(j_1 + t)}(\mathbf{E}\epsilon_{i_1}^{(j_1)}\epsilon_{i_1+s}^{(j_1+t)}|\mathcal{H}_{i_1, (l, l)}^{(j_1), (l,l - 1)} - \mathbf{E}\epsilon_{i_1}^{(j_1)}\epsilon_{i_1+s}^{(j_1+t)}|\mathcal{H}_{i_1, (l, l)}^{(j_1), (l - 1,l - 1)})
\Vert_{M/2}\\
\leq \frac{C}{\mathcal{K}}\times (1 + \vert t\vert)
\end{align*}
and we have
\begin{equation}
\begin{aligned}
\Vert\sum_{i_1 = 1\vee (1 - s)}^{n\wedge (n - s)}\sum_{j_1 = 1\vee (1 - t)}^{m\wedge (m - t)}
c_{i_1,v_1}^{(j_1)}c_{i_1 + s,v_2}^{(j_1 + t)}(\epsilon_{i_1}^{(j_1)}\epsilon_{i_1+s}^{(j_1+t)} - \mathbf{E}\epsilon_{i_1}^{(j_1)}\epsilon_{i_1+s}^{(j_1+t)})
\Vert_{M/2}\leq \frac{C(1 + \vert s\vert + \vert t\vert)}{\mathcal{K}}\\
\text{and }\Vert
\sum_{i_1 = 1}^n\sum_{j_1 = 1}^m\sum_{i_2 = 1}^n\sum_{j_2 = 1}^m c_{i_1,v_1}^{(j_1)}c_{i_2,v_2}^{(j_2)}
\times K\left(\frac{i_1 - i_2}{\mathcal{B}}\right)K\left(\frac{j_1 - j_2}{\mathcal{B}}\right)(\epsilon_{i_1}^{(j_1)}\epsilon_{i_2}^{(j_2)} - \mathbf{E}\epsilon_{i_1}^{(j_1)}\epsilon_{i_2}^{(j_2)})
\Vert_{M/2}\\
\leq \frac{4C}{\mathcal{K}}\sum_{s = 0}^{n - 1}\sum_{t = 0}^{m - 1}K\left(\frac{s}{\mathcal{B}}\right)K\left(\frac{t}{\mathcal{B}}\right)(1 + s + t)
\end{aligned}
\end{equation}
For $\sum_{s = 0}^{n - 1}K\left(\frac{s}{\mathcal{B}}\right)\leq K(0) + \int_{[0,\infty)}K\left(\frac{s}{\mathcal{B}}\right)\mathrm{d}s\leq C\mathcal{B}$ and
\begin{align*}
\sum_{s = 0}^{n - 1} sK\left(\frac{s}{\mathcal{B}}\right)\leq \int_{[0,\infty)}(s+1)K\left(\frac{s}{\mathcal{B}}\right)ds\leq C\mathcal{B}^2
\end{align*}
we have
\begin{equation}
\begin{aligned}
\Vert
\sum_{i_1 = 1}^n\sum_{j_1 = 1}^m\sum_{i_2 = 1}^n\sum_{j_2 = 1}^m c_{i_1,v_1}^{(j_1)}c_{i_2,v_2}^{(j_2)}
\times K\left(\frac{i_1 - i_2}{\mathcal{B}}\right)K\left(\frac{j_1 - j_2}{\mathcal{B}}\right)(\epsilon_{i_1}^{(j_1)}\epsilon_{i_2}^{(j_2)} - \mathbf{E}\epsilon_{i_1}^{(j_1)}\epsilon_{i_2}^{(j_2)})
\Vert_{M/2}\\
\leq \frac{C\mathcal{B}^3}{\mathcal{K}}
\end{aligned}
\end{equation}
In particular,
\begin{equation}
\begin{aligned}
\Vert\ \max_{v_1, v_2 = 1,\cdots, V}\vert
\sum_{i_1 = 1}^n\sum_{j_1 = 1}^m\sum_{i_2 = 1}^n\sum_{j_2 = 1}^m c_{i_1,v_1}^{(j_1)}c_{i_2,v_2}^{(j_2)}\epsilon_{i_1}^{(j_1)}\epsilon_{i_2}^{(j_2)}
\times K\left(\frac{i_1 - i_2}{\mathcal{B}}\right)K\left(\frac{j_1 - j_2}{\mathcal{B}}\right)\\
-\sum_{i_1 = 1}^n\sum_{j_1 = 1}^m\sum_{i_2 = 1}^n\sum_{j_2 = 1}^m c_{i_1,v_1}^{(j_1)}c_{i_2,v_2}^{(j_2)}\mathbf{E}\epsilon_{i_1}^{(j_1)}\epsilon_{i_2}^{(j_2)}
\vert\ \Vert_{M / 2}\\
\leq CS(\mathcal{B}) + V^{4/M}\max_{v_1,v_2 = 1,\cdots,V}\Vert
\sum_{i_1 = 1}^n\sum_{j_1 = 1}^m\sum_{i_2 = 1}^n\sum_{j_2 = 1}^m c_{i_1,v_1}^{(j_1)}c_{i_2,v_2}^{(j_2)}
\times K\left(\frac{i_1 - i_2}{\mathcal{B}}\right)K\left(\frac{j_1 - j_2}{\mathcal{B}}\right)\\
\times(\epsilon_{i_1}^{(j_1)}\epsilon_{i_2}^{(j_2)} - \mathbf{E}\epsilon_{i_1}^{(j_1)}\epsilon_{i_2}^{(j_2)})
\Vert_{M/2}\\
\leq CS(\mathcal{B}) + C\frac{V^{4/M}\mathcal{B}^3}{\mathcal{K}}
\end{aligned}
\end{equation}
and we prove \eqref{eq.variance}.
\end{proof}
\section{Proofs of theorems in section \ref{section.theoretical_justification}}
\begin{proof}[proof of Theorem \ref{theorem.mu}]
For
\begin{align*}
\vert
\frac{1}{T_{n,m}}\sum_{i =  - \mathcal{K}}^{ \mathcal{K}}\sum_{j = - \mathcal{K}}^{\mathcal{K}}X_{i + p_v}^{(j + q_v)}
G\left(\frac{i}{\mathcal{K}}\right) G\left(\frac{j}{\mathcal{K}}\right) - \mu(x_v, y_v)
\vert\\
\leq \frac{1}{T_{n,m}}\vert \sum_{i =  - \mathcal{K}}^{ \mathcal{K}}\sum_{j = - \mathcal{K}}^{\mathcal{K}}
\left(\mu\left(\frac{i + p_v}{n},\frac{j + q_v}{m}\right) - \mu(x_v, y_v)\right)
G\left(\frac{i}{\mathcal{K}}\right) G\left(\frac{j}{\mathcal{K}}\right)\vert\\
+ \frac{1}{T_{n,m}}\vert\sum_{i = -\mathcal{K}}^{\mathcal{K}}\sum_{j = -\mathcal{K}}^{\mathcal{K}}\epsilon_{i + p_v}^{(j + q_v)}
G\left(\frac{i}{\mathcal{K}}\right)G\left(\frac{j}{\mathcal{K}}\right)\vert
\end{align*}
and
\begin{align*}
\frac{1}{T_{n,m}}\vert \sum_{i =  - \mathcal{K}}^{ \mathcal{K}}\sum_{j = - \mathcal{K}}^{\mathcal{K}}
\left(\mu\left(\frac{i + p_v}{n},\frac{j + q_v}{m}\right) - \mu(x_v, y_v)\right)
G\left(\frac{i}{\mathcal{K}}\right) G\left(\frac{j}{\mathcal{K}}\right)\vert\\
\leq C\max_{i =  - \mathcal{K},\cdots, \mathcal{K}}\vert\frac{i + p_v}{n} - x_v\vert + C\max_{j =  - \mathcal{K},\cdots, \mathcal{K}}\vert \frac{j + q_v}{m} - y_v\vert
\end{align*}
here $C$ is a constant satisfying $C\geq \max_{x,y\in[0,1]}\vert\frac{\partial \mu}{\partial x}(x,y)\vert$ and
$C\geq \max_{x,y\in[0,1]}\vert\frac{\partial \mu}{\partial y}(x,y)\vert$. For
\begin{align*}
\vert\frac{i + p_v}{n} - x_v\vert\leq \frac{\vert i\vert}{n} + \vert\frac{p_v - nx_v}{n}\vert\leq \frac{\mathcal{K} + 1}{n}\ \text{and }
\vert \frac{j + q_v}{m} - y_v\vert\leq \frac{\mathcal{K}}{m} + \frac{\vert q_v - my_v\vert }{m}\leq \frac{\mathcal{K} + 1}{m}
\end{align*}
we have
\begin{equation}
\begin{aligned}
\max_{v = 1,\cdots,V}\frac{1}{T_{n,m}}\vert \sum_{i =  - \mathcal{K}}^{ \mathcal{K}}\sum_{j = - \mathcal{K}}^{\mathcal{K}}
\left(\mu\left(\frac{i + p_v}{n},\frac{j + q_v}{m}\right) - \mu(x_v, y_v)\right)
G\left(\frac{i}{\mathcal{K}}\right) G\left(\frac{j}{\mathcal{K}}\right)\vert\\
 = O\left(\frac{\mathcal{K}}{n} + \frac{\mathcal{K}}{m}\right)
\end{aligned}
\label{eq.mean_field_difference}
\end{equation}
On the other hand, from lemma \ref{lemma.consistent_linear_combination},
\begin{equation}
\begin{aligned}
\Vert\sum_{i = -\mathcal{K}}^{\mathcal{K}}\sum_{j = -\mathcal{K}}^{\mathcal{K}}\epsilon_{i + p_v}^{(j + q_v)}
G\left(\frac{i}{\mathcal{K}}\right)G\left(\frac{j}{\mathcal{K}}\right)\Vert_M\leq C\sqrt{\sum_{i = -\mathcal{K}}^{\mathcal{K}}\sum_{j = -\mathcal{K}}^{\mathcal{K}}
G^2\left(\frac{i}{\mathcal{K}}\right)G^2\left(\frac{j}{\mathcal{K}}\right)
}
\end{aligned}
\end{equation}
For
\begin{align*}
c\mathcal{K}\leq\int_{[0,\mathcal{K}]}G(x/\mathcal{K})dx\leq \sum_{i = -\mathcal{K}}^\mathcal{K}G\left(\frac{i}{\mathcal{K}}\right)\leq 1 + 2\int_{[0,\mathcal{K}]}G(x/\mathcal{K})dx\leq C\mathcal{K}
\end{align*}
with $c,C>0$ being constants. Similarly we have
\begin{align*}
c\mathcal{K}\leq \int_{[0,\mathcal{K}]}G^2(x/\mathcal{K})dx\leq \sum_{i = -\mathcal{K}}^\mathcal{K}G^2\left(\frac{i}{\mathcal{K}}\right)\leq 1 + 2\int_{[0,\mathcal{K}]}G^2(x/\mathcal{K})dx\leq C\mathcal{K}
\end{align*}
Therefore $\Vert\frac{1}{T_{n,m}}\sum_{i = -\mathcal{K}}^{\mathcal{K}}\sum_{j = -\mathcal{K}}^{\mathcal{K}}\epsilon_{i + p_v}^{(j + q_v)}
G\left(\frac{i}{\mathcal{K}}\right)G\left(\frac{j}{\mathcal{K}}\right)\Vert_M = O\left(\frac{1}{\mathcal{K}}\right)$, and
\begin{equation}
\begin{aligned}
\Vert\ \max_{v = 1,\cdots, V}\vert\widehat{\mu}(x_v, y_v) - \mu(x_v, y_v)\vert\ \Vert_M = O\left(\frac{\mathcal{K}}{n} + \frac{\mathcal{K}}{m} + \frac{V^{1/M}}{\mathcal{K}}\right)
\end{aligned}
\end{equation}

Define $c_{i,v}^{(j)}$ as in algorithm \ref{algorithm.bootstrap}, then $\sum_{i = 1}^n\sum_{j = 1}^m c_{i,v}^{(j)2} = 1$. Besides,
\begin{align*}
\vert c_{i,v}^{(j)}\vert = \frac{G\left(\frac{i - p_v}{\mathcal{K}}\right)}{\sqrt{\sum_{i = p_v - \mathcal{K}}^{p_v + \mathcal{K}}G^2\left(\frac{i - p_v}{\mathcal{K}}\right)}}\times \frac{G\left(\frac{j - q_v}{\mathcal{K}}\right)}{\sqrt{\sum_{j = q_v - \mathcal{K}}^{q_v + \mathcal{K}}G^2\left(\frac{j - q_v}{\mathcal{K}}\right)}}
\end{align*}
and $\frac{G\left(\frac{i - p_v}{\mathcal{K}}\right)}{\sqrt{\sum_{i = p_v - \mathcal{K}}^{p_v + \mathcal{K}}G^2\left(\frac{i - p_v}{\mathcal{K}}\right)}}\leq \frac{C}{\sqrt{\mathcal{K}}}$ with a constant $C$. From theorem \ref{theorem.Gaussian}, define $R_v = \sum_{i  =1}^n\sum_{j = 1}^m c_{i,v}^{(j)}\epsilon_i^{(j)}$, then
\begin{equation}
\sup_{x\in\mathbf{R}}\vert Prob\left(\max_{v = 1,\cdots, V}\vert R_v\vert\leq x\right)
- Prob\left(\max_{v = 1,\cdots, V}\vert \xi_v\vert\leq x\right)\vert = o(1)
\end{equation}
Besides,
\begin{align*}
\vert\frac{T_{n,m}}{B_{n,m}}(\widehat{\mu}(x_v, y_v) - \mu(x_v, y_v)) - R_v\vert\\
\leq \frac{1}{B_{n,m}}\sum_{i = -\mathcal{K}}^\mathcal{K}\sum_{j = -\mathcal{K}}^\mathcal{K}\vert\mu\left(\frac{i + p_v}{n},\frac{j + q_v}{m}\right) - \mu(x_v, y_v)\vert\times
G\left(\frac{i}{\mathcal{K}}\right)G\left(\frac{j}{\mathcal{K}}\right)\\
\leq C\frac{T_{n,m}}{B_{n,m}}\times\left(\frac{\mathcal{K}}{n} + \frac{\mathcal{K}}{m}\right) = O\left(\frac{\mathcal{K}^2}{n} + \frac{\mathcal{K}^2}{m}\right)
\end{align*}
Define $\Delta = C(\frac{\mathcal{K}^2}{n} + \frac{\mathcal{K}^2}{m})$ for a sufficiently large constant $C$. From lemma \ref{lemma.Gaussian_property} in the appendix and assumption 4, we know that $\mathbf{E}\xi^2_v\geq c > 0$ for a constant $c$ and
\begin{equation}
\begin{aligned}
Prob\left(\max_{v = 1,\cdots, V}\vert\frac{T_{n,m}}{B_{n,m}}(\widehat{\mu}(x_v, y_v) - \mu(x_v, y_v))\vert\leq x\right)
- Prob\left(\max_{v = 1,\cdots, V}\vert\xi_v\vert\leq x\right)\\
\leq  Prob\left(\max_{v = 1,\cdots, V}\vert R_v\vert\leq x + \Delta\right) - Prob\left(\max_{v = 1,\cdots, V}\vert\xi_v\vert\leq x + \Delta\right)\\
+ C\Delta\left(1 + \sqrt{\log(V)} + \sqrt{\vert\log(\Delta)\vert}\right)\\
\text{and }Prob\left(\max_{v = 1,\cdots, V}\vert\frac{T_{n,m}}{B_{n,m}}(\widehat{\mu}(x_v, y_v) - \mu(x_v, y_v))\vert\leq x\right)
- Prob\left(\max_{v = 1,\cdots, V}\vert\xi_v\vert\leq x\right)\\
\geq Prob\left(\max_{v = 1,\cdots, V}\vert R_v\vert\leq x - \Delta\right) - Prob\left(\max_{v = 1,\cdots, V}\vert\xi_v\vert\leq x - \Delta\right)\\
- C\Delta\left(1 + \sqrt{\log(V)} + \sqrt{\vert\log(\Delta)\vert}\right)\\
\Rightarrow \sup_{x\in\mathbf{R}}\vert Prob\left(\max_{v = 1,\cdots, V}\vert\frac{T_{n,m}}{B_{n,m}}(\widehat{\mu}(x_v, y_v) - \mu(x_v, y_v))\vert\leq x\right)
- Prob\left(\max_{v = 1,\cdots, V}\vert\xi_v\vert\leq x\right)\vert\to 0
\end{aligned}
\end{equation}
and we prove eq.\eqref{eq.Gaussian_mu}.

For the second result, define $\tau_v = \sigma_v^{1/3}$, from assumption 4,
\begin{align*}
\sigma^2_v = \sum_{i_1 = 1}^n\sum_{j_1 = 1}^m\sum_{i_2 = 1}^n\sum_{j_2 = 1}^m c_{i_1, v}^{(j_1)}c_{i_2, v}^{(j_2)}\mathbf{E}\epsilon_{i_1}^{(j_1)}\epsilon_{i_2}^{(j_2)}\geq c\sum_{i = 1}^n\sum_{j = 1}^m  c_{i, v}^{(j)2} = c\\
\text{and } \sigma_v\leq \Vert \sum_{i = 1}^n\sum_{j = 1}^m c_{i,v}^{(j)}\epsilon_{i}^{(j)}\Vert_M\leq C\sqrt{\sum_{i = 1}^n\sum_{j = 1}^m  c_{i, v}^{(j)2}} = C
\end{align*}
so $0 < (\sqrt{c})^{1/3}\leq \tau_v = \sigma_v^{1/3}\leq C^{1/3}$. Also notice that
\begin{align*}
\frac{T_{n,m}}{\tau_v B_{n,m}}\left(\widehat{\mu}(x_v, y_v) - \mu(x_v, y_v)\right)\\
= \frac{1}{\tau_v B_{n,m}}\sum_{i = -\mathcal{K}}^{\mathcal{K}}\sum_{j = -\mathcal{K}}^{\mathcal{K}}(\mu\left(\frac{i + p_v}{n}, \frac{j + q_v}{m}\right) - \mu(x_v, y_v))G\left(\frac{i}{\mathcal{K}}\right)G\left(\frac{j}{\mathcal{K}}\right)\\
+ \frac{1}{\tau_v B_{n,m}}\sum_{i = -\mathcal{K}}^{\mathcal{K}}\sum_{j = -\mathcal{K}}^{\mathcal{K}}\epsilon_{i + p_v}^{(j + q_v)}G\left(\frac{i}{\mathcal{K}}\right) G\left(\frac{j}{\mathcal{K}}\right)
\end{align*}
From theorem \ref{theorem.Gaussian},
\begin{equation}
\begin{aligned}
\sup_{x\in\mathbf{R}}\vert Prob\left(\max_{v = 1,\cdots, V}
\frac{\vert \sum_{i = -\mathcal{K}}^{\mathcal{K}}\sum_{j = -\mathcal{K}}^{\mathcal{K}}\epsilon_{i + p_v}^{(j + q_v)}G\left(\frac{i}{\mathcal{K}}\right) G\left(\frac{j}{\mathcal{K}}\right)\vert}{\tau_v B_{n, m}}\leq x\right)\\
- Prob\left(\max_{v = 1,\cdots, V}\vert\frac{\xi_v}{\sigma^{1/3}_v}\vert\leq x\right)\vert = o(1)
\end{aligned}
\end{equation}
From eq.\eqref{eq.mean_field_difference}, define $\delta_v = \frac{1}{\tau_v B_{n,m}}\sum_{i = -\mathcal{K}}^{\mathcal{K}}\sum_{j = -\mathcal{K}}^{\mathcal{K}}(\mu\left(\frac{i + p_v}{n}, \frac{j + q_v}{m}\right) - \mu(x_v, y_v))G\left(\frac{i}{\mathcal{K}}\right)G\left(\frac{j}{\mathcal{K}}\right)$ and $\delta = \max_{v = 1,\cdots, V}\vert\delta_v\vert $, we have
$\delta = O\left(\frac{\mathcal{K}^2}{n} + \frac{\mathcal{K}^2}{m}\right)$, and from lemma \ref{lemma.Gaussian_property},
\begin{equation}
\begin{aligned}
\sup_{x\in\mathbf{R}}\vert
Prob\left(\max_{v = 1,\cdots, V}\vert \frac{T_{n,m}}{\tau_v B_{n,m}}\left(\widehat{\mu}(x_v, y_v) - \mu(x_v, y_v)\right)\vert\leq x\right)
- Prob\left(\max_{v = 1,\cdots, V}\vert\frac{\xi_v}{\sigma^{1/3}_v}\vert\leq x\right)
\vert\\
\leq \sup_{x\in\mathbf{R}}\vert Prob\left(\max_{v = 1,\cdots, V}
\frac{\vert \sum_{i = -\mathcal{K}}^{\mathcal{K}}\sum_{j = -\mathcal{K}}^{\mathcal{K}}\epsilon_{i + p_v}^{(j + q_v)}G\left(\frac{i}{\mathcal{K}}\right) G\left(\frac{j}{\mathcal{K}}\right)\vert}{\tau_v B_{n, m}}\leq x\right)\\
- Prob\left(\max_{v = 1,\cdots, V}\vert\frac{\xi_v}{\sigma^{1/3}_v}\vert\leq x\right)\vert\\
+ C\delta(1 + \sqrt{\log(V)} +\sqrt{\vert\log(\delta)\vert}) = o(1)
\end{aligned}
\end{equation}
Form eq.\eqref{eq.delta_hat},
\begin{align*}
\max_{v = 1,\cdots, V}\vert\widehat{\sigma}^2_v - \sigma^2_v\vert = O_p(\mathcal{S}(\mathcal{B}) + \frac{V^{4/M}\mathcal{B}^3}{\mathcal{K}} + \frac{V^{2/M}\mathcal{BK}}{n}
+\frac{V^{2/M}\mathcal{BK}}{m})\\
\Rightarrow \max_{v = 1,\cdots, V}\vert (\widehat{\sigma}^2)^{1/6} - (\sigma_v^2)^{1/6}\vert = O_p(\mathcal{S}(\mathcal{B}) + \frac{V^{4/M}\mathcal{B}^3}{\mathcal{K}} + \frac{V^{2/M}\mathcal{BK}}{n}
+\frac{V^{2/M}\mathcal{BK}}{m})
\end{align*}
so for sufficiently large $n,m$, $\widehat{\sigma}_v > \sqrt{c}/2$ and $\widehat{\sigma}_v^{1/3}\geq (\sqrt{c}/2)^{1/3}$ with probability close to 1, and correspondingly from theorem \ref{theorem.Gaussian}
\begin{align*}
\vert\ \max_{v = 1,\cdots, V}\frac{T_{n,m}}{\widehat{\tau}_v B_{n,m}}\vert\widehat{\mu}(x_v, y_v) - \mu(x_v, y_v)\vert - \max_{v = 1,\cdots, V}\vert \frac{T_{n,m}}{\tau_v B_{n,m}}\left(\widehat{\mu}(x_v, y_v) - \mu(x_v, y_v)\right)\vert \ \vert\\
\leq \max_{v = 1,\cdots, V}\frac{T_{n,m}}{B_{n,m}}\vert\widehat{\mu}(x_v, y_v) - \mu(x_v, y_v)\vert\times \max_{v = 1,\cdots, V}\frac{2\vert\widehat{\tau}_v - \tau_v\vert}{c^\prime}\\
= O_p\left(V^{1/M}\mathcal{S}(\mathcal{B}) + \frac{V^{5/M}\mathcal{B}^3}{\mathcal{K}} + \frac{V^{3/M}\mathcal{BK}}{n}
+\frac{V^{3/M}\mathcal{BK}}{m}\right)
\end{align*}
Set $\lambda = V^{1/M}\mathcal{S}(\mathcal{B}) + \frac{V^{5/M}\mathcal{B}^3}{\mathcal{K}} + \frac{V^{3/M}\mathcal{BK}}{n}
+\frac{V^{3/M}\mathcal{BK}}{m}$, then for any given $0< \phi < 1$, $\exists $ a constant $C_\phi$ such that
\begin{align*}
Prob\left(\vert\ \max_{v = 1,\cdots, V}\frac{T_{n,m}}{\widehat{\tau}_v B_{n,m}}\vert\widehat{\mu}(x_v, y_v) - \mu(x_v, y_v)\vert - \max_{v = 1,\cdots, V}\vert \frac{T_{n,m}}{\tau_v B_{n,m}}\left(\widehat{\mu}(x_v, y_v) - \mu(x_v, y_v)\right)\vert \ \vert > C_\phi\lambda\right) < \phi
\end{align*}
and
\begin{align*}
Prob\left(\max_{v = 1,\cdots, V}\frac{T_{n,m}}{\widehat{\tau}_v B_{n,m}}\vert\widehat{\mu}(x_v, y_v) - \mu(x_v, y_v)\vert\leq x\right)\\
\leq \phi
 + Prob\left(\max_{v = 1,\cdots, V}\vert \frac{T_{n,m}}{\tau_v B_{n,m}}\left(\widehat{\mu}(x_v, y_v) - \mu(x_v, y_v)\right)\vert\leq x + C_\phi\lambda\right)\\
 Prob\left(\max_{v = 1,\cdots, V}\frac{T_{n,m}}{\widehat{\tau}_v B_{n,m}}\vert\widehat{\mu}(x_v, y_v) - \mu(x_v, y_v)\vert\leq x\right)\\
\geq -\phi
 + Prob\left(\max_{v = 1,\cdots, V}\vert \frac{T_{n,m}}{\tau_v B_{n,m}}\left(\widehat{\mu}(x_v, y_v) - \mu(x_v, y_v)\right)\vert\leq x - C_\phi\lambda\right)
\end{align*}
From lemma \ref{lemma.Gaussian_property}, we prove  \eqref{eq.Gaussian_weight}.
\end{proof}

\begin{proof}[proof of theorem \ref{theorem.bootstrap}]
In Algorithm \ref{algorithm.bootstrap}, define 
\begin{align*}
   T_v^* =  \frac{T_{n,m}}{B_{n,m}}(\widehat{\mu}^*(x_v, y_v) - \widehat{\mu}(x_v, y_v)) = \sum_{i = 1}^n\sum_{j = 1}^m c_{i,v}^{(j)}\epsilon_i^{(j)*},
\end{align*}
then $T_v^*$ are joint normal random variables in the bootstrap world. Moreover,
\begin{align*}
\mathbf{E}^*T_v^* =  \sum_{i = 1}^n\sum_{j = 1}^m c_{i,v}^{(j)}\mathbf{E}^*\epsilon_i^{(j)*} = 0,\\
\text{and } \mathbf{E}^*T_{v_1}^*T_{v_2}^*  = \sum_{i_1 = 1}^n\sum_{j_1 = 1}^m \sum_{i_2 = 1}^n\sum_{j_2 = 1}^m c_{i_1,v_1}^{(j_1)}c_{i_2,v_2}^{(j_2)}\mathbf{E}^*\epsilon_{i_1}^{(j_1)*}\epsilon_{i_2}^{(j_2)*}\\
= \sum_{i_1 = p_{v_1} - \mathcal{K}}^{p_{v_1} + \mathcal{K}}\sum_{j_1 = q_{v_1} - \mathcal{K}}^{q_{v_1} + \mathcal{K}}
\sum_{i_2 = p_{v_2} - \mathcal{K}}^{p_{v_2} + \mathcal{K}}\sum_{j_2 = q_{v_2} - \mathcal{K}}^{q_{v_2} + \mathcal{K}}
c_{i_1, v_1}^{(j_1)}c_{i_2, v_2}^{(j_2)}\widehat{\epsilon}_{i_1}^{(j_1)}\widehat{\epsilon}_{i_2}^{(j_2)}K\left(\frac{i_1 - i_2}{\mathcal{B}}\right)
K\left(\frac{j_1 - j_2}{\mathcal{B}}\right).
\end{align*}
According to the definition of $c_{i,v}^{(j)}$, i.e., eq.\eqref{eq.def_cij}, define $\Delta_i^{(j)} = \widehat{\epsilon}_{i}^{(j)} - \epsilon_i^{(j)} = \mu\left(\frac{i}{n}, \frac{j}{m}\right) - \widehat{\mu}\left(\frac{i}{n}, \frac{j}{m}\right)$ for $i = \mathcal{K} + 1,\cdots, n - \mathcal{K}$ and $j = \mathcal{K} + 1,\cdots, m - \mathcal{K}$. Notice that $c_{i,v}^{(j)} = 0$ for $i\leq \mathcal{K}$ or $i > n - \mathcal{K}$ or $j\leq \mathcal{K}$ or $j > m - \mathcal{K}$, so we have
\begin{align*}
\vert \mathbf{E}^*T_{v_1}^*T_{v_2}^* - \mathbf{E}\xi_{v_1}\xi_{v_2}\vert\leq
\vert
\sum_{i_1 = 1}^n\sum_{j_1 = 1}^m\sum_{i_2 = 1}^n\sum_{j_2 = 1}^m c_{i_1, v_1}^{(j_1)}c_{i_2, v_2}^{(j_2)}
\epsilon_{i_1}^{(j_1)}\epsilon_{i_2}^{(j_2)}K\left(\frac{i_1 - i_2}{\mathcal{B}}\right)
K\left(\frac{j_1 - j_2}{\mathcal{B}}\right)\\
-\sum_{i_1 = 1}^n\sum_{j_1 = 1}^m\sum_{i_2 = 1}^n\sum_{j_2 = 1}^m c_{i_1, v_1}^{(j_1)}c_{i_2, v_2}^{(j_2)}
\mathbf{E}\epsilon_{i_1}^{(j_1)}\epsilon_{i_2}^{(j_2)}
\vert\\
+\vert
\sum_{i_1 = p_{v_1} - \mathcal{K}}^{p_{v_1} + \mathcal{K}}\sum_{j_1 = q_{v_1} - \mathcal{K}}^{q_{v_1} + \mathcal{K}}
\sum_{i_2 = p_{v_2} - \mathcal{K}}^{p_{v_2} + \mathcal{K}}\sum_{j_2 = q_{v_2} - \mathcal{K}}^{q_{v_2} + \mathcal{K}}
c_{i_1, v_1}^{(j_1)}c_{i_2, v_2}^{(j_2)}\Delta_{i_1}^{(j_1)}\epsilon_{i_2}^{(j_2)}K\left(\frac{i_1 - i_2}{\mathcal{B}}\right)
K\left(\frac{j_1 - j_2}{\mathcal{B}}\right)
\vert\\
+\vert
\sum_{i_1 = p_{v_1} - \mathcal{K}}^{p_{v_1} + \mathcal{K}}\sum_{j_1 = q_{v_1} - \mathcal{K}}^{q_{v_1} + \mathcal{K}}
\sum_{i_2 = p_{v_2} - \mathcal{K}}^{p_{v_2} + \mathcal{K}}\sum_{j_2 = q_{v_2} - \mathcal{K}}^{q_{v_2} + \mathcal{K}}
c_{i_1, v_1}^{(j_1)}c_{i_2, v_2}^{(j_2)}\epsilon_{i_1}^{(j_1)}\Delta_{i_2}^{(j_2)}K\left(\frac{i_1 - i_2}{\mathcal{B}}\right)
K\left(\frac{j_1 - j_2}{\mathcal{B}}\right)
\vert\\
+\vert
\sum_{i_1 = p_{v_1} - \mathcal{K}}^{p_{v_1} + \mathcal{K}}\sum_{j_1 = q_{v_1} - \mathcal{K}}^{q_{v_1} + \mathcal{K}}
\sum_{i_2 = p_{v_2} - \mathcal{K}}^{p_{v_2} + \mathcal{K}}\sum_{j_2 = q_{v_2} - \mathcal{K}}^{q_{v_2} + \mathcal{K}}
c_{i_1, v_1}^{(j_1)}c_{i_2, v_2}^{(j_2)}\Delta_{i_1}^{(j_1)}\Delta_{i_2}^{(j_2)}K\left(\frac{i_1 - i_2}{\mathcal{B}}\right)
K\left(\frac{j_1 - j_2}{\mathcal{B}}\right)
\vert
\end{align*}
From theorem \ref{theorem.covariances},
\begin{equation}
\begin{aligned}
\max_{v_1,v_2 = 1,\cdots, V}\vert
\sum_{i_1 = 1}^n\sum_{j_1 = 1}^m\sum_{i_2 = 1}^n\sum_{j_2 = 1}^m c_{i_1, v_1}^{(j_1)}c_{i_2, v_2}^{(j_2)}
\epsilon_{i_1}^{(j_1)}\epsilon_{i_2}^{(j_2)}K\left(\frac{i_1 - i_2}{\mathcal{B}}\right)
K\left(\frac{j_1 - j_2}{\mathcal{B}}\right)\\
-\sum_{i_1 = 1}^n\sum_{j_1 = 1}^m\sum_{i_2 = 1}^n\sum_{j_2 = 1}^m c_{i_1, v_1}^{(j_1)}c_{i_2, v_2}^{(j_2)}
\mathbf{E}\epsilon_{i_1}^{(j_1)}\epsilon_{i_2}^{(j_2)}
\vert = O_p\left(S(\mathcal{B}) + \frac{V^{4/M}\mathcal{B}^3}{\mathcal{K}}\right)
\end{aligned}
\end{equation}
with $S(\mathcal{B})$ being defined in theorem \ref{theorem.covariances}. On the other hand, by definition
\begin{align*}
\Delta_i^{(j)} =
\frac{1}{T_{n,m}}\sum_{s = -\mathcal{K}}^{\mathcal{K}}\sum_{t = -\mathcal{K}}^{\mathcal{K}}
\left(\mu\left(\frac{i}{n}, \frac{j}{m}\right) - \mu\left(\frac{i + s}{n}, \frac{j + t}{m}\right)\right)
G\left(\frac{s}{\mathcal{K}}\right)G\left(\frac{t}{\mathcal{K}}\right)\\
-\frac{1}{T_{n,m}}\sum_{s = -\mathcal{K}}^{\mathcal{K}}\sum_{t = -\mathcal{K}}^{\mathcal{K}}\epsilon_{i + s}^{(j + t)}G\left(\frac{s}{\mathcal{K}}\right)G\left(\frac{t}{\mathcal{K}}\right)
\end{align*}
Define $\eta_{i}^{(j)} = \frac{1}{T_{n,m}}\sum_{s = -\mathcal{K}}^{\mathcal{K}}\sum_{t = -\mathcal{K}}^{\mathcal{K}}
\left(\mu\left(\frac{i}{n}, \frac{j}{m}\right) - \mu\left(\frac{i + s}{n}, \frac{j + t}{m}\right)\right)
G\left(\frac{s}{\mathcal{K}}\right)G\left(\frac{t}{\mathcal{K}}\right)$ and

\noindent$\zeta_i^{(j)} = \frac{1}{T_{n,m}}\sum_{s = -\mathcal{K}}^{\mathcal{K}}\sum_{t = -\mathcal{K}}^{\mathcal{K}}\epsilon_{i + s}^{(j + t)}G\left(\frac{s}{\mathcal{K}}\right)G\left(\frac{t}{\mathcal{K}}\right)$, then $\Delta_i^{(j)} = \eta_{i}^{(j)} - \zeta_i^{(j)}$; and
\begin{align*}
\vert
\sum_{i_1 = p_{v_1} - \mathcal{K}}^{p_{v_1} + \mathcal{K}}\sum_{j_1 = q_{v_1} - \mathcal{K}}^{q_{v_1} + \mathcal{K}}
\sum_{i_2 = p_{v_2} - \mathcal{K}}^{p_{v_2} + \mathcal{K}}\sum_{j_2 = q_{v_2} - \mathcal{K}}^{q_{v_2} + \mathcal{K}}
c_{i_1, v_1}^{(j_1)}c_{i_2, v_2}^{(j_2)}\Delta_{i_1}^{(j_1)}\epsilon_{i_2}^{(j_2)}K\left(\frac{i_1 - i_2}{\mathcal{B}}\right)
K\left(\frac{j_1 - j_2}{\mathcal{B}}\right)
\vert\\
\leq  \vert
\sum_{i_1 = p_{v_1} - \mathcal{K}}^{p_{v_1} + \mathcal{K}}\sum_{j_1 = q_{v_1} - \mathcal{K}}^{q_{v_1} + \mathcal{K}}
\sum_{i_2 = p_{v_2} - \mathcal{K}}^{p_{v_2} + \mathcal{K}}\sum_{j_2 = q_{v_2} - \mathcal{K}}^{q_{v_2} + \mathcal{K}}
c_{i_1, v_1}^{(j_1)}c_{i_2, v_2}^{(j_2)}\eta_{i_1}^{(j_1)}\epsilon_{i_2}^{(j_2)}K\left(\frac{i_1 - i_2}{\mathcal{B}}\right)
K\left(\frac{j_1 - j_2}{\mathcal{B}}\right)
\vert\\
+\vert
\sum_{i_1 = p_{v_1} - \mathcal{K}}^{p_{v_1} + \mathcal{K}}\sum_{j_1 = q_{v_1} - \mathcal{K}}^{q_{v_1} + \mathcal{K}}
\sum_{i_2 = p_{v_2} - \mathcal{K}}^{p_{v_2} + \mathcal{K}}\sum_{j_2 = q_{v_2} - \mathcal{K}}^{q_{v_2} + \mathcal{K}}
c_{i_1, v_1}^{(j_1)}c_{i_2, v_2}^{(j_2)}\zeta_{i_1}^{(j_1)}\epsilon_{i_2}^{(j_2)}K\left(\frac{i_1 - i_2}{\mathcal{B}}\right)
K\left(\frac{j_1 - j_2}{\mathcal{B}}\right)
\vert
\end{align*}
From lemma \ref{lemma.consistent_linear_combination},
\begin{align*}
\Vert
\sum_{i_1 = p_{v_1} - \mathcal{K}}^{p_{v_1} + \mathcal{K}}\sum_{j_1 = q_{v_1} - \mathcal{K}}^{q_{v_1} + \mathcal{K}}
\sum_{i_2 = p_{v_2} - \mathcal{K}}^{p_{v_2} + \mathcal{K}}\sum_{j_2 = q_{v_2} - \mathcal{K}}^{q_{v_2} + \mathcal{K}}
c_{i_1, v_1}^{(j_1)}c_{i_2, v_2}^{(j_2)}\eta_{i_1}^{(j_1)}\epsilon_{i_2}^{(j_2)}K\left(\frac{i_1 - i_2}{\mathcal{B}}\right)
K\left(\frac{j_1 - j_2}{\mathcal{B}}\right)
\Vert_M\\
\leq C\sqrt{
\sum_{i_2 = p_{v_2} - \mathcal{K}}^{p_{v_2} + \mathcal{K}}\sum_{j_2 = q_{v_2} - \mathcal{K}}^{q_{v_2} + \mathcal{K}}
c_{i_2, v_2}^{(j_2)2}\left(
\sum_{i_1 = p_{v_1} - \mathcal{K}}^{p_{v_1} + \mathcal{K}}\sum_{j_1 = q_{v_1} - \mathcal{K}}^{q_{v_1} + \mathcal{K}}
c_{i_1, v_1}^{(j_1)}\eta_{i_1}^{(j_1)}K\left(\frac{i_1 - i_2}{\mathcal{B}}\right)
K\left(\frac{j_1 - j_2}{\mathcal{B}}\right)
\right)^2
}
\end{align*}
Notice that $\vert\eta_i^{(j)}\vert\leq \frac{C\mathcal{K}}{n} + \frac{C\mathcal{K}}{m}$ for a constant $C$, from Cauchy inequality
\begin{align*}
\left(
\sum_{i_1 = p_{v_1} - \mathcal{K}}^{p_{v_1} + \mathcal{K}}\sum_{j_1 = q_{v_1} - \mathcal{K}}^{q_{v_1} + \mathcal{K}}
c_{i_1, v_1}^{(j_1)}\eta_{i_1}^{(j_1)}K\left(\frac{i_1 - i_2}{\mathcal{B}}\right)
K\left(\frac{j_1 - j_2}{\mathcal{B}}\right)
\right)^2
\leq \left(\sum_{i_1 = p_{v_1} - \mathcal{K}}^{p_{v_1} + \mathcal{K}}\sum_{j_1 = q_{v_1} - \mathcal{K}}^{q_{v_1} + \mathcal{K}}
c_{i_1, v_1}^{(j_1)2}\right)\\
\times\left(\sum_{i_1 = p_{v_1} - \mathcal{K}}^{p_{v_1} + \mathcal{K}}\sum_{j_1 = q_{v_1} - \mathcal{K}}^{q_{v_1} + \mathcal{K}}
\eta_{i_1}^{(j_1)2}K^2\left(\frac{i_1 - i_2}{\mathcal{B}}\right)
K^2\left(\frac{j_1 - j_2}{\mathcal{B}}\right)
\right)\\
\leq C\left(\frac{\mathcal{K}^2}{n^2} + \frac{\mathcal{K}^2}{m^2}\right)\times\left(\sum_{i = -\infty}^\infty K^2\left(\frac{i}{\mathcal{B}}\right)\right)^2
\end{align*}
For $K$ is decreasing on $[0,\infty)$,
$\sum_{i = 0}^\infty K^2\left(\frac{i}{\mathcal{B}}\right)\leq 1 + \int_{0}^\infty K^2\left(\frac{x}{\mathcal{B}}\right)\mathrm{d}x\leq C\mathcal{B}$, so
\begin{equation}
\begin{aligned}
\Vert
\sum_{i_1 = p_{v_1} - \mathcal{K}}^{p_{v_1} + \mathcal{K}}\sum_{j_1 = q_{v_1} - \mathcal{K}}^{q_{v_1} + \mathcal{K}}
\sum_{i_2 = p_{v_2} - \mathcal{K}}^{p_{v_2} + \mathcal{K}}\sum_{j_2 = q_{v_2} - \mathcal{K}}^{q_{v_2} + \mathcal{K}}
c_{i_1, v_1}^{(j_1)}c_{i_2, v_2}^{(j_2)}\eta_{i_1}^{(j_1)}\epsilon_{i_2}^{(j_2)}K\left(\frac{i_1 - i_2}{\mathcal{B}}\right)
K\left(\frac{j_1 - j_2}{\mathcal{B}}\right)
\Vert_M\\
= O\left(\frac{\mathcal{B}\mathcal{K}}{n} + \frac{\mathcal{B}\mathcal{K}}{m}\right)\\
\Rightarrow \max_{v_1,v_2 = 1,\cdots, V}\vert \sum_{i_1 = p_{v_1} - \mathcal{K}}^{p_{v_1} + \mathcal{K}}\sum_{j_1 = q_{v_1} - \mathcal{K}}^{q_{v_1} + \mathcal{K}}
\sum_{i_2 = p_{v_2} - \mathcal{K}}^{p_{v_2} + \mathcal{K}}\sum_{j_2 = q_{v_2} - \mathcal{K}}^{q_{v_2} + \mathcal{K}}
c_{i_1, v_1}^{(j_1)}c_{i_2, v_2}^{(j_2)}\eta_{i_1}^{(j_1)}\epsilon_{i_2}^{(j_2)}\\
\times K\left(\frac{i_1 - i_2}{\mathcal{B}}\right)
K\left(\frac{j_1 - j_2}{\mathcal{B}}\right)\vert
= O_p\left(\frac{V^{2/M}\times \mathcal{B}\mathcal{K}}{n} + \frac{V^{2/M}\times\mathcal{B}\mathcal{K}}{m}\right)
\end{aligned}
\end{equation}
On the other hand, for $\left(K\left(\frac{i_1 - i_2}{\mathcal{B}}\right)\right)_{i_1, i_2 = 1,\cdots, n}$ and
$\left(K\left(\frac{j_1 - j_2}{\mathcal{B}}\right)\right)_{j_1, j_2 = 1,\cdots, m}$ are Toeplitz matrix and
$\sum_{k = -\infty}^\infty K\left(\frac{k}{\mathcal{B}}\right)\leq 1 + 2\int_{0}^\infty K\left(\frac{x}{\mathcal{B}}\right) \mathrm{d}x\leq C\mathcal{B}$, from section 0.9.7 in \cite{MR2978290}
\begin{align*}
\vert\sum_{i_1 = p_{v_1} - \mathcal{K}}^{p_{v_1} + \mathcal{K}}\sum_{j_1 = q_{v_1} - \mathcal{K}}^{q_{v_1} + \mathcal{K}}
\sum_{i_2 = p_{v_2} - \mathcal{K}}^{p_{v_2} + \mathcal{K}}\sum_{j_2 = q_{v_2} - \mathcal{K}}^{q_{v_2} + \mathcal{K}}
c_{i_1, v_1}^{(j_1)}c_{i_2, v_2}^{(j_2)}\zeta_{i_1}^{(j_1)}\epsilon_{i_2}^{(j_2)}K\left(\frac{i_1 - i_2}{\mathcal{B}}\right)
K\left(\frac{j_1 - j_2}{\mathcal{B}}\right)\vert\\
\leq
\sum_{i_1 = p_{v_1} - \mathcal{K}}^{p_{v_1} + \mathcal{K}}\sum_{i_2 = p_{v_2} - \mathcal{K}}^{p_{v_2} + \mathcal{K}}K\left(\frac{i_1 - i_2}{\mathcal{B}}\right)
\vert\sum_{j_1 = q_{v_1} - \mathcal{K}}^{q_{v_1} + \mathcal{K}}\sum_{j_2 = q_{v_2} - \mathcal{K}}^{q_{v_2} + \mathcal{K}}c_{i_1, v_1}^{(j_1)}c_{i_2, v_2}^{(j_2)}\zeta_{i_1}^{(j_1)}\epsilon_{i_2}^{(j_2)}
K\left(\frac{j_1 - j_2}{\mathcal{B}}\right)\vert\\
\leq C\sum_{i_1 = p_{v_1} - \mathcal{K}}^{p_{v_1} + \mathcal{K}}\sum_{i_2 = p_{v_2} - \mathcal{K}}^{p_{v_2} + \mathcal{K}}K\left(\frac{i_1 - i_2}{\mathcal{B}}\right)\times \mathcal{B}
\sqrt{\sum_{j_1 = q_{v_1} - \mathcal{K}}^{q_{v_1} + \mathcal{K}}c_{i_1, v_1}^{(j_1)2}\zeta_{i_1}^{(j_1)2}}
\times \sqrt{\sum_{j_2 = q_{v_2} - \mathcal{K}}^{q_{v_2} + \mathcal{K}}c_{i_2, v_2}^{(j_2)2}\epsilon_{i_2}^{(j_2)2}}\\
\leq C^\prime \mathcal{B}^2\sqrt{\sum_{i_1 = p_{v_1} - \mathcal{K}}^{p_{v_1} + \mathcal{K}}\sum_{j_1 = q_{v_1} - \mathcal{K}}^{q_{v_1} + \mathcal{K}}c_{i_1, v_1}^{(j_1)2}\zeta_{i_1}^{(j_1)2}}\times \sqrt{
\sum_{i_2 = p_{v_2} - \mathcal{K}}^{p_{v_2} + \mathcal{K}}\sum_{j_2 = q_{v_2} - \mathcal{K}}^{q_{v_2} + \mathcal{K}}c_{i_2, v_2}^{(j_2)2}\epsilon_{i_2}^{(j_2)2}
}
\end{align*}
Therefore, from Cauchy inequality
\begin{align*}
\Vert\sum_{i_1 = p_{v_1} - \mathcal{K}}^{p_{v_1} + \mathcal{K}}\sum_{j_1 = q_{v_1} - \mathcal{K}}^{q_{v_1} + \mathcal{K}}
\sum_{i_2 = p_{v_2} - \mathcal{K}}^{p_{v_2} + \mathcal{K}}\sum_{j_2 = q_{v_2} - \mathcal{K}}^{q_{v_2} + \mathcal{K}}
c_{i_1, v_1}^{(j_1)}c_{i_2, v_2}^{(j_2)}\zeta_{i_1}^{(j_1)}\epsilon_{i_2}^{(j_2)}K\left(\frac{i_1 - i_2}{\mathcal{B}}\right)
K\left(\frac{j_1 - j_2}{\mathcal{B}}\right)\Vert_{M / 2}\\
\leq C\mathcal{B}^2\Vert\sqrt{\sum_{i_1 = p_{v_1} - \mathcal{K}}^{p_{v_1} + \mathcal{K}}\sum_{j_1 = q_{v_1} - \mathcal{K}}^{q_{v_1} + \mathcal{K}}c_{i_1, v_1}^{(j_1)2}\zeta_{i_1}^{(j_1)2}}\times \sqrt{
\sum_{i_2 = p_{v_2} - \mathcal{K}}^{p_{v_2} + \mathcal{K}}\sum_{j_2 = q_{v_2} - \mathcal{K}}^{q_{v_2} + \mathcal{K}}c_{i_2, v_2}^{(j_2)2}\epsilon_{i_2}^{(j_2)2}}\Vert_{M / 2}\\
\leq C\mathcal{B}^2\sqrt{\sum_{i_1 = p_{v_1} - \mathcal{K}}^{p_{v_1} + \mathcal{K}}\sum_{j_1 = q_{v_1} - \mathcal{K}}^{q_{v_1} + \mathcal{K}}c_{i_1, v_1}^{(j_1)2}\Vert\zeta_{i_1}^{(j_1)}\Vert_{M}^2}\times \sqrt{\sum_{i_2 = p_{v_2} - \mathcal{K}}^{p_{v_2} + \mathcal{K}}\sum_{j_2 = q_{v_2} - \mathcal{K}}^{q_{v_2} + \mathcal{K}}c_{i_2, v_2}^{(j_2)2}\Vert\epsilon_{i_2}^{(j_2)}\Vert_M^2}
\end{align*}
Form lemma \ref{lemma.consistent_linear_combination}, $\Vert\zeta_{i}^{(j)}\Vert_{M}\leq \frac{C}{\mathcal{K}}$, so
\begin{equation}
\begin{aligned}
\Vert\sum_{i_1 = p_{v_1} - \mathcal{K}}^{p_{v_1} + \mathcal{K}}\sum_{j_1 = q_{v_1} - \mathcal{K}}^{q_{v_1} + \mathcal{K}}
\sum_{i_2 = p_{v_2} - \mathcal{K}}^{p_{v_2} + \mathcal{K}}\sum_{j_2 = q_{v_2} - \mathcal{K}}^{q_{v_2} + \mathcal{K}}
c_{i_1, v_1}^{(j_1)}c_{i_2, v_2}^{(j_2)}\zeta_{i_1}^{(j_1)}\\
\times\epsilon_{i_2}^{(j_2)}K\left(\frac{i_1 - i_2}{\mathcal{B}}\right)
K\left(\frac{j_1 - j_2}{\mathcal{B}}\right)\Vert_{M / 2}
\leq \frac{C\mathcal{B}^2}{\mathcal{K}}\\
\Rightarrow \max_{v_1,v_2 = 1,\cdots, V}\vert\sum_{i_1 = p_{v_1} - \mathcal{K}}^{p_{v_1} + \mathcal{K}}\sum_{j_1 = q_{v_1} - \mathcal{K}}^{q_{v_1} + \mathcal{K}}
\sum_{i_2 = p_{v_2} - \mathcal{K}}^{p_{v_2} + \mathcal{K}}\sum_{j_2 = q_{v_2} - \mathcal{K}}^{q_{v_2} + \mathcal{K}}
c_{i_1, v_1}^{(j_1)}c_{i_2, v_2}^{(j_2)}\zeta_{i_1}^{(j_1)}\\
\times \epsilon_{i_2}^{(j_2)}K\left(\frac{i_1 - i_2}{\mathcal{B}}\right)
K\left(\frac{j_1 - j_2}{\mathcal{B}}\right)\vert
= O_p\left(V^{4/M}\times \frac{\mathcal{B}^2}{\mathcal{K}}\right)
\end{aligned}
\end{equation}
Similarly
\begin{align*}
\max_{v_1,v_2 = 1,\cdots, V}\vert
\sum_{i_1 = p_{v_1} - \mathcal{K}}^{p_{v_1} + \mathcal{K}}\sum_{j_1 = q_{v_1} - \mathcal{K}}^{q_{v_1} + \mathcal{K}}
\sum_{i_2 = p_{v_2} - \mathcal{K}}^{p_{v_2} + \mathcal{K}}\sum_{j_2 = q_{v_2} - \mathcal{K}}^{q_{v_2} + \mathcal{K}}
c_{i_1, v_1}^{(j_1)}c_{i_2, v_2}^{(j_2)}\epsilon_{i_1}^{(j_1)}\Delta_{i_2}^{(j_2)}K\left(\frac{i_1 - i_2}{\mathcal{B}}\right)
K\left(\frac{j_1 - j_2}{\mathcal{B}}\right)\\
\leq \max_{v_1,v_2 = 1,\cdots, V}\vert
\sum_{i_1 = p_{v_1} - \mathcal{K}}^{p_{v_1} + \mathcal{K}}\sum_{j_1 = q_{v_1} - \mathcal{K}}^{q_{v_1} + \mathcal{K}}
\sum_{i_2 = p_{v_2} - \mathcal{K}}^{p_{v_2} + \mathcal{K}}\sum_{j_2 = q_{v_2} - \mathcal{K}}^{q_{v_2} + \mathcal{K}}
c_{i_1, v_1}^{(j_1)}c_{i_2, v_2}^{(j_2)}\epsilon_{i_1}^{(j_1)}\eta_{i_2}^{(j_2)}K\left(\frac{i_1 - i_2}{\mathcal{B}}\right)
K\left(\frac{j_1 - j_2}{\mathcal{B}}\right)
\vert\\
+ \max_{v_1,v_2 = 1,\cdots, V}\vert
\sum_{i_1 = p_{v_1} - \mathcal{K}}^{p_{v_1} + \mathcal{K}}\sum_{j_1 = q_{v_1} - \mathcal{K}}^{q_{v_1} + \mathcal{K}}
\sum_{i_2 = p_{v_2} - \mathcal{K}}^{p_{v_2} + \mathcal{K}}\sum_{j_2 = q_{v_2} - \mathcal{K}}^{q_{v_2} + \mathcal{K}}
c_{i_1, v_1}^{(j_1)}c_{i_2, v_2}^{(j_2)}\epsilon_{i_1}^{(j_1)}\zeta_{i_2}^{(j_2)}K\left(\frac{i_1 - i_2}{\mathcal{B}}\right)
K\left(\frac{j_1 - j_2}{\mathcal{B}}\right)\vert\\
=O_p\left(\frac{V^{2/M}\times \mathcal{B}\mathcal{K}}{n} + \frac{V^{2/M}\times\mathcal{B}\mathcal{K}}{m} + V^{4/M}\times \frac{\mathcal{B}^2}{\mathcal{K}}\right)
\end{align*}
Finally, notice that
\begin{align*}
\vert
\sum_{i_1 = p_{v_1} - \mathcal{K}}^{p_{v_1} + \mathcal{K}}\sum_{j_1 = q_{v_1} - \mathcal{K}}^{q_{v_1} + \mathcal{K}}
\sum_{i_2 = p_{v_2} - \mathcal{K}}^{p_{v_2} + \mathcal{K}}\sum_{j_2 = q_{v_2} - \mathcal{K}}^{q_{v_2} + \mathcal{K}}
c_{i_1, v_1}^{(j_1)}c_{i_2, v_2}^{(j_2)}\Delta_{i_1}^{(j_1)}\Delta_{i_2}^{(j_2)}K\left(\frac{i_1 - i_2}{\mathcal{B}}\right)
K\left(\frac{j_1 - j_2}{\mathcal{B}}\right)
\vert\\
\leq C\sum_{i_1 = p_{v_1} - \mathcal{K}}^{p_{v_1} + \mathcal{K}}\sum_{i_2 = p_{v_2} - \mathcal{K}}^{p_{v_2} + \mathcal{K}}
K\left(\frac{i_1 - i_2}{\mathcal{B}}\right)\mathcal{B}\sqrt{\sum_{j_1 = q_{v_1} - \mathcal{K}}^{q_{v_1} + \mathcal{K}}c_{i_1, v_1}^{(j_1)2}\Delta_{i_1}^{(j_1)2}}\sqrt{\sum_{j_2 = q_{v_2} - \mathcal{K}}^{q_{v_2} + \mathcal{K}}c_{i_2, v_2}^{(j_2)2}\Delta_{i_2}^{(j_2)2}}\\
\leq C^\prime \mathcal{B}^2\sqrt{\sum_{i_1 = p_{v_1} - \mathcal{K}}^{p_{v_1} + \mathcal{K}}\sum_{j_1 = q_{v_1} - \mathcal{K}}^{q_{v_1} + \mathcal{K}}c_{i_1, v_1}^{(j_1)2}\Delta_{i_1}^{(j_1)2}}\times \sqrt{\sum_{i_2 = p_{v_2} - \mathcal{K}}^{p_{v_2} + \mathcal{K}}\sum_{j_2 = q_{v_2} - \mathcal{K}}^{q_{v_2} + \mathcal{K}}c_{i_2, v_2}^{(j_2)2}\Delta_{i_2}^{(j_2)2}}\\
\Rightarrow \Vert \sum_{i_1 = p_{v_1} - \mathcal{K}}^{p_{v_1} + \mathcal{K}}\sum_{j_1 = q_{v_1} - \mathcal{K}}^{q_{v_1} + \mathcal{K}}
\sum_{i_2 = p_{v_2} - \mathcal{K}}^{p_{v_2} + \mathcal{K}}\sum_{j_2 = q_{v_2} - \mathcal{K}}^{q_{v_2} + \mathcal{K}}
c_{i_1, v_1}^{(j_1)}c_{i_2, v_2}^{(j_2)}\Delta_{i_1}^{(j_1)}\Delta_{i_2}^{(j_2)}K\left(\frac{i_1 - i_2}{\mathcal{B}}\right)
K\left(\frac{j_1 - j_2}{\mathcal{B}}\right)\Vert_{M / 2}\\
\leq C\mathcal{B}^2\sqrt{\sum_{i_1 = p_{v_1} - \mathcal{K}}^{p_{v_1} + \mathcal{K}}\sum_{j_1 = q_{v_1} - \mathcal{K}}^{q_{v_1} + \mathcal{K}}c_{i_1, v_1}^{(j_1)2}\Vert\Delta_{i_1}^{(j_1)}\Vert_M^2}\times \sqrt{\sum_{i_2 = p_{v_2} - \mathcal{K}}^{p_{v_2} + \mathcal{K}}\sum_{j_2 = q_{v_2} - \mathcal{K}}^{q_{v_2} + \mathcal{K}}c_{i_2, v_2}^{(j_2)2}\Vert\Delta_{i_2}^{(j_2)}\Vert_M^2}
\end{align*}
and $\Vert\Delta_{i}^{(j)}\Vert_M\leq \vert\eta_i^{(j)}\vert + \Vert\zeta_i^{(j)}\Vert_M\leq \frac{C\mathcal{K}}{n} + \frac{C\mathcal{K}}{m} + \frac{C}{\mathcal{K}}$, so
\begin{equation}
\begin{aligned}
\max_{v_1, v_2 = 1,\cdots, V}\vert
\sum_{i_1 = p_{v_1} - \mathcal{K}}^{p_{v_1} + \mathcal{K}}\sum_{j_1 = q_{v_1} - \mathcal{K}}^{q_{v_1} + \mathcal{K}}
\sum_{i_2 = p_{v_2} - \mathcal{K}}^{p_{v_2} + \mathcal{K}}\sum_{j_2 = q_{v_2} - \mathcal{K}}^{q_{v_2} + \mathcal{K}}
c_{i_1, v_1}^{(j_1)}c_{i_2, v_2}^{(j_2)}\\
\times\Delta_{i_1}^{(j_1)}\Delta_{i_2}^{(j_2)}K\left(\frac{i_1 - i_2}{\mathcal{B}}\right)
K\left(\frac{j_1 - j_2}{\mathcal{B}}\right)
\vert = O_p\left(\frac{V^{4/M}\mathcal{B}^2\mathcal{K}^2}{n^2} + \frac{V^{4/M}\mathcal{B}^2\mathcal{K}^2}{m^2} + \frac{V^{4/M}\mathcal{B}^2}{\mathcal{K}^2}\right)
\end{aligned}
\end{equation}
In particular,
\begin{equation}
\begin{aligned}
\max_{v_1,v_2 = 1,\cdots, V}\vert\sum_{i_1 = p_{v_1} - \mathcal{K}}^{p_{v_1} + \mathcal{K}}\sum_{j_1 = q_{v_1} - \mathcal{K}}^{q_{v_1} + \mathcal{K}}
\sum_{i_2 = p_{v_2} - \mathcal{K}}^{p_{v_2} + \mathcal{K}}\sum_{j_2 = q_{v_2} - \mathcal{K}}^{q_{v_2} + \mathcal{K}}
c_{i_1, v_1}^{(j_1)}c_{i_2, v_2}^{(j_2)}\widehat{\epsilon}_{i_1}^{(j_1)}\widehat{\epsilon}_{i_2}^{(j_2)}K\left(\frac{i_1 - i_2}{\mathcal{B}}\right)
K\left(\frac{j_1 - j_2}{\mathcal{B}}\right)\\
-\sum_{i_1 = 1}^n\sum_{j_1 = 1}^m\sum_{i_2 = 1}^n\sum_{j_2 = 1}^m c_{i_1, v_1}^{(j_1)}c_{i_2, v_2}^{(j_2)}
\mathbf{E}\epsilon_{i_1}^{(j_1)}\epsilon_{i_2}^{(j_2)}
\vert\\
= O_p\left(\mathcal{S}(\mathcal{B}) + \frac{V^{4/M}\mathcal{B}^3}{\mathcal{K}} + \frac{V^{2/M}\mathcal{BK}}{n}
+\frac{V^{2/M}\mathcal{BK}}{m}\right)
\end{aligned}
\label{eq.delta_hat}
\end{equation}
Define $\delta = \mathcal{S}(\mathcal{B}) + \frac{V^{4/M}\mathcal{B}^3}{\mathcal{K}} + \frac{V^{2/M}\mathcal{BK}}{n}
+\frac{V^{2/M}\mathcal{BK}}{m}$. Notice that $\delta^{1/3}\log(\mathcal{K})\to 0$  as $\min(n,m)\to\infty$,
from lemma \ref{lemma.Gaussian_property}, we prove eq.\eqref{eq.bootstrap_consistency}.

According to Algorithm \ref{algorithm.bootstrap}, define
\begin{align*}
    \widehat{T}_v^{\dagger *} = \frac{T_{n,m}}{B_{n,m}\times\widehat{\tau}_v} (\widehat{\mu}^*(x_v,y_v) - \widehat{\mu}(x_v, y_v))
    = \frac{1}{\widehat{\tau}_v}\sum_{i = 1}^n\sum_{j = 1}^m c_{i,v}^{(j)}\epsilon_i^{(j)*},
\end{align*}
then ${T_1^{\dagger*}}, \cdots, {T_V^{\dagger*}}$ are joint normal random variables with mean $0$ and covariance
$$
\mathbf{E}^*T_{v_1}^{\dagger*}T_{v_2}^{\dagger*} = \frac{1}{\widehat{\tau}_{v_1}\widehat{\tau}_{v_2}}
\sum_{i_1 = p_{v_1} - \mathcal{K}}^{p_{v_1} + \mathcal{K}}\sum_{j_1 = q_{v_1} - \mathcal{K}}^{q_{v_1} + \mathcal{K}}\sum_{i_2 = p_{v_2} - \mathcal{K}}^{p_{v_2} + \mathcal{K}}\sum_{j_2 = q_{v_2} - \mathcal{K}}^{q_{v_2} + \mathcal{K}}
c_{i_1, v_1}^{(j_1)}c_{i_2, v_2}^{(j_2)}\widehat{\epsilon}_{i_1}^{(j_1)}\widehat{\epsilon}_{i_2}^{(j_2)}K\left(\frac{i_1 - i_2}{\mathcal{B}}\right)
K\left(\frac{j_1 - j_2}{\mathcal{B}}\right)
$$
in the bootstrap world (i.e., conditional on $X_i^{(j)}, i = 1,\cdots, n, j = 1,\cdots, m$). Since
\begin{align*}
\vert \mathbf{E}^*T_{v_1}^{\dagger*}T_{v_2}^{\dagger*} - \mathbf{E}\xi_{v_1}^\dagger\xi_{v_2}^\dagger\vert\\
\leq \frac{1}{\widehat{\tau}_{v_1}\widehat{\tau}_{v_2}}\vert\sum_{i_1 = p_{v_1} - \mathcal{K}}^{p_{v_1} + \mathcal{K}}\sum_{j_1 = q_{v_1} - \mathcal{K}}^{q_{v_1} + \mathcal{K}}
\sum_{i_2 = p_{v_2} - \mathcal{K}}^{p_{v_2} + \mathcal{K}}\sum_{j_2 = q_{v_2} - \mathcal{K}}^{q_{v_2} + \mathcal{K}}
c_{i_1, v_1}^{(j_1)}c_{i_2, v_2}^{(j_2)}\widehat{\epsilon}_{i_1}^{(j_1)}\widehat{\epsilon}_{i_2}^{(j_2)}K\left(\frac{i_1 - i_2}{\mathcal{B}}\right)
K\left(\frac{j_1 - j_2}{\mathcal{B}}\right)\\
-\sum_{i_1 = 1}^n\sum_{j_1 = 1}^m\sum_{i_2 = 1}^n\sum_{j_2 = 1}^m c_{i_1, v_1}^{(j_1)}c_{i_2, v_2}^{(j_2)}
\mathbf{E}\epsilon_{i_1}^{(j_1)}\epsilon_{i_2}^{(j_2)}\vert\\
+ \vert \sum_{i_1 = 1}^n\sum_{j_1 = 1}^m\sum_{i_2 = 1}^n\sum_{j_2 = 1}^m c_{i_1, v_1}^{(j_1)}c_{i_2, v_2}^{(j_2)}
\mathbf{E}\epsilon_{i_1}^{(j_1)}\epsilon_{i_2}^{(j_2)}\vert\times \frac{\vert \widehat{\tau}_{v_1}\widehat{\tau}_{v_2} - \tau_{v_1}\tau_{v_2} \vert}{\widehat{\tau}_{v_1}\widehat{\tau}_{v_2}\tau_{v_1}\tau_{v_2}}
\end{align*}
With probability tending to $1$, we have
\begin{equation}
\begin{aligned}
\vert \mathbf{E}^*T_{v_1}^{\dagger*}T_{v_2}^{\dagger*} - \mathbf{E}\xi_{v_1}^\dagger\xi_{v_2}^\dagger\vert\\
\leq c^\prime\vert\sum_{i_1 = p_{v_1} - \mathcal{K}}^{p_{v_1} + \mathcal{K}}\sum_{j_1 = q_{v_1} - \mathcal{K}}^{q_{v_1} + \mathcal{K}}
\sum_{i_2 = p_{v_2} - \mathcal{K}}^{p_{v_2} + \mathcal{K}}\sum_{j_2 = q_{v_2} - \mathcal{K}}^{q_{v_2} + \mathcal{K}}
c_{i_1, v_1}^{(j_1)}c_{i_2, v_2}^{(j_2)}\widehat{\epsilon}_{i_1}^{(j_1)}\widehat{\epsilon}_{i_2}^{(j_2)}K\left(\frac{i_1 - i_2}{\mathcal{B}}\right)
K\left(\frac{j_1 - j_2}{\mathcal{B}}\right)\\
-\sum_{i_1 = 1}^n\sum_{j_1 = 1}^m\sum_{i_2 = 1}^n\sum_{j_2 = 1}^m c_{i_1, v_1}^{(j_1)}c_{i_2, v_2}^{(j_2)}
\mathbf{E}\epsilon_{i_1}^{(j_1)}\epsilon_{i_2}^{(j_2)}\vert
+ C\vert\widehat{\tau}_{v_1} - \tau_{v_1}\vert + C\vert\widehat{\tau}_{v_2} - \tau_{v_2}\vert\\
\Rightarrow \max_{v_1,v_2 = 1,\cdots, V}\vert \mathbf{E}^*T_{v_1}^{\dagger*}T_{v_2}^{\dagger*} - \mathbf{E}\xi_{v_1}^\dagger\xi_{v_2}^\dagger\vert
= O_p\left(\mathcal{S}(\mathcal{B}) + \frac{V^{4/M}\mathcal{B}^3}{\mathcal{K}} + \frac{V^{2/M}\mathcal{BK}}{n}
+\frac{V^{2/M}\mathcal{BK}}{m}\right)
\end{aligned}
\end{equation}
and we prove \eqref{eq.bootstrap_consistency_hetero}.

\end{proof}

\scriptsize
\bibliographystyle{Chicago}
\bibliography{combinepdf}
\end{document}